\documentclass[graybox]{svmult}       

\smartqed  

\usepackage{bm}

\RequirePackage[OT1]{fontenc}
\RequirePackage{amsmath}

\RequirePackage[colorlinks,citecolor=blue,urlcolor=blue]{hyperref}

\usepackage{type1cm}        
\usepackage{makeidx}         
\usepackage{graphicx}        
\usepackage{multicol}        
\usepackage[bottom]{footmisc}

\usepackage{newtxtext}       %
\usepackage{newtxmath}       

\makeindex             

\usepackage{graphicx}
\usepackage{amsmath}
\usepackage{graphics}
\usepackage{epsfig}
\usepackage{amssymb}
\usepackage{color}
\usepackage{float}
\usepackage[section]{placeins}
\usepackage{subcaption}
\usepackage{tikz}

\usepackage[margin=1.2in]{geometry}

\newcommand{\be}{\begin{eqnarray}}
\newcommand{\ee}{\end{eqnarray}}
\newcommand{\bea}{\begin{eqnarray*}}
\newcommand{\eea}{\end{eqnarray*}}


\begin{document}

\title*{Non-lattice  covering and quanitization of high dimensional sets}

\author{ Jack Noonan and Anatoly Zhigljavsky
}

\institute{J. Noonan \at
               School of Mathematics, Cardiff University, Cardiff, CF244AG, UK \\
              \email{Noonanj1@cardiff.ac.uk}
           \and
       A. Zhigljavsky \at
              School of Mathematics, Cardiff University, Cardiff, CF244AG, UK \\
              \email{ZhigljavskyAA@cardiff.ac.uk}           
}


\maketitle

\abstract{The main problem considered in this paper is construction and theoretical study of efficient $n$-point coverings of a $d$-dimensional cube $[-1,1]^d$. Targeted values of $d$ are between 5 and 50; $n$ can be in hundreds or thousands and the designs (collections of points) are nested.
This paper is a continuation of our paper \cite{us}, where we have theoretically investigated several simple schemes and numerically studied many more.
In this paper, we extend the theoretical constructions  of  \cite{us} for studying the designs which were found to be  superior to the ones theoretically investigated  in \cite{us}. We  also extend our constructions for new construction schemes which provide even better coverings (in the class of nested designs) than the ones numerically found in \cite{us}. In view of a close connection of the problem of quantization to the problem of covering, we extend our theoretical approximations and practical recommendations to the problem of construction of efficient quantization designs in a cube $[-1,1]^d$. In the last section, we discuss the problems of covering and quantization in a $d$-dimensional simplex; practical significance of this problem has been communicated to the authors by Professor Michael  Vrahatis, a co-editor of the present volume.
}

\section{Introduction}
\label{sec:introd}

The problem of the main importance  in this paper is the following problem of covering a cube $[-1,1]^d$ by $n$ balls.
Let  $Z_1, \ldots, Z_n$ be a collection of points in $\mathbb{R}^d$ and ${\cal B}_d(Z_j,r)= \{Z: \|Z-Z_j\|\leq r\}$ be the Euclidean balls of radius $r$ centered at $Z_j$ $(j=1, \ldots, n)$.
The dimension $d$, the number of balls $n$ and their radius $r$  could be arbitrary.

We are interested in choosing the locations of the centers of the balls $Z_1, \ldots, Z_n$ so that the union of the balls $\cup_{j}{\cal B}_d(Z_j,r)$ covers the largest possible proportion of the cube $[-1,1]^d$. More precisely, we are interested in choosing a  collection of points (called `design') $\mathbb{Z}_n=\{Z_1, \ldots, Z_n\}$  so that
\be
\label{eq:cover1}
\mbox{$C_d(\mathbb{Z}_n,r):= $vol$([-1,1]^d \cap {\cal B}_d(\mathbb{Z}_n,r))/2^d$}
\ee
is as large as possible (given $n$, $r$ and the freedom we are able to use in choosing $Z_1, \ldots, Z_n$). Here ${\cal B}_d(\mathbb{Z}_n,r)$ is the union of the balls
\be
\label{eq:cover1a}
{\cal B}_d(\mathbb{Z}_n,r)= \bigcup_{j=1}^n {\cal B}_d(Z_j,r)
\ee
and $C_d(\mathbb{Z}_n,r)$ is the proportion of the cube $[-1,1]^d$ covered by ${\cal B}_d(\mathbb{Z}_n,r)$. If $Z_j \in \mathbb{Z}_n$ are random
then we shall consider $\mathbb{E}_{\mathbb{Z}_n}C_d(\mathbb{Z}_n,r)$, the expected value of the proportion \eqref{eq:cover1}; for simplicity of notation, we will drop $\mathbb{E}_{\mathbb{Z}_n} $ while referring to  $\mathbb{E}_{\mathbb{Z}_n}C_d(\mathbb{Z}_n,r)$.

For a  design $\mathbb{Z}_n$, its covering radius is defined by CR$(\mathbb{Z}_n)= \max_{X \in {\cal  C}_{d}} \min_{Z_j \in \mathbb{Z}_n}\|X-Z_j\|$. In computer experiments, covering radius is called minimax-distance criterion, see \cite{johnson1990minimax} and \cite{pronzato2012design}; in the theory of low-discrepancy sequences, covering radius is called dispersion, see \cite[Ch. 6]{niederreiter1992random}.

The problem of optimal covering of a cube by $n$ balls has very high importance  for the theory of global optimization and many branches of numerical mathematics.
In particular, the  $n$-point designs $\mathbb{Z}_n$ with smallest CR provide the following:
  (a) the $n$-point min-max optimal quadratures, see \cite[Ch.3,Th.1.1]{sukharev2012minimax},
(b)
 min-max $n$-point global optimization methods in the set of all adaptive $n$-point optimization strategies, see \cite[Ch.4,Th.2.1]{sukharev2012minimax}, and
 (c) worst-case  $n$-point multi-objective global optimization methods in the set of all adaptive $n$-point algorithms, see
 \cite{vzilinskas2013worst}.
  In all three cases, the class of (objective) functions is the class of Liptshitz functions, where the  Liptshitz constant may be unknown.
The results (a) and (b) are the celebrated results of A.G.Sukharev obtained in the late nineteen-sixties, see e.g. \cite{sukharev1971optimal}, and (c) is a recent result of A. \v{Z}ilinskas.

If $d$ is not small (say, $d>5$) then computation of the covering radius CR$(\mathbb{Z}_n) $ for any non-trivial design~$\mathbb{Z}_n$ is a very difficult computational problem. This explains why the problem of construction of optimal $n$-point designs with smallest covering radius is notoriously difficult, see for example recent surveys \cite{toth20172,toth1993packing}.
If $r=$CR$(\mathbb{Z}_n)$, then
 $C_d(\mathbb{Z}_n,r)$ defined in \eqref{eq:cover1} is equal to 1, and  the whole cube ${\cal C}_{d}$ gets covered by the balls.
 However,  we are only interested in reaching the values like 0.95 or 0.99, when only a large part of the ball is covered.

We will say that ${\cal B}_d(\mathbb{Z}_n,r)$ makes a $(1-\gamma)$-covering of $[-1,1]^d$ if
\be
\label{eq:quant}
C_d(\mathbb{Z}_n,r)=1-\gamma \,;
\ee
the corresponding value of $r$ will be called $(1\!-\!\gamma)$-covering radius and denoted $r_{1-\gamma}$ or $r_{1-\gamma}(\mathbb{Z}_n)$.
  If $\gamma=0$ then the $(1\!-\!\gamma)$-covering becomes the full covering {and 1-covering radius  $r_{1}(\mathbb{Z}_n)$ becomes the covering radius CR$(\mathbb{Z}_n)$. The problem of construction of efficient designs with smallest possible $(1\!-\!\gamma)$-covering radius (with some small $\gamma>0$) will be referred to as the problem of weak covering.\\

  Let us give two strong arguments why  the problem of weak covering could be even more practically important than the problem of
  full covering.

  \begin{itemize}
    \item  Numerical checking of weak covering (with an approximate value of $\gamma$) is straightforward while
    numerical checking of the full covering is practically impossible, if $d$ is large enough.

    \item   For a given design $\mathbb{Z}_n$,  $C_d(\mathbb{Z}_n,r)$ defined in \eqref{eq:cover1} and considered as a function of $r$, is a
 cumulative distribution function (c.d.f.) of the random variable (r.v.)   $ \varrho(U,\mathbb{Z}_n)= \min_{Z_i \in \mathbb{Z}_n} \|U-Z_i\| $,
 where $U$ is a random vector uniformly distributed on $[-1,1]^d$,  see  \eqref{eq:rho_s} below. The covering radius CR$(\mathbb{Z}_n)$ is the upper bound of this r.v. while in view of \eqref{eq:quant}, $r_{1-\gamma}(\mathbb{Z}_n)$ is the $(1-\gamma)$-quantile.
 Many practically important characteristics of designs such as quantization error considered in Section~\ref{sec:quantization} are expressed in terms of the whole c.d.f.
 $C_d(\mathbb{Z}_n,r)$ and their dependence on the upper bound CR$(\mathbb{Z}_n)$ is marginal.
 As shown in Section~\ref{sec:q_eff}, numerical studies indicate that comparison of designs on the base of their weak coverage properties is very similar to quantization error comparisons, but this  may not be true for comparisons with respect to CR$(\mathbb{Z}_n)$. This phenomenon is similar to the  well-known fact in the
 theory of space covering by lattices (see  an excellent book  \cite{Conway} and surveys \cite{toth20172,toth1993packing}), where
 best lattice coverings of space are often poor quantizers and vice-versa. Moreover, Figures~\ref{Main_pic1}-\ref{Main_pic2} below show that CR$(\mathbb{Z}_n)$  may give a totally inadequate impression about the c.d.f. $C_d(\mathbb{Z}_n,r)$ and could be much larger than  $r_{1-\gamma}(\mathbb{Z}_n)$ with very small $\gamma>0$.

  \end{itemize}

 In Figures~\ref{Main_pic1}--\ref{Main_pic2} we consider two simple designs for which we plot their c.d.f. $C_d(\cdot,r)$, black line, and also   indicate the location of the $r_{1}$=CR and $r_{0.999}$ by vertical red and green line respectively.
   In Figure~\ref{Main_pic1}, we  take $d=10$, $n=512$ and use a $2^{d-1} $ design of maximum resolution concentrated at the
  points\footnote{For simplicity of notation, vectors in $\mathbb{R}^d$ are represented as rows.}
  $(\pm 1/2, \ldots, \pm 1/2) \in \mathbb{R}^d$ as design $\mathbb{Z}_n$; this design  is a particular case of Design 4 of Section~\ref{sec:numeric2} and can be defined for any $d>2$.
  In Figure~\ref{Main_pic2}, we keep $d=10$ but take the full factorial $2^{d} $ design with $m=2^d$ points, again concentrated at the points $(\pm 1/2, \ldots, \pm 1/2)$; denote this design $\mathbb{Z}_m^\prime$.

  For both designs, it is very easy to analytically compute their covering radii (for any $d>2$):
  CR$(\mathbb{Z}_n)=\sqrt{d+8}/2 $ and CR$(\mathbb{Z}_m^\prime)=\sqrt{d}/2  $; for $d=10$ this gives
  CR$(\mathbb{Z}_n) \simeq 2.1213 $ and CR$(\mathbb{Z}_m^\prime) \simeq 1.58114. $
The values of  $r_{0.999}$ are: $r_{0.999}(\mathbb{Z}_n) \simeq 1.3465 $ and $r_{0.999}(\mathbb{Z}_m^\prime) \simeq 1.2708$.
Their values have been computed using very accurate approximations developed in \cite{us1}; we claim   3 correct decimal places in both values of $r_{0.999}$. We will return to this example in Section~\ref{sec:com_2}.

\begin{figure}[!h]
\centering
\begin{minipage}{.5\textwidth}
  \centering
  \includegraphics[width=1\linewidth]{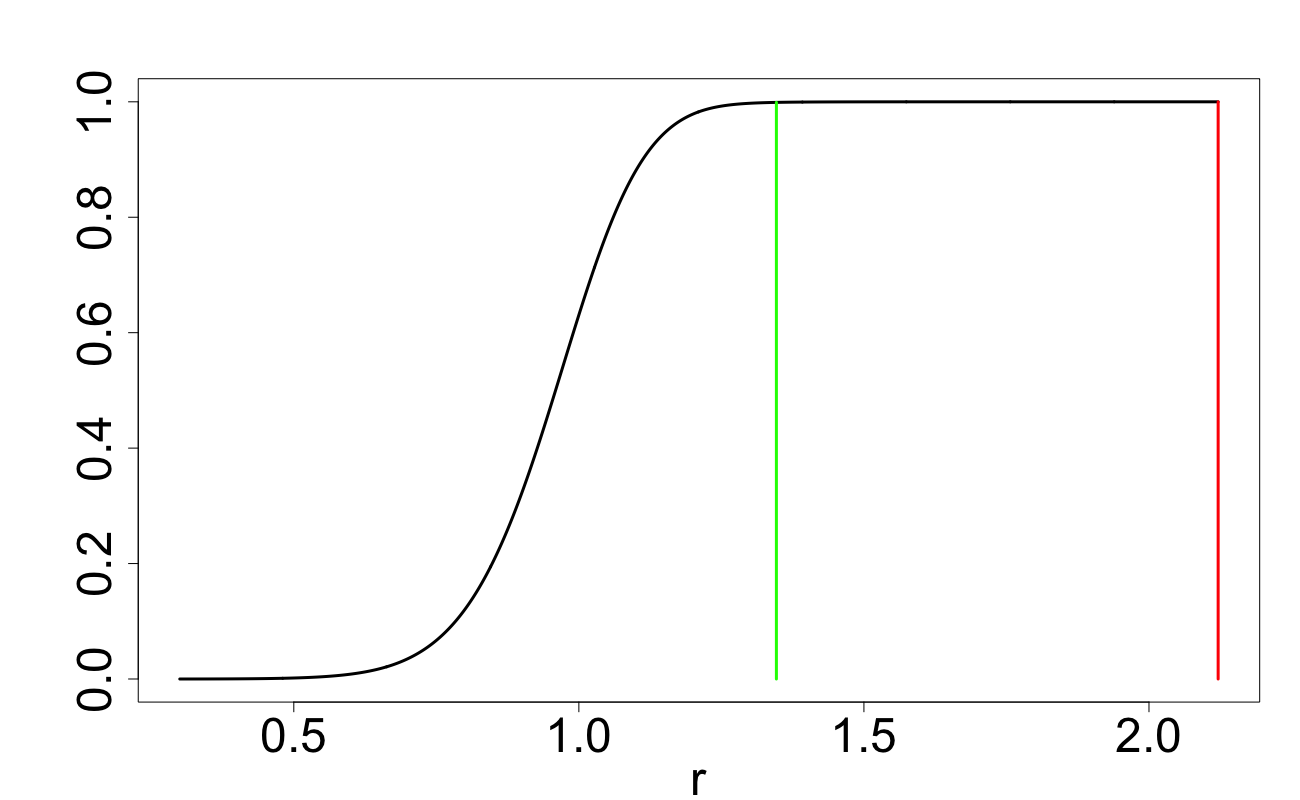}
  \caption{$C_d({\mathbb{Z}_{n}},r_{})$ with $r_{0.999}$ and $r_{1}$: $d=10$, \\
 $\mathbb{Z}_{n}$ is a $2^{d-1}$-factorial design with $n=2^{d-1}$}
  \label{Main_pic1}
\end{minipage}%
\begin{minipage}{.5\textwidth}
  \centering
  \includegraphics[width=1\linewidth]{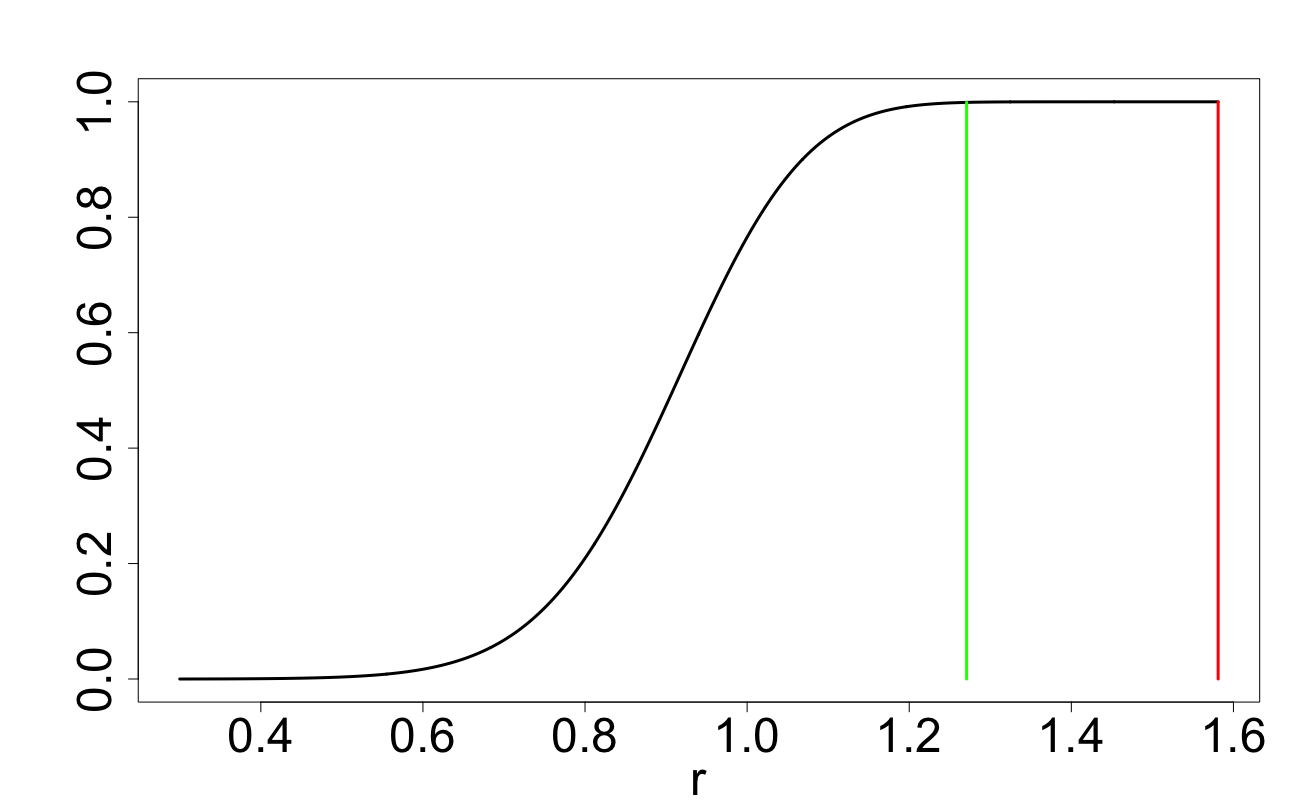}
 \caption{$C_d({\mathbb{Z}_{m}^\prime},r_{})$ with $r_{0.999}$ and $r_{1}$: $d=10$, \\
 $\mathbb{Z}_{m}$ is a $2^{d}$-factorial design}  \label{Main_pic2}
\end{minipage}
\end{figure}

Of course, for any $\mathbb{Z}_n=\{Z_1, \ldots, Z_n\}$ we can reach $C_d(\mathbb{Z}_n,r)= 1$ by means of increase of $r$. Likewise, for any given~$r$
we can
reach $C_d(\mathbb{Z}_n,r)= 1$ by sending $n \to \infty$. However, we are not interested in very large values of $n$ and try to get the coverage
of the most part of the cube ${\cal C}_{d}$ with the radius  $r$ as small as possible.  We will keep in mind the following typical values of $d$ and $n$ which we will use for illustrating our results:  $d=5,10,20,50$; $n=2^k $ with $k=6, \ldots, 11$ (we have chosen $n$ as a power of~2 since this a favorable number for Sobol's sequence (Design 3) as well as  Design 4 defined in Section~\ref{sec:numeric2}).\\

The structure of the rest of the paper is as follows.
In Section~\ref{sec:2weak} we
discuss the concept of weak covering in more detail and introduce three generic designs which we will concentrate our attention on.
In Sections~\ref{sec:des1}, \ref{sec:des2a} and \ref{sec:des2b}
we derive approximations for  the expected volume of intersection the cube $[-1,1]^d$ with
  $n$  balls centred at the points of these designs. In Section~\ref{sec:numeric}, we  provide numerical results showing that  the developed approximations
are very accurate.
In Section~\ref{sec:quantization}, we derive approximations for the mean squared quantization error for chosen families of designs and
numerically demonstrate that the developed approximations are very accurate. In Section~\ref{sec:numeric2}, we numerically compare
covering and quantization properties of different designs including scaled Sobol's
sequence and a family of very efficient designs defined only for very specific values of $n$.
In Section~\ref{sec:simplex} we try to answer the question raised by Michael  Vrahatis by
 numerically investigating the importance of the  effect of scaling points away from the boundary (we call it $\delta$-effect)  for covering and quantization in a $d$-dimensional simplex.
 In Appendix, Section~\ref{sec:appB}, we formulate a  simple but important lemma about the distribution and moments of a certain random variable.\\

Our main theoretical contributions in this paper are:

\begin{itemize}
\item derivation of accurate approximations \eqref{2eq:inters2f_corrected}
and \eqref{alpha_zero_corrected}
for the probability
$P_{U,\delta,\alpha,r}$ defined in \eqref{eq:prob_u};

  \item  derivation of accurate approximations \eqref{eq:accurate_app}, \eqref{eq:accurate_app2} and  \eqref{eq:prod7a} for the expected volume of intersection of the cube $[-1,1]^d$ with
  $n$  balls centred at the points of the selected designs;
   \item derivation of accurate approximations
\eqref{final_quant}, \eqref{final_quant_alpha_zero} and \eqref{final_quant_alpha_zero2b} for the mean squared quantization error for the selected designs.
\end{itemize}

We have performed a large-scale numerical study and provided a number of figures and tables. The following are the key  messages containing in these figures and tables.

\begin{itemize}
  \item Figures~\ref{Main_pic1}--\ref{Main_pic2}: weak covering could be much more practically useful than the full covering;
  \item Figures~\ref{Overall1}--\ref{One_ball_approx1_alph_zero2}: developed approximations for  the probability
$P_{U,\delta,\alpha,r}$ defined in \eqref{eq:prob_u} are very accurate;
\item  Figures~\ref{pos_alpha1}--\ref{end_pic}: (a) developed approximations for $C_d(\mathbb{Z}_n,r) $  are very accurate,  (b) there is a very strong $\delta$-effect for all three types of designs, and (c) this $\delta$-effect gets stronger as $d$ increases;
    \item Tables~\ref{alpha_table} and \ref{alpha_table2} and Figures~\ref{alpha_comparison1}-\ref{alpha_comparison2}: smaller values of $\alpha$ are beneficial in Design 1 but  Design 2 (where $\alpha=0$) becomes inefficient when $n$ gets close to $2^d$;
        \item Figures~\ref{worsen_figure111}--\ref{quant_without2}: developed approximations for the quantization error are very accurate
        and there is a very strong $\delta$-effect for all three types of designs used for quantization;
        \item Tables~\ref{table_d_10}--\ref{Table_d_20} and Figures~\ref{function_of_r1}-\ref{function_of_r2}: (a) Designs 2a and especially 2b provide very high quality coverage for suitable $n$,
        (b) properly $\delta$-tuned deterministic non-nested Design 4 provides superior  covering, (c)
         coverage properties of $\delta$-tuned  low-discrepancy sequences are much better than of the original low-discrepancy sequences, and
         (d) coverage properties of unadjusted low-discrepancy sequences is very low, if dimension $d$ is not small;
         \item Tables~\ref{d_10_quantization} and \ref{d_20_quantization}, Figures~\ref{dom1} and~\ref{dom2}: very similar conclusions to the above but made with respect to the quantization error;

         \item Figures~\ref{simplex_1_pic}--\ref{quant_s2_end}: the $\delta$-effect for covering and quantization schemes in a simplex is definitely present (this effect is more apparent in quantization) but it is much weaker than in a cube.

\end{itemize}

\section{Weak covering}

\label{sec:2weak}

In this section, we consider the problem of weak covering defined and discussed in Section~\ref{sec:introd}.
The main
 characteristic of interest will be $C_d(\mathbb{Z}_n,r)$,
the proportion of the cube   covered by the union of balls ${\cal B}_d(\mathbb{Z}_n,r)$; it is defined in \eqref{eq:cover1}.
We start the section with  short discussion on comparison of designs based on their covering properties. 

\subsection{Comparison of designs from the view-point of weak covering}
\label{sec:com_2}

   Two different designs will be differentiated in terms of covering performance as follows. Fix $d$ and let $\mathbb{Z}_n$ and $\mathbb{Z}_n^{\prime}$ be two $n$-point designs. For $(1-\gamma)$-covering with $\gamma\geq 0$, if $C_d(\mathbb{Z}_n,r)=C_d(\mathbb{Z}_n^{\prime},r^\prime)=1-\gamma$ and $r<r^\prime$, then the design $\mathbb{Z}_n$ provides a more efficient $(1-\gamma)$-covering and is therefore preferable. Moreover, the natural scaling for the radius is $r_n = n ^{1/d} r$ and therefore we can compare an $n$-point design $\mathbb{Z}_n$ with an $m$-point design $\mathbb{Z}_m^{\prime}$ as follows:
if $C_d(\mathbb{Z}_n,r)=C_d(\mathbb{Z}_m^{\prime},r^\prime)=1-\gamma$ and $n ^{1/d} r<m ^{1/d} r^\prime$,  then we say that  the design $\mathbb{Z}_n$ provides a more efficient $(1-\gamma)$-covering than the design $\mathbb{Z}_m^{\prime}$.

As an example, consider the designs used  for plotting Figures~\ref{Main_pic1} and \ref{Main_pic2} in Section~\ref{sec:introd}:
$\mathbb{Z}_n$ with $n=2^{d-1}$ and  $\mathbb{Z}_m^\prime$ with $m=2^d$.  For the full covering, we have for any $d$:
$$
n ^{1/d} r_1(\mathbb{Z}_n)= {2^{-1/d} \sqrt{d+8}}  > { \sqrt{d}}= r_1(\mathbb{Z}_m^\prime)    m ^{1/d}
$$
so that the design  $\mathbb{Z}_m^\prime$  is better than $\mathbb{Z}_n$  for the full covering for any $d$ and the difference between normalized covering radii is quite  significant. For example, for $d=10$ we have
$$
n ^{1/d} r_1(\mathbb{Z}_n) \simeq 3.9585 \;\;\; {\rm and} \;\;  r_1(\mathbb{Z}_m^\prime)    m ^{1/d} \simeq 3.1623
$$

For 0.999-covering, however, the situation is reverse, at least for $d=10$, where we have:
$$
n ^{1/d} r_{0.999}(\mathbb{Z}_n) \simeq   2.5126  < 2.5416  \simeq  r_1(\mathbb{Z}_m^\prime)    m ^{1/d}
$$
and therefore the design $\mathbb{Z}_n$ is better for 0.999-covering than the design $\mathbb{Z}_m^\prime$ for $d=10$.

\subsection{Reduction to the probability of  covering a point by one ball}

In the designs $\mathbb{Z}_n$, which are of most interest to us, the points $Z_j \in \mathbb{Z}_n$ are i.i.d. random vectors in $ \mathbb{R}^d$
with a specified  distribution. Let us show that for these designs,  we can reduce computation of $C_d(\mathbb{Z}_n,r)$ to the probability of covering $[-1,1]^d$ by one ball.

Let $Z_1, \ldots, Z_n$ be i.i.d. random vectors in $ \mathbb{R}^d$ and $
{\cal B}_d(\mathbb{Z}_n,r)$ be as defined in \eqref{eq:cover1a}.
Then, for given $U = (u_1, \ldots, u_d) \in \mathbb{R}^d$,
\be
\mathbb{P} \left\{ U \in {\cal B}_d(\mathbb{Z}_n,r)  \right\}&=& 1-\prod_{j=1}^n \mathbb{P} \left\{ U \notin {\cal B}_d({Z}_j,r)  \right\} \nonumber \\
&=& 1-\prod_{j=1}^n\left(1-\mathbb{P} \left\{ U \in {\cal B}_d({Z}_j,r)  \right\} \right) \nonumber\\
&=& 1-\bigg(1-\mathbb{P}_{_Z} \left\{ \|U - {Z} \| \leq r \right\} \bigg)^n\, .
\label{eq:prod}
\ee

$C_d(\mathbb{Z}_n,r)$, defined in \eqref{eq:cover1},
 is simply
\be \label{eq:prod5}
C_d(\mathbb{Z}_n,r) = \mathbb{E}_{_U} \mathbb{P} \left\{ U \in {\cal B}_d(\mathbb{Z}_n,r)  \right\}\, ,
\ee
where the expectation is taken with respect to the uniformly distributed  $U \in [-1,1]^d$.
For numerical convenience, we shall simplify the expression \eqref{eq:prod} by using the approximation
\be\label{eq:inters2fb}
(1-t)^n \simeq e^{-nt}\, ,
\ee
where $t=\mathbb{P}_{_Z} \left\{ \|U - {Z} \| \leq r \right\}$. This approximation is very accurate for small   values of $t$ and moderate values of $nt$, which
is always the case of our interest.
Combining \eqref{eq:prod}, \eqref{eq:prod5} and \eqref{eq:inters2fb}, we obtain the approximation
\be
\label{eq:CA}
C_d(\mathbb{Z}_n,r) \simeq 1- \mathbb{E}_{_U} \exp ( -n \cdot \mathbb{P}_{_Z} \left\{ \|U - {Z} \| \leq r \right\}   )\, .
\ee
In the next section we will formulate three schemes that will be of theoretical interest in this paper. For each scheme and hence different distribution of $Z$, we shall derive accurate approximations for $\mathbb{P}_{_Z} \left\{ \|U - {Z} \| \leq r \right\} $ and therefore, using \eqref{eq:CA}, for $C_d(\mathbb{Z}_n,r)$.

\subsection{Designs of theoretical interest}
The three designs that will be the focus of theoretical investigation in this paper are:\\

{\bf Design 1.} {\it
$Z_1, \ldots, Z_n \in \mathbb{Z}_n$ are i.i.d. random vectors on $[-\delta, \delta]^d$ with independent components  distributed according to the following Beta$_\delta(\alpha,\alpha)$ distribution with density:
}

\be\label{Beta_density}
p_{\alpha,\delta}(t)= \frac{(2\delta)^{1-2\alpha}}{\mbox{Beta$(\alpha,\alpha)$}} [\delta^2-t^2]^{\alpha-1}\, , \;\;-\delta<t<\delta\, ,\text{ for some }\alpha> 0\;\;{\rm and }\; 0\leq \delta\leq 1.
\ee

{\bf Design 2a.} {\it
$Z_1, \ldots, Z_n \in \mathbb{Z}_n$ are i.i.d. random vectors  obtained by sampling with replacement from the vertices of the cube $[-\delta, \delta]^d$.
}

{\bf Design 2b.} {\it
$Z_1, \ldots, Z_n \in \mathbb{Z}_n$ are  random vectors  obtained by sampling without replacement from the vertices of the cube $[-\delta, \delta]^d$.}  \\

All three designs above  are nested so that $ \mathbb{Z}_n \subset  \mathbb{Z}_{n+1}$ for all eligible $n$. Designs 1 and 2a are defined for all $n=1,2, \ldots$ whereas Design 2b is  defined for  $n=1,2, \ldots, 2^d$. The appealing property of any design whose points $Z_i$ are i.i.d. is the possibility of using   \eqref{eq:prod}; this is the case of Designs 1 and 2a. For Design~2b, we will need to make some adjustments, see Section~\ref{sec:des2b}.

In the case of $\alpha=1$ in Design 1, the distribution Beta$_\delta(\alpha,\alpha)$  becomes uniform on $[-\delta,\delta]^d$.
This case has been comprehensively studied in \cite{us} with a number of approximations for $C_d(\mathbb{Z}_n,r)$ being developed.
The approximations developed  in Section~\ref{sec:des1} are generalizations of the approximations of \cite{us}.
Numerical results of \cite{us} indicated that Beta-distribution with $\alpha<1$ provides more efficient covering schemes; this explains
 the importance of the approximations of Section~\ref{sec:des1}.
Design 2a  is the limiting form of Design 1 as $\alpha \rightarrow 0$. Theoretical approximations developed below for $C_d(\mathbb{Z}_n,r)$ for Design 2a are, however, more precise than the limiting cases of approximations  obtained for $C_d(\mathbb{Z}_n,r)$  in case of Design 1.
 For numerical comparison, in Section \ref{sec:numeric} we shall also consider several other designs.\\

\section{Approximation of $C_d(\mathbb{Z}_n,r)$ for Design 1 }
\label{sec:des1}

As a result of \eqref{eq:CA}, our main quantity of interest in this section will be the probability
\be \label{eq:prob_u}
P_{U,\delta,\alpha,r}:=\mathbb{P}_{_Z} \left\{ \| U\!-\!Z \|\! \leq \! { r } \right\}\!= \! \mathbb{P}_{_Z} \left\{ \| U\!-\!Z \|^2 \leq  { r^2 } \right\}\!= \! \mathbb{P} \left\{\sum_{j=1}^d (u_j\!-\!z_j)^2 \leq  { r }^2 \right\}   \;\;
\ee
in the case when $Z$ has  the Beta-distribution with density \eqref{Beta_density}. We shall develop a simple approximation based on the Central Limit Theorem (CLT) and then subsequently refine it using the general expansion in the CLT for sums of independent non-identical r.v.

\subsection{Normal approximation for $P_{U,\delta,\alpha,r}$ }\label{sec:normal_approx}

Let $\eta_{u,\delta,\alpha} = (z-u)^2$, where $z$ has density \eqref{Beta_density}. In view of Lemma~\ref{lem:1}, the r.v. $\eta_{u,\delta,\alpha} $ is concentrated on the interval $ [(\max(0, \delta-|u|) )^2,(\delta+|u| )^2]$ and its first three  central moments  are:
\be
\label{1eq:inters1c}
\mu_{u}^{(1)} &=& \mathbb{E}\eta_{u,\delta,\alpha} =u^2+ \frac {{{\delta}}^{2}}{2\,{\alpha}+1} \, ,\\
\label{1eq:inters1c2}
\mu_{u}^{(2)}&=&{\rm var} (\eta_{u,\delta,\alpha}) =
{\frac {4\delta^{2}}{2\,{\alpha}+1}}
 \left[u^2+
 {\frac {{{\delta}}^{2}{\alpha}}{ \left( 2\,{\alpha}+1
 \right)  \left( 2\,{\alpha}+3 \right) }}
  \right]  \, ,\\
\label{1eq:inters1c3}
\mu_{u}^{(3)}&=&\mathbb{E} \left[\eta_{u,\delta,\alpha} - \mu_{u}^{(1)}\right]^3 =
{\frac {48{\alpha}\,{{\delta}}^{4}}{ \left( 2\,{\alpha}+1
 \right) ^{2} \left( 2\,{\alpha}+3 \right) }} \left[  u^2+ {\frac {{{\delta}}^{2} \left( 2\,{\alpha} -1\right) }{3
 \left( 2\,{\alpha}+5 \right)  \left( 2\,{\alpha}+1 \right) }}
 \right]
\, .
\ee

For a given $U=(u_1, \ldots, u_d) \in \mathbb{R}^d$, consider the r.v.
\bea
\| U-Z \|^2 =\sum_{i=1}^d \eta_{u_i,\delta,\alpha} =\sum_{j=1}^d (u_j-z_j)^2\, ,
\eea
where we assume that $Z=(z_1, \ldots, z_d) $ is a random vector with i.i.d. components $z_i$ with density \eqref{Beta_density}.
From \eqref{1eq:inters1c}, its mean is
\bea
\mu=\mu_{d,\delta,\alpha,U}:=\mathbb{E}\| U-Z \|^2  =\|U\|^2 +\frac {{d{\delta}}^{2}}{2\,{\alpha}+1}\, .
\eea
Using independence of $z_1, \ldots, z_d$ and  \eqref{1eq:inters1c2}, we obtain
\bea
 {\sigma}_{d,\delta,\alpha,U}^2 :={\rm var}(\| U-Z \|^2 )  =
{\frac {4\delta^{2}}{2\,{\alpha}+1}}
 \left[\|U\|^2+
 {\frac {{d{\delta}}^{2}{\alpha}}{ \left( 2\,{\alpha}+1
 \right)  \left( 2\,{\alpha}+3 \right) }}
  \right]  \, ,
\eea
and from independence of $z_1, \ldots, z_d$ and  \eqref{1eq:inters1c3} we get
\bea
  {\mu}_{d,\delta,\alpha,U}^{(3)} := \mathbb{E}\left[\| U-Z \|^2- \mu\right]^3  = \sum_{j=1}^d     \mu_{u_j}^{(3)} = {\frac {48{\alpha}\,{{\delta}}^{4}}{ \left( 2\,{\alpha}+1
 \right) ^{2} \left( 2\,{\alpha}+3 \right) }} \left[  \|U\|^2+ {\frac {{d{\delta}}^{2} \left( 2\,{\alpha} -1\right) }{3
 \left( 2\,{\alpha}+5 \right)  \left( 2\,{\alpha}+1 \right) }}
 \right]
\, .\;\;\;\;\;\;
\eea

If $d$ is large enough then the conditions of the CLT for $\| U-Z \|^2$ are approximately met and  the distribution of $\| U-Z \|^2 $
 is approximately normal with mean $\mu_{d,\delta,\alpha,U}$ and variance ${\sigma}_{d,\delta,\alpha,U}^2$. That is, we can approximate the probability
$P_{U,\delta,\alpha,r}= \mathbb{P}_{_Z} \left\{ \| U\!-\!Z \|\! \leq \! { r } \right\}$
by
\be
\label{1eq:inters2f1}
P_{U,\delta,\alpha,r}\!\cong \Phi \left(\frac{{ r }^2-\mu_{d,\delta,\alpha,U}}{{\sigma}_{d,\delta,\alpha,U}} \right) \, ,
\ee
where $\Phi (\cdot)$   is the c.d.f. of the standard normal distribution:
$$
\Phi (t) = \int_{-\infty}^t \varphi(v)dv\;\;{\rm with}\;\; \varphi(v)=\frac{1}{\sqrt{2\pi}} e^{-v^2/2}\, .
$$
The approximation \eqref{1eq:inters2f1} has acceptable accuracy if the probability $P_{U,\delta,\alpha,r}$ is not very small; for example, it falls inside a $2\sigma$-confidence interval generated by the standard normal distribution.
In the next section, we improve approximations \eqref{1eq:inters2f1}  by
using an Edgeworth-type  expansion in the CLT for sums of independent  non-identically distributed r.v.

\subsection{Refined approximation for $P_{U,\delta,\alpha,r}$ }
\label{sec:quantuty2I}

General expansion in the central limit theorem for sums of independent non-identical r.v. has been derived
by V.Petrov, see
Theorem 7 in Chapter 6 in  \cite{petrov2012sums}, see also Proposition 1.5.7 in \cite{rao1987asymptotic}.
 The first three terms of this expansion have been specialized by V.Petrov in
Section 5.6 in \cite{petrov}.
By using only the first term in this expansion,
we obtain the following approximation for the distribution function of $\| U-Z \|^2 $:
\be\label{original_petrov}
\mathbb{P}\left(\frac{\| U-Z \|^2-\mu_{d,\delta,\alpha,U}}{\sigma_{d,\delta,\alpha,U}} \leq x \right) \cong \Phi(x) + \frac{ \mu_{d,\delta,\alpha,U}^{(3)}}{6  {\sigma}_{d,\delta,\alpha,U}^3 }(1-x^2)\varphi(x),
\ee
leading to the following improved  form of \eqref{1eq:inters2f1}:

\be
\label{1eq:inters2f_corrected_1}
P_{U,\delta,\alpha,r} \cong \Phi(t) + \frac{\alpha \delta \left[ \|U\|^2+\frac{d\delta^2(2\alpha-1)}{3(2\alpha+5)(2\alpha+1)} \right]}{ (2\alpha+3)(2\alpha+1)^{1/2}\left [\|U\|^2+\frac{d\delta^2\alpha}{(2\alpha+1)(2\alpha+3)}\right]^{3/2} }(1-t^2)\varphi(t) \, ,
\ee
where
\bea
t := \frac{{ r }^2-\mu_{d,\delta,\alpha,U}}{{\sigma}_{d,\delta,\alpha,U}}
=  \frac{\sqrt{2\alpha+1}(r^2- \|U\|^2 -\frac{d\delta^2}{2\alpha+1})}{2\delta\sqrt{ \|U\|^2 +\frac{{d\delta^2\alpha}}{(2\alpha+1)(2\alpha+3)}} }\,.
\eea

For $\alpha=1$, we obtain
\bea
\label{eq:inters2f_corrected1}
P_{U,\delta,\alpha,r} \cong \Phi(t) + \frac{ \delta \left[ \|U\|^2+{d\delta^2}/{63} \right]}{ 5\sqrt{3}\left [\|U\|^2+{d\delta^2}/{15}\right]^{3/2} }(1-t^2)\varphi(t) \,\; \;{\rm  with }\;\; t
=  \frac{\sqrt{3}(r^2- \|U\|^2 -{d\delta^2}/{3})}{2\delta\sqrt{ \|U\|^2 +{{d\delta^2}}/{15}} }\,,
\eea
which coincides with formula (16) of \cite{us}.

A very attractive feature of the approximations \eqref{1eq:inters2f1} and  \eqref{1eq:inters2f_corrected_1} is their dependence on $U$ through $\|U\|$ only. We could have specialized for our case the next terms in Petrov's approximation but these terms no longer depend on  $\|U\|$ only and hence the next terms are much more complicated. Moreover,  adding one or two extra terms from Petrov's expansion to the approximation \eqref{1eq:inters2f_corrected_1}  does not fix the problem entirely for all $U$, $\delta$, $\alpha$ and $r$. Instead,
we propose a slight adjustment to the r.h.s of \eqref{1eq:inters2f_corrected_1} to improve this approximation, especially for small dimensions. Specifically, we suggest  the approximation
\be
\label{2eq:inters2f_corrected}
P_{U,\delta,\alpha,r} \cong \Phi(t) +c_{d,\alpha} \frac{\alpha \delta \left[ \|U\|^2+\frac{d\delta^2(2\alpha-1)}{3(2\alpha+5)(2\alpha+1)} \right]}{ (2\alpha+3)(2\alpha+1)^{1/2}\left [\|U\|^2+\frac{d\delta^2\alpha}{(2\alpha+1)(2\alpha+3)}\right]^{3/2} }(1-t^2)\varphi(t) \, ,\ee
where $c_{d,\alpha}= 1+ {3}/{(\alpha d)}$.

Below, there are  figures of two types. In Figures~\ref{Overall1}--\ref{Overall2}, we plot  $P_{U,\delta,\alpha,r} $ over a wide range of $r$ ensuring that values of $P_{U,\delta,\alpha,r} $  lie in the whole range $[0,1]$. In Figures~\ref{One_ball_approx1}--\ref{One_ball_approx2}, we plot $P_{U,\delta,\alpha,r}$ over a much smaller range of $r$ with $P_{U,\delta,\alpha,r}$ lying roughly in the range $[0,0.02]$. For the purpose of using formula \eqref{eq:prod}, we need to assess the accuracy of all approximations for smaller values of $P_{U,\delta,\alpha,r}$ and hence the second type of plots are more useful.
In these figures, the solid black line depicts $P_{U,\delta,\alpha,r} $ obtained via Monte Carlo methods where for simplicity we have set $U = (1/2,1/2,\ldots,1/2)$ and $\delta=1/2$. Approximations \eqref{1eq:inters2f1} and \eqref{2eq:inters2f_corrected} are depicted with a dotted blue and dash green line respectively. From numerous simulations and these figures, we can conclude the following. Whilst the basic normal approximation \eqref{1eq:inters2f1} seems adequate in the whole range of values of $r$, for particularly small probabilities, that we are most interested in, approximation \eqref{2eq:inters2f_corrected} is much superior and appears to be very accurate for all values of $\alpha$.

\begin{figure}[h]
\centering
\begin{minipage}{.5\textwidth}
  \centering
  \includegraphics[width=1\linewidth]{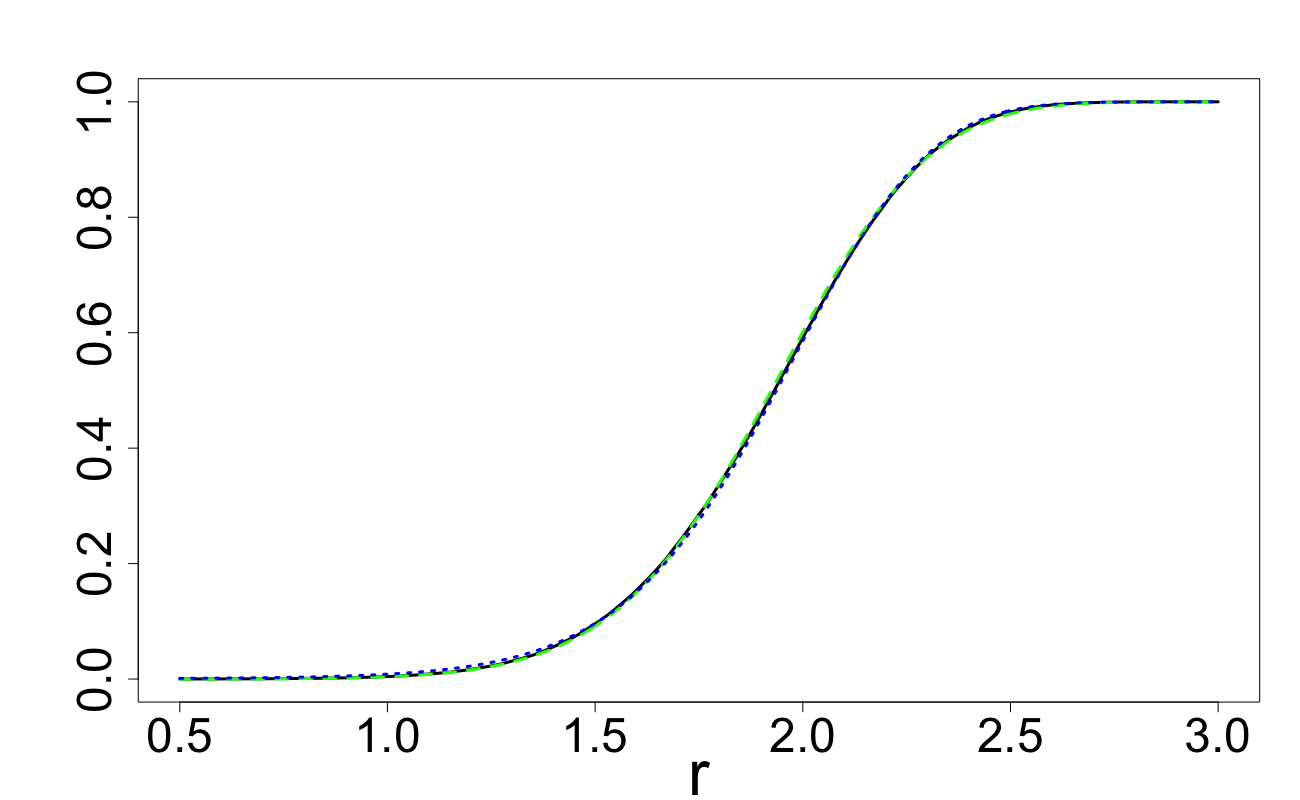}
  \caption{$P_{U,\delta,\alpha,r} $ and approximations: $d=10$,\\ $\alpha=0.5$.}
  \label{Overall1}
\end{minipage}%
\begin{minipage}{.5\textwidth}
  \centering
  \includegraphics[width=1\linewidth]{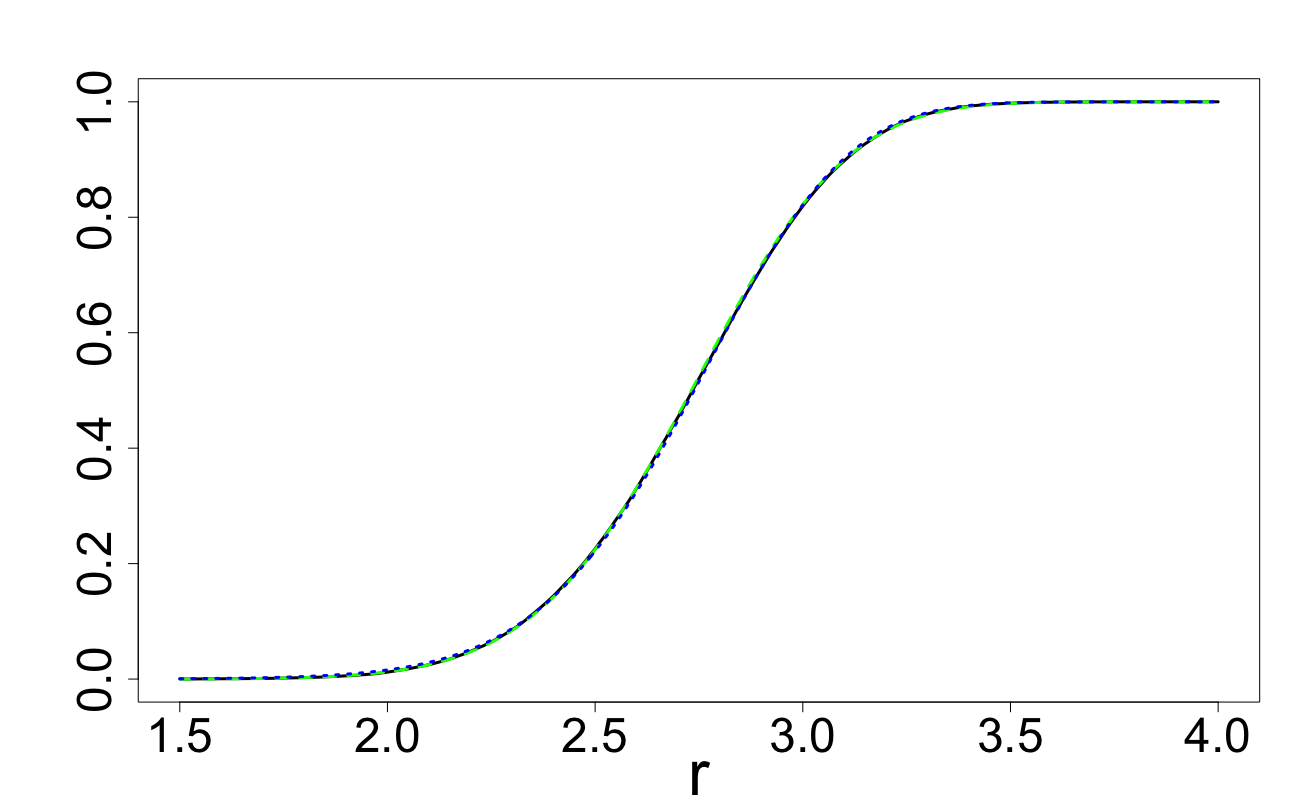}
  \caption{$P_{U,\delta,\alpha,r} $ and approximations: $d=20$,\\ $\alpha=0.5$.}
  \label{Overall2}
\end{minipage}
\end{figure}
\begin{figure}[h]
\centering
\begin{minipage}{.5\textwidth}
  \centering
  \includegraphics[width=1\linewidth]{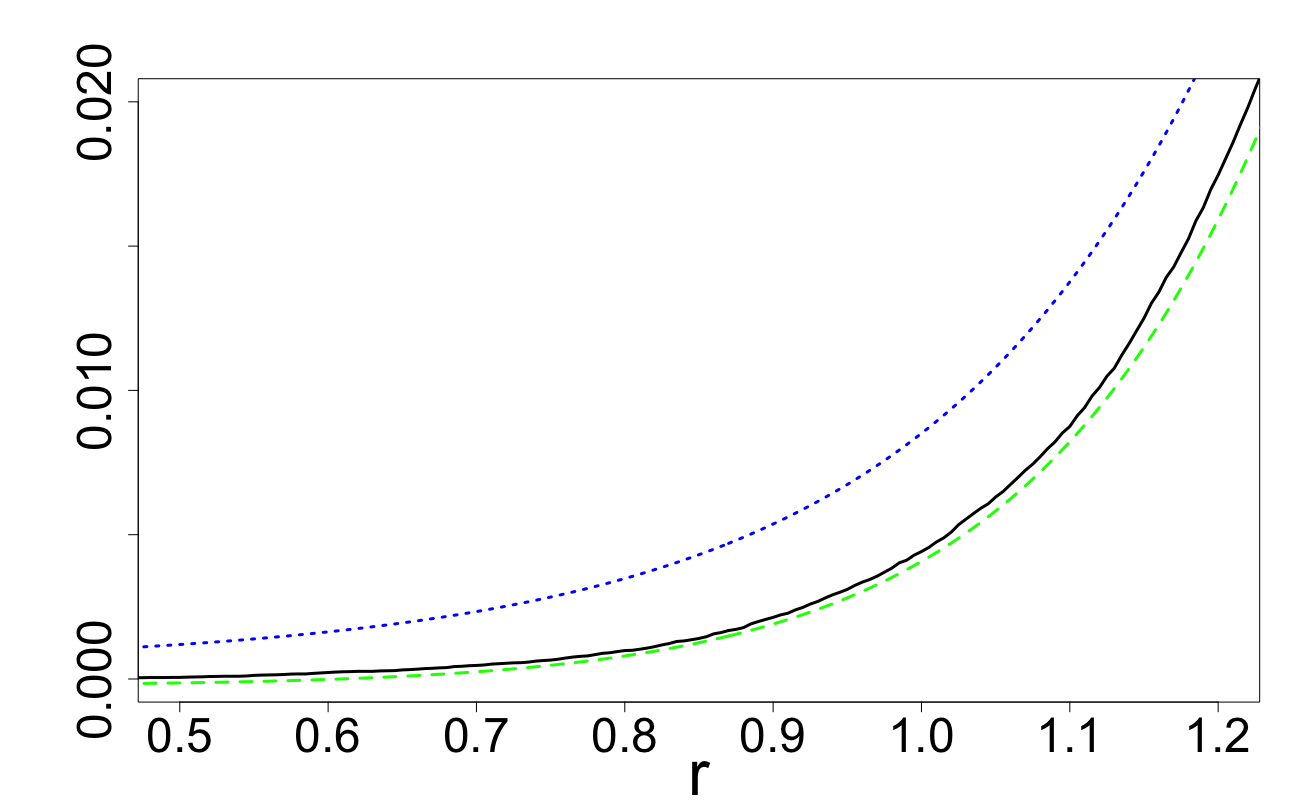}
  \caption{$P_{U,\delta,\alpha,r} $ and approximations: $d=10$,\\ $\alpha=0.5$.}
  \label{One_ball_approx1}
\end{minipage}%
\begin{minipage}{.5\textwidth}
  \centering
  \includegraphics[width=1\linewidth]{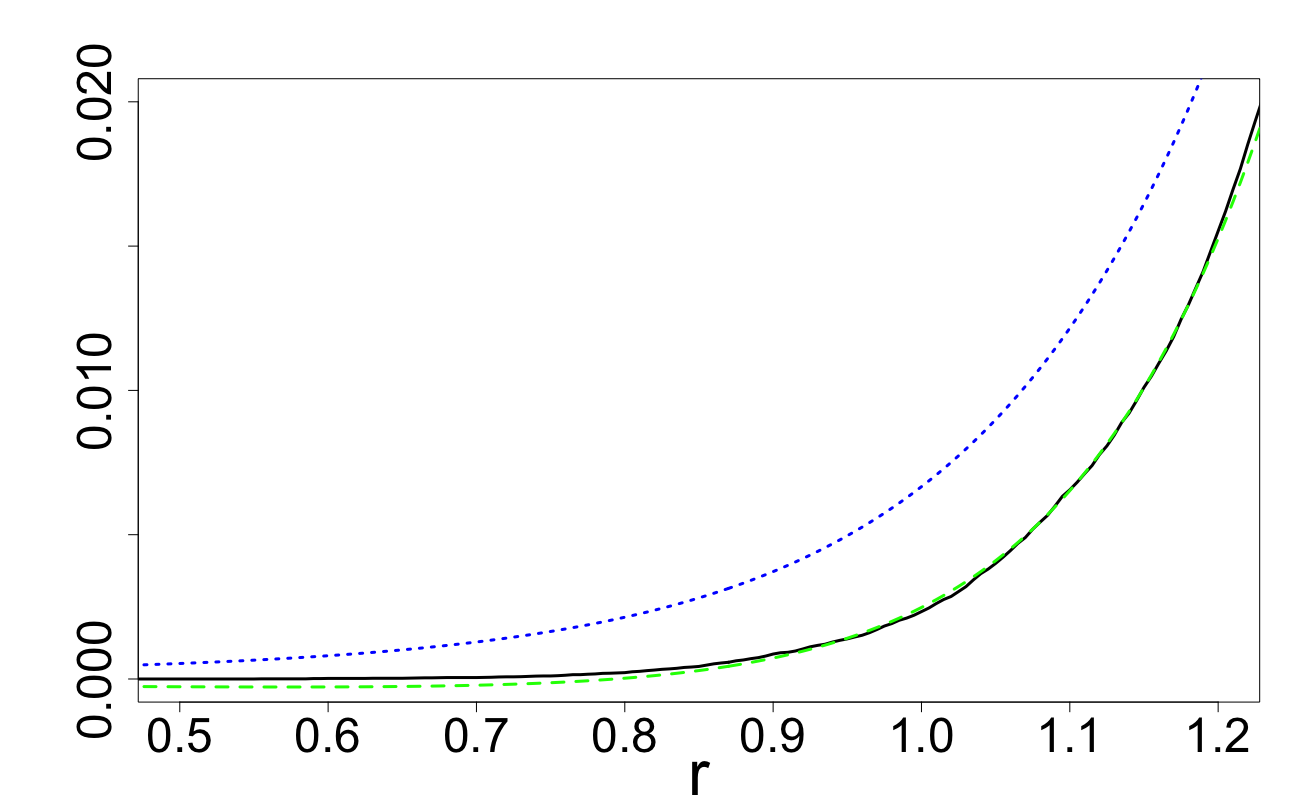}
  \caption{$P_{U,\delta,\alpha,r} $ and approximations: $d=10$,\\ $\alpha=1$.}
\end{minipage}
\end{figure}
\begin{figure}[H]
\centering
\begin{minipage}{.5\textwidth}
  \centering
  \includegraphics[width=1\linewidth]{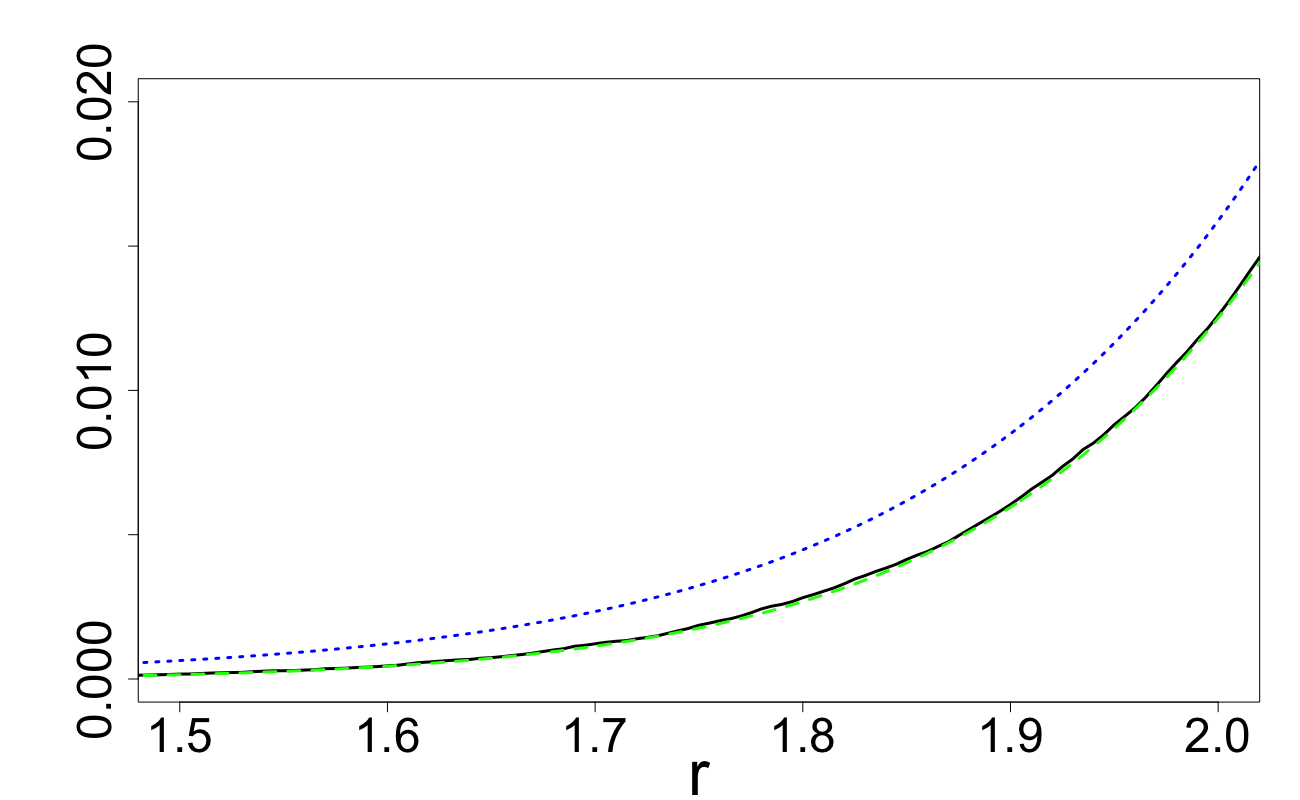}
  \caption{$P_{U,\delta,\alpha,r}$ and approximations: $d=20$,\\ $\alpha=0.5$.}
  \label{}
\end{minipage}%
\begin{minipage}{.5\textwidth}
  \centering
  \includegraphics[width=1\linewidth]{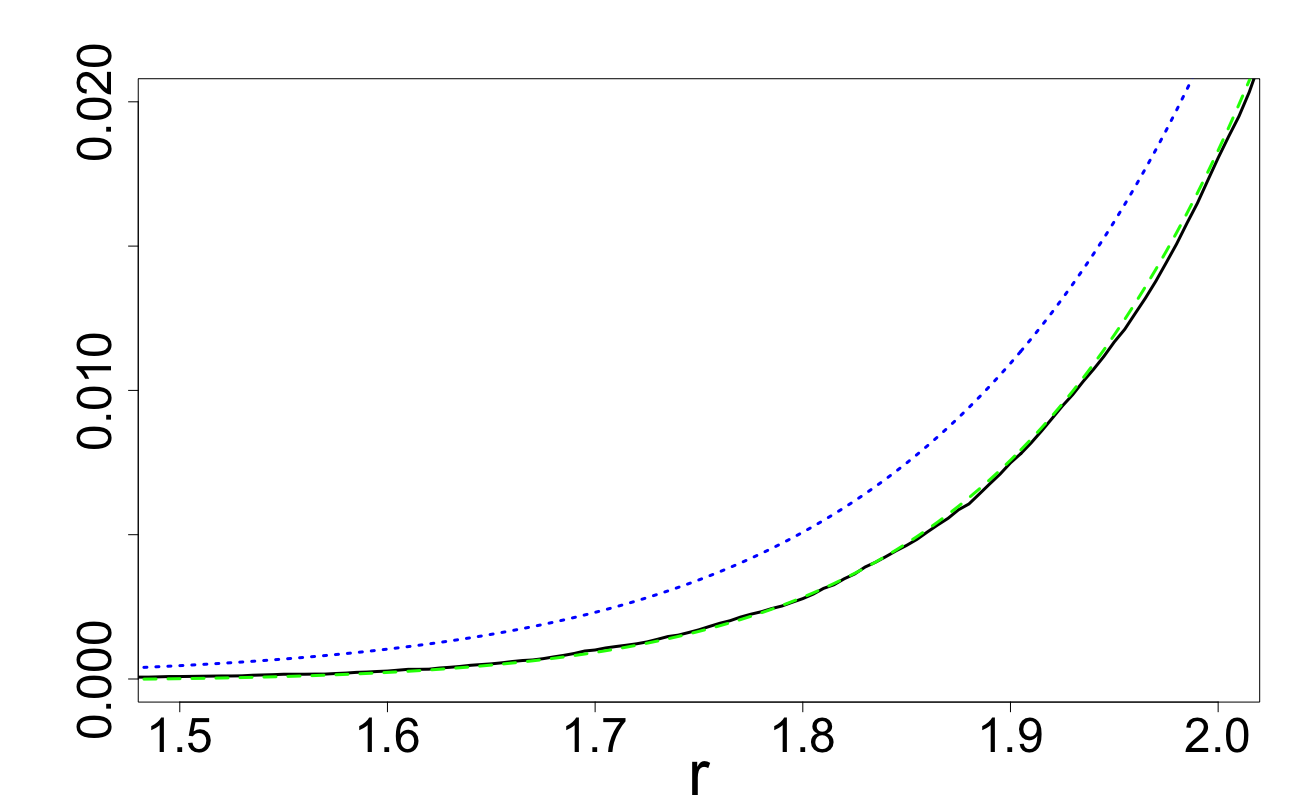}
  \caption{$P_{U,\delta,\alpha,r} $ and approximations: $d=20$,\\ $\alpha=1$.}
  \label{One_ball_approx2}
\end{minipage}
\end{figure}

\subsection{Approximation for $C_d(\mathbb{Z}_n,r)$ for Design 1}\label{approx_alpha_gt_0}

Consider now  $C_d(\mathbb{Z}_n,r)$ for Design 1,
as expressed via
$P_{U,\delta,\alpha,r}$ in \eqref{eq:CA}.
As $U$ is uniform on $[-1,1]^d$,
$\mathbb{E} \| U \|^2= d/3$ and $ {\rm var}( \| U \|^2)= {4d}/{45}  .$ Moreover, if $d$ is large enough then $\| U \|^2= \sum_{j=1}^d u_j^2$
is approximately normal.

We will combine the expressions \eqref{eq:CA}  with approximations \eqref{1eq:inters2f1} and \eqref{2eq:inters2f_corrected}   as well as with the normal approximation for the distribution of $\| U \|^2$, to arrive at two final approximations for $C_d(\mathbb{Z}_n,r)$ that differ in complexity.
If the original normal approximation \eqref{1eq:inters2f1} of $P_{U,\delta,\alpha,r}$ is used then we obtain:

\be
C_d(\mathbb{Z}_n,r) \simeq 1- \int_{-\infty}^{\infty} \psi_{1,\alpha}(s) \varphi(s)d s   \label{eq:prod5a}\;\;
\ee
with
\bea
 \!  \psi_{1,\alpha}(s)=\exp \left\{-n
\Phi(c_s) \right\} , \,\ c_s= \frac{(2\alpha+1)^{1/2}\left(r^2 \!-\!s' \!-\! \frac{d\delta^2}{2\alpha+1} \right)}{2\delta\sqrt{ s'+\kappa  }}, \,\, s'=s\sqrt{\frac{4d}{45}}+d/3 , \,\, \kappa = \frac{d\delta^2\alpha}{(2\alpha\!+\!1)(2\alpha\!+\!3)}\, .\;\;
\eea

 If the approximation \eqref{2eq:inters2f_corrected} is used,  we obtain:
\be
\label{eq:accurate_app}
C_d(\mathbb{Z}_n,r) \simeq 1- \int_{-\infty}^{\infty} \psi_{2,\alpha}(s) \varphi(s)d s, \;\; \label{eq:prod6a}\;\;
\ee
with
\bea
   \psi_{2,\alpha}(s)=\exp \left\{-n \left(
\Phi(c_s) + c_{d,\alpha} \frac{\alpha \delta \left[ s'+\frac{d\delta^2(2\alpha-1)}{3(2\alpha+5)(2\alpha+1)} \right]}{ (2\alpha+3)(2\alpha+1)^{1/2}\left [s'+\kappa\right]^{3/2} }(1-c_s^2)\varphi(c_s)   \right) \right\} \,.
\eea

For $\alpha=1$, we get
\be \label{eq:psi2-1}
   \psi_{2,1}(s)=\exp \left\{-n \left(
\Phi(c_s) + c_{d,\alpha} \frac{ \delta \left[ s'+\frac{d\delta^2}{63} \right]}{ 5\sqrt{3}\left [s'+\frac{d\delta^2}{15}\right]^{3/2} }(1-c_s^2)\varphi(c_s)   \right) \right\} \,
\ee
and the approximation    \eqref{eq:accurate_app} coincides with the approximation (26) in \cite{us}.
The accuracy of approximations  \eqref{eq:prod5a}   and   \eqref{eq:accurate_app}  will be assessed in Section~\ref{sec:numeric1}.

\section{Approximating $C_d(\mathbb{Z}_n,r)$ for Design 2a}
\label{sec:des2a}

Our main quantity  of interest in this section will be the probability $P_{U,\delta,0,r}$ defined in \eqref{eq:prob_u}
in the case where components $z_i$ of the vector $Z=(z_1, \ldots, z_d) \in \mathbb{R}^d$ are i.i.d.r.v with ${\rm Pr}(z_i=\delta)= {\rm Pr}(z_i=-\delta)= 1/2$; this is a limiting case of $P_{U,\delta,\alpha,r}$ as $\alpha \to 0$.

\subsection{Normal approximation for $P_{U,\delta,0,r}$ }

Using the same approach that led to approximation \eqref{1eq:inters2f1} in Section~\ref{sec:normal_approx}, the initial normal approximation for $P_{U,\delta,0,r}$ is:
\be
\label{1eq:alph_zero}
P_{U,\delta,0,r}\!\cong \Phi \left(\frac{{ r }^2-\mu_{d,\delta,U}}{{\sigma}_{d,\delta,U}} \right) \, ,
\ee
where, from Lemma~\ref{lem:1}, we have
\bea
\mu_{d,\delta,U}=\|U\|^2 +d{\delta}^{2} \,\, \text{ and }  \,\,\  {\sigma}_{d,\delta,U}^2  =
4\delta^{2}\|U\|^2\, .
\eea

\subsection{Refined approximation for $P_{U,\delta,0,r}$ }\label{sec:General_expansion}

From \eqref{eq:mom_0},
 we have $   {\mu}_{d,\delta,\alpha,U}^{(3)}=0$ and therefore the last term in the rhs of \eqref{original_petrov} with $\alpha=0$ is no longer present. By  taking an additional term in the general expansion, see V.Petrov in
Section 5.6 in \cite{petrov}, we obtain the following approximation for the distribution function of $\| U-Z \|^2$:
\be\label{new_petrov}
\mathbb{P}\left(\frac{\| U-Z \|^2-\mu_{d,\delta,U}}{\sigma_{d,\delta,U}} \leq x \right) \cong \Phi(x)  - (x^3-3x)\frac{\kappa_{d,\delta,0,U}^{(4)}}{24{\sigma}_{d,\delta,0,U}^4 }\varphi(x),
\ee
where $\kappa_{d,\delta,0,U}^{(4)}$ is the sum of $d$ fourth cumulants of the centred r.v. $(z-u)^2$, where $z$ is concentrated at two points $\pm \delta$ with
${\rm Pr}(z=\pm\delta)=  1/2$. From \eqref{eq:mom_0},
\bea
\kappa_{d,\delta,0,U}^{(4)}:= \sum_{j=1}^d     (\mu_{u_j}^{(4)}-3[\mu_{u_j}^{(2)}]^2 ) =  -32\delta^4 \sum_{i=1}^{d}u_i^4 \, .
\eea

Unlike \eqref{original_petrov}, the rhs of \eqref{new_petrov} does not depends solely on $\|U\|^2$. However, the quantities  $\|U\|^2$ and   $\sum_{i=1}^{d}u_i^4$ are strongly correlated; one can show that for all $d$
\bea
 {\rm corr}\left(\|U\|^2,  \sum_{i=1}^{d}u_i^4 \right) = \frac{3\sqrt{5}}{7  } \cong 0.958 \,.
\eea

This suggests  (by rounding the correlation above to 1) the following approximation:
\bea
\sum_{i=1}^{d}u_i^4 \cong \frac{4\sqrt{d}}{15} \left ( \frac{\|U\|^2-d/3}{\sqrt{\frac{4d}{45}}} \right) + \frac{d}{5}\, .
\eea
With this approximation,  the rhs of \eqref{new_petrov}  depends only on $\|U\|^2$. As a result, the following refined form of \eqref{1eq:alph_zero} is:
\bea\label{alpha_zero}
P_{U,\delta,0,r} \cong \Phi(t)  + (t^3-3t)\frac{2(\|U\|^2-d/3)/\sqrt{5} +d/5}{12\|U\|^4} \varphi(t),
\eea
where
\bea
t := \frac{{ r }^2-\mu_{d,\delta,0,U}}{{\sigma}_{d,\delta,0,U}}
=  \frac{(r^2- \|U\|^2 -d\delta^2)}{2\delta\|U\| }\,.
\eea

Similarly to approximation \eqref{2eq:inters2f_corrected}, we propose a slight adjustment to the r.h.s of the approximation above:
\be\label{alpha_zero_corrected}
P_{U,\delta,0,r} \cong \Phi(t)  + \left(1+\frac3d \right)(t^3-3t)\frac{2(\|U\|^2-d/3)/\sqrt{5} +d/5}{12\|U\|^4} \varphi(t).
\ee

In the same style as at the end of Section~\ref{sec:quantuty2I}, below there are figures of two types. In Figures~\ref{Overall3}--\ref{Overall4}, we plot  $P_{U,\delta,0,r} $ over a wide range of $r$ ensuring that values of $P_{U,\delta,0,r} $  lie in the range $[0,1]$. In Figures~\ref{One_ball_approx1_alph_zero}--\ref{One_ball_approx1_alph_zero2}, we plot $P_{U,\delta,0,r}$ over a much smaller range of $r$ with $P_{U,\delta,0,r}$ lying in the range $[0,0.02]$.
In these figures, the solid black line depicts $P_{U,\delta,\alpha,r} $ obtained via Monte Carlo methods where we have set $\delta=1/2$ and $U$ is a point sampled uniformly on $[-1,1]^d$; for reproducibility, in the caption of each figure we state the random seed used in R. Approximations \eqref{1eq:alph_zero} and \eqref{alpha_zero_corrected} are depicted with a dotted blue and dash green line respectively. From these figures, we can conclude the same outcome as in Section~\ref{sec:quantuty2I}. Whilst the approximation \eqref{1eq:alph_zero} is rather good overall, for small probabilities the approximation \eqref{alpha_zero_corrected} is much superior and is very accurate. Note that since
random vectors $Z_j$ are taking values on a finite set, which is the set of points $(\pm \delta, \ldots, \pm \delta)$, the probability $P_{U,\delta,0,r}$ considered as a function of $r$, is a piece-wise constant function.

\begin{figure}[h]
\centering
\begin{minipage}{.5\textwidth}
  \centering
    \includegraphics[width=1\linewidth]{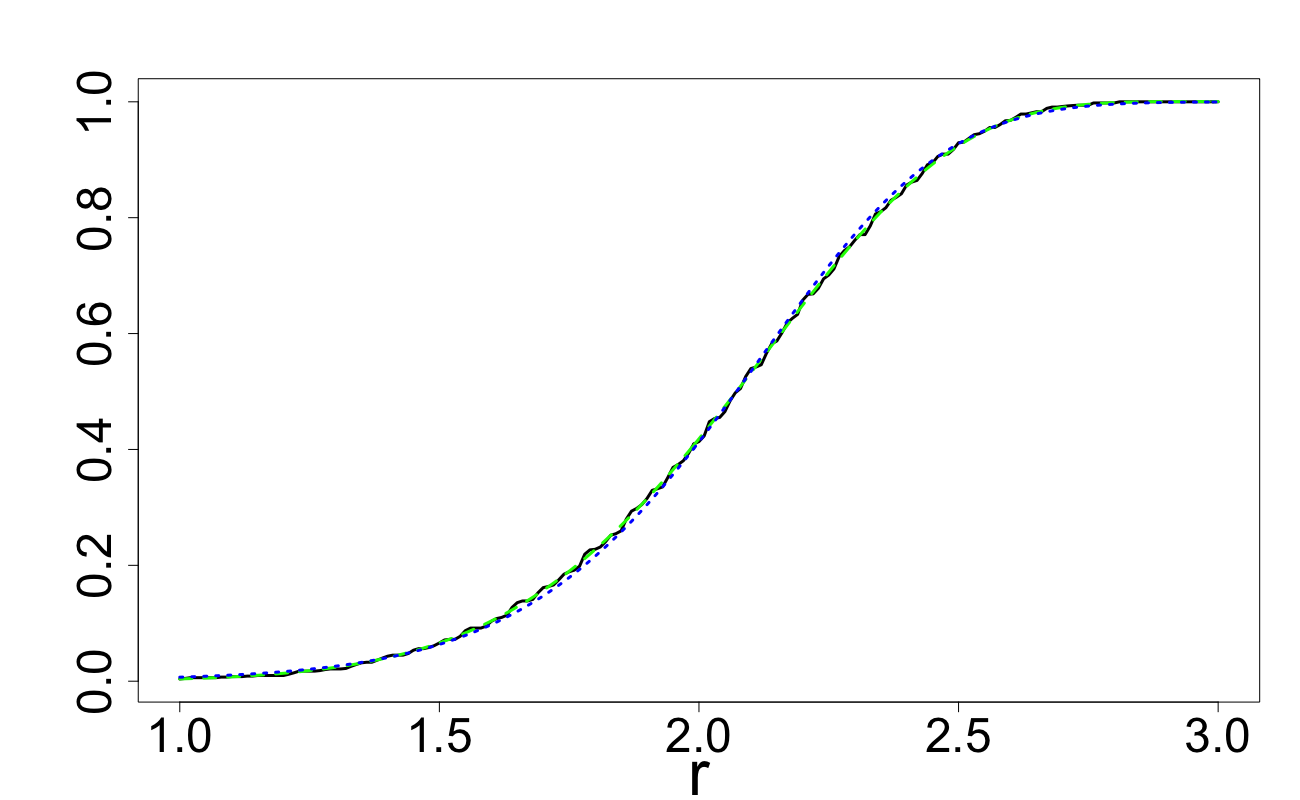}
  \caption{$P_{U,\delta,0,r} $ and approximations: $d=10$,\\ $seed=10$.}
 \label{Overall3}
\end{minipage}%
\begin{minipage}{.5\textwidth}
  \centering
\includegraphics[width=1\linewidth]{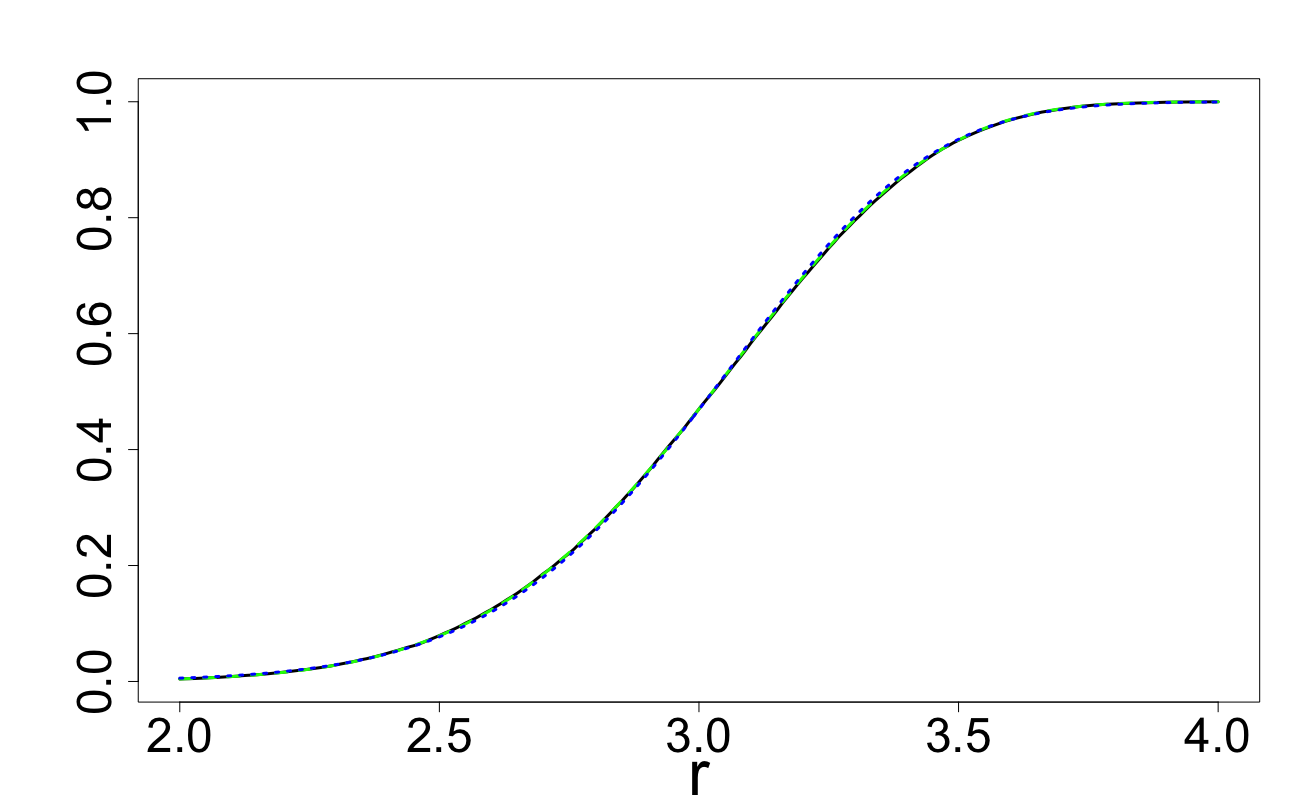}
  \caption{$P_{U,\delta,0,r}$ and approximations: $d=20$,\\ $seed =10$.}
  \label{Overall4}
\end{minipage}
\end{figure}
\begin{figure}[h]
\centering
\begin{minipage}{.5\textwidth}
  \centering
    \includegraphics[width=1\linewidth]{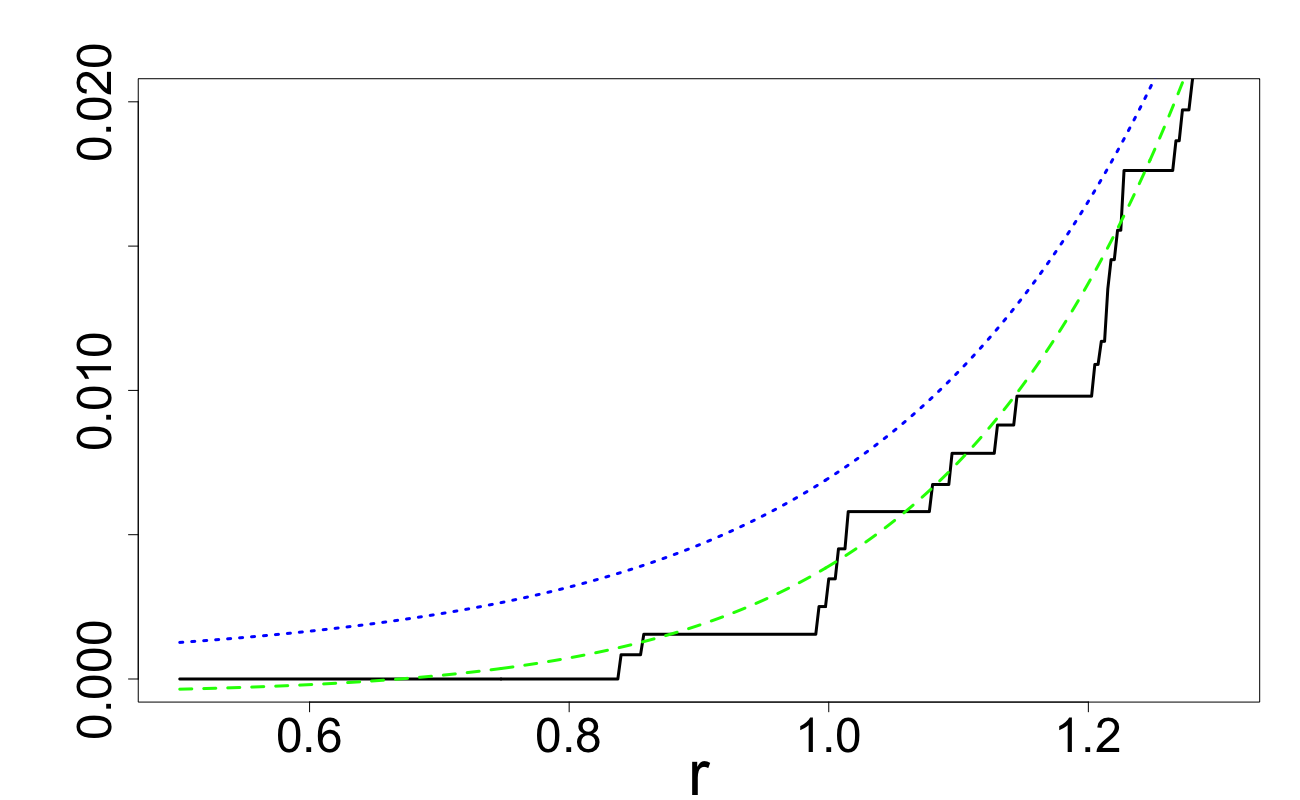}
  \caption{$P_{U,\delta,0,r} $ and approximations: $d=10$,\\$seed=10$.}
 \label{One_ball_approx1_alph_zero}
\end{minipage}%
\begin{minipage}{.5\textwidth}
  \centering
\includegraphics[width=1\linewidth]{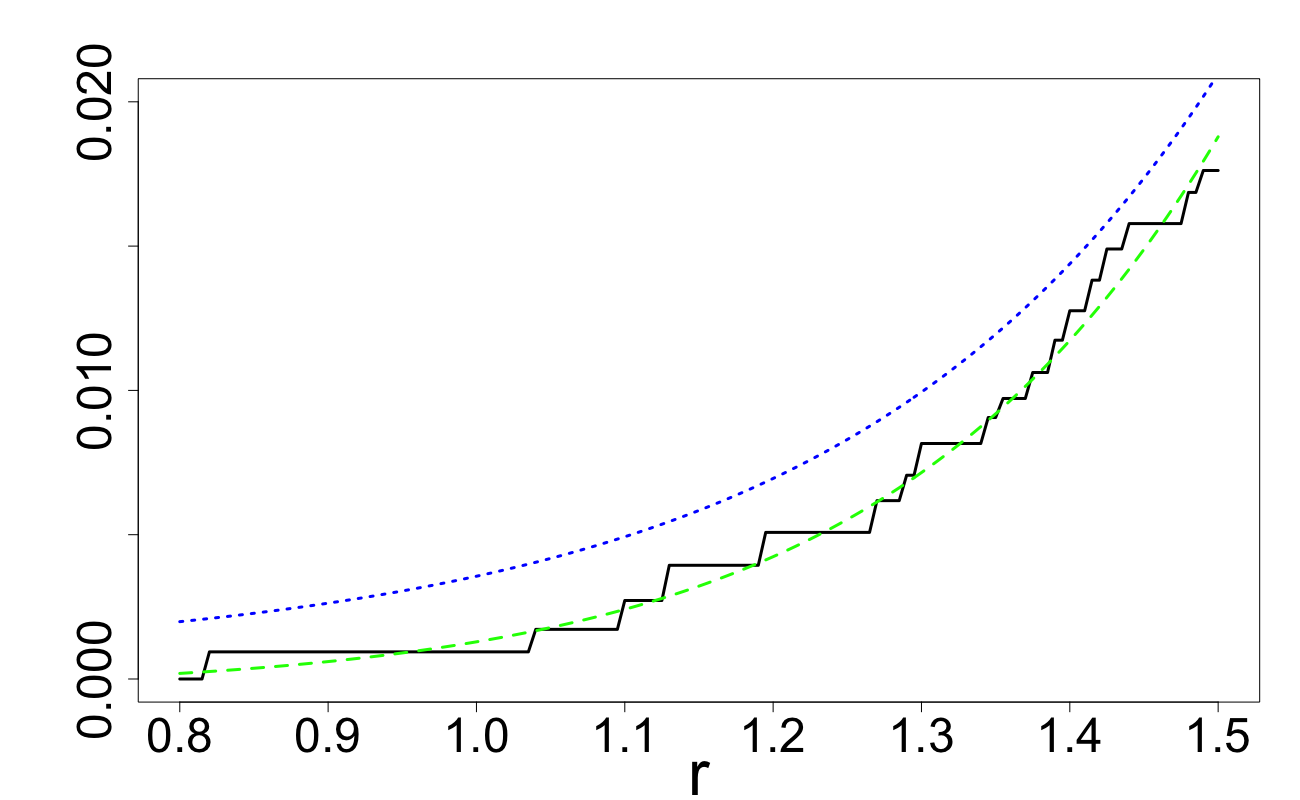}
  \caption{$P_{U,\delta,0,r}$ and approximations: $d=10$,\\$seed =15$.}
\end{minipage}
\end{figure}
\begin{figure}[!h]
\centering
\begin{minipage}{.5\textwidth}
  \centering
  \includegraphics[width=1\linewidth]{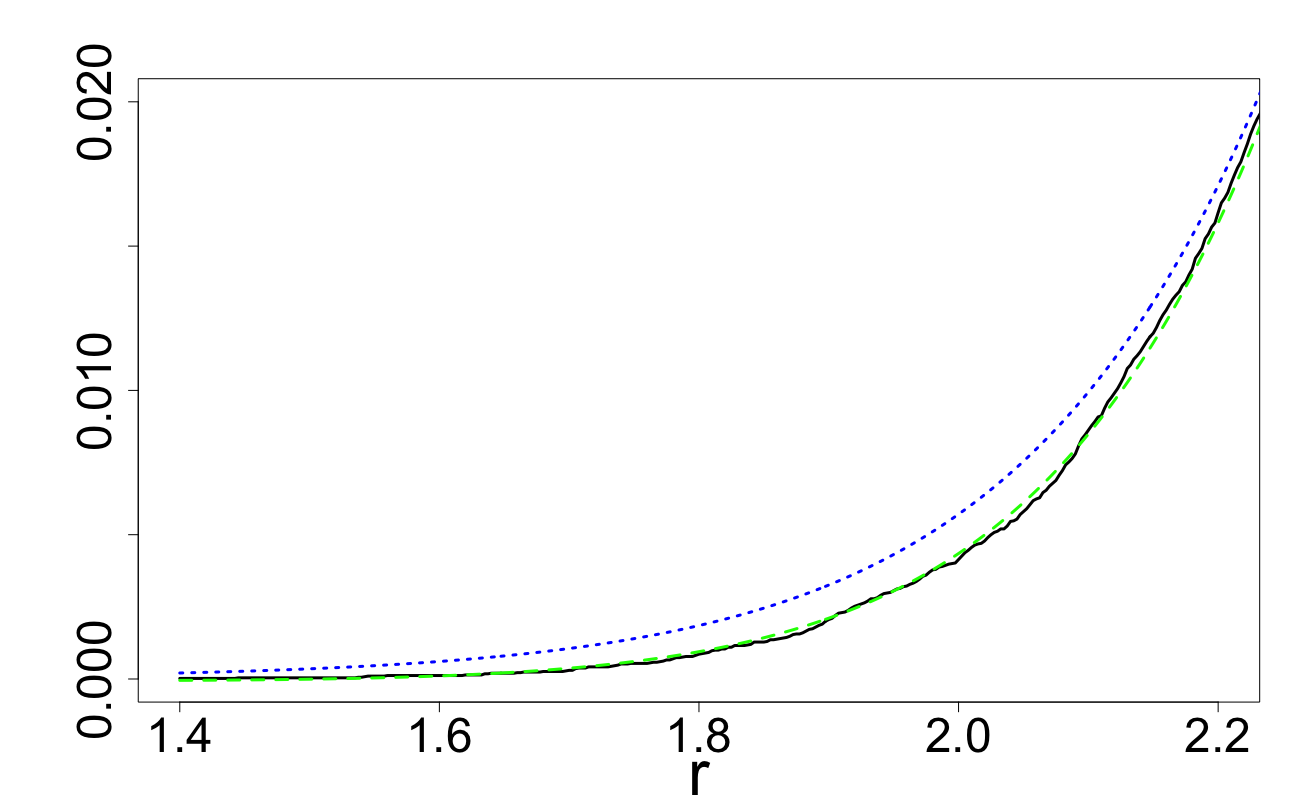}
  \caption{$P_{U,\delta,0,r} $ and approximations: $d=20$,\\ $seed=10$.}
  \label{}
\end{minipage}%
\begin{minipage}{.5\textwidth}
  \centering
  \includegraphics[width=1\linewidth]{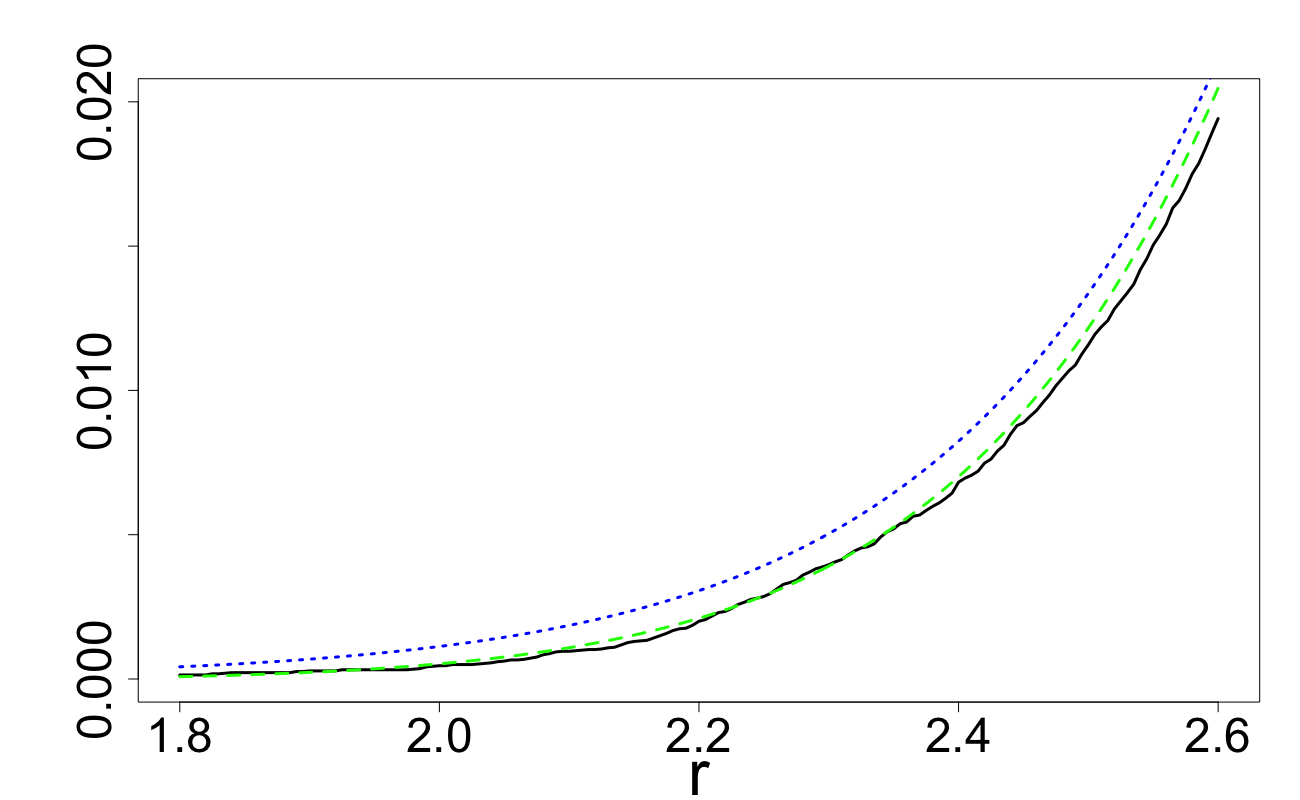}
  \caption{$P_{U,\delta,0,r} $ and approximations: $d=20$,\\ $seed=15$.}
   \label{One_ball_approx1_alph_zero2}
\end{minipage}
\end{figure}

\subsection{Approximation for $C_d(\mathbb{Z}_n,r)$}

Consider now  $C_d(\mathbb{Z}_n,r)$ for Design 2a,
as expressed via
$P_{U,\delta,\alpha,r}$ in \eqref{eq:CA}. Using the normal approximation for $\| U \|^2$ as made in the beginning of Section~\ref{approx_alpha_gt_0},
we will combine the expressions \eqref{eq:CA}  with approximations \eqref{1eq:alph_zero} and \eqref{alpha_zero_corrected} to arrive at two approximations for $C_d(\mathbb{Z}_n,r)$ that differ in complexity.

If the original normal approximation \eqref{1eq:alph_zero} of $P_{U,\delta,0,r}$ is used then we obtain:

\be \label{normal_alp_0}
C_d(\mathbb{Z}_n,r) \simeq 1- \int_{-\infty}^{\infty} \psi_{3,n}(s) \varphi(s)d s,\;\;
\ee
with
\bea
   \psi_{3,n}(s)=\exp \left\{-n
\Phi(c_s) \right\} \,, \,\,\,\, c_s= \frac{\left(r^2 -s' - {d\delta^2} \right)}{2\delta\sqrt{ s'  }}, \,\, s'=s\sqrt{\frac{4d}{45}}+d/3 \,\, .
\eea

 If the approximation \eqref{alpha_zero_corrected} is used,  we obtain:
\be
\label{eq:accurate_app2}
C_d(\mathbb{Z}_n,r) \simeq 1- \int_{-\infty}^{\infty} \psi_{4,n}(s) \varphi(s)d s, \;\; \label{eq:prod6a}\;\;
\ee
with
\be \label{eq:psi4}
   \psi_{4,n}(s)=\exp \left\{-n \left( \Phi(c_s)  +  \left(1+\frac3d \right)(c_s^3-3c_s)\frac{2(s'-d/3) /\sqrt{5} +d/5}{12(s')^2} \varphi(c_s)\right)\right\} \,.
\ee
and
\bea
\,\ c_s= \frac{\left(r^2 -s' - {d\delta^2} \right)}{2\delta\sqrt{ s'  }}, \,\, s'=s\sqrt{\frac{4d}{45}}+d/3 \,\, .
\eea

The accuracy of approximations  \eqref{normal_alp_0}   and   \eqref{eq:accurate_app2}  will be assessed in Section~\ref{sec:numeric1}.

\section{Approximating $C_d(\mathbb{Z}_n,r)$ for Design 2b}
\label{sec:des2b}

Designs whose points $Z_i$ have been sampled from a finite discrete set  without replacement have dependence, for example Design 2b, and therefore formula \eqref{eq:prod} cannot be used.

In this section, we suggest a way of modifying the approximations developed in Section~\ref{sec:des2a} for Design 2a. This will amount to
approximating sampling without replacement by a suitable sampling with replacement.

\subsection{Establishing a connection between  sampling with and without replacement: general case }\label{sec:with_replacement}

Let $\cal{S}$  be a discrete set with $k$ distinct elements, where  $k$ is reasonably large. In  case of Design 2b, the set $\cal{S}$ consists of  $k=2^d$ vertices  of the cube $[-\delta,\delta]^d$.
Let $\mathbb{Z}_n=\{Z_1, \ldots, Z_n\}$   denote an $n-$point design whose points $Z_i$ have been sampled
 without replacement  from  $\cal{S}$; $n < k$. Also,
 let $\mathbb{Z}_m^\prime=\{Z_1^\prime, \ldots, Z_m^\prime\}$   denote an associated $m-$point design whose points $Z_i^\prime$ are  sampled
 with replacement  from the same discrete set $\cal{S}$; $Z_1^\prime, \ldots, Z_m^\prime$ are i.i.d. random vectors with values in $\cal{S}$.
 Our aim in this section is to establish an approximate correspondence between $n$ and $m$.

 When sampling $m$ times with replacement, denote by $X_i$ the number of times the $i^{th}$ element of $\cal{S}$ appears.
 Then the vector $(X_1,X_2,\ldots,X_k)$ has the multinomial distribution with number of trials $m$ and event probabilities  $(1/k,1/k, \ldots,1/k) $ with each individual $X_i$ having the Binomial distribution $ Binomial(m,1/k)$. Since ${\rm corr}(X_i,X_j) = -1/k^2$ when $i\neq j$, for large $k$ the correlation between random variables $X_1,X_2,\ldots,X_k$ is very small and will be neglected. Introduce the random variables:
\bea
Y_i = \begin{cases}
1, \text{ if } X_i=0\\
0, \text{ if } X_i> 0.
\end{cases}
\eea
Then the random variable $N_0 = \sum_{i=1}^{k}Y_i$ represents the number of elements of $\cal{S}$ not selected. Given the weak correlation between $X_i$, we approximately have $N_0\sim Binomial(k,P(X_1=0))$. Using the fact $P(X_1=0) = (1-1/k)^{m}$, the expected number of unselected elements when sampling with replacement is approximately 
$\;
\mathbb{E}N_0 \cong k(1-1/k)^{m} \,.
$
Since, when sampling without replacement from $\cal S$ we have chosen $N_0=k-n$ elements, to choose the value of $m$ we equate  $\mathbb{E}N_0 $ to $ k-n$. By solving the equation $$  k-n= k\left(1-\frac1{k}\right)^{m}\;\;\; \left[ \;\cong\; \mathbb{E}N_0 \;\right]$$ for $m$ we obtain
\be\label{value_of_n}
 m = \frac{\log(k-n)-\log(k)}{\log(k-1)-\log(k)} \, .
\ee

\subsection{Approximation of $C_d(\mathbb{Z}_n,r)$ for Design 2b.}

Consider now  $C_d(\mathbb{Z}_n,r)$ for Design 2b. By applying the approximation developed in the previous section, the quantity $C_d(\mathbb{Z}_n,r)$ can be approximated by $C_d(\mathbb{Z}_{m},r)$ for Design 2a with $m$ given in \eqref{value_of_n}: \\

{\bf Approximation of $C_d(\mathbb{Z}_n,r)$ for Design 2b.} {\it  We approximate it   by $C_d(\mathbb{Z}_{m},r)$ where $m$ is given in \eqref{value_of_n} and $C_d(\mathbb{Z}_{m},r)$ is approximated by  \eqref{eq:accurate_app2} with $n$ substituted by $m$ from \eqref{value_of_n}.  }

Specifying this, we obtain:

\be
C_d(\mathbb{Z}_n,r) \simeq 1- \int_{-\infty}^{\infty} \psi_{4,m}(s) \varphi(s)d s, \;\; \label{eq:prod7a}\;\;
\ee
where
\be \label{eq:m_n}
m =m_{n,d}= \frac{\log(2^d-n)-d\log(2)}{\log(2^d-1)-d\log(2)} \,
\ee
and the function  $\psi_{4,\cdot}(\cdot)$  is defined in \eqref{eq:psi4}.
The accuracy of the approximation  \eqref{eq:prod7a}   will be assessed in Section~\ref{sec:numeric1}.

\section{Numerical study}

\label{sec:numeric}

\subsection{Assessing accuracy of approximations of $C_d(\mathbb{Z}_n,r)$   and studying their dependence on $\delta$ }

\label{sec:numeric1}

In this section, we present the results of a large-scale numerical study assessing the accuracy of approximations \eqref{eq:prod5a}, \eqref{eq:accurate_app}, \eqref{normal_alp_0}, \eqref{eq:accurate_app2} and \eqref{eq:prod7a}. In Figures~\ref{pos_alpha1}--\ref{end_pic}, by using a solid black line we depict $C_d(\mathbb{Z}_n,r) $ obtained by Monte Carlo methods, where the value of $r$ has been chosen such that the maximum coverage across $\delta$ is approximately $0.9$.
In Figures~\ref{pos_alpha1}--\ref{approx2}, dealing with Design 1, approximations \eqref{eq:prod5a} and \eqref{eq:accurate_app}  are depicted with a dotted blue and
dashed green lines respectively.
In Figures~\ref{approx1}--\ref{alpha_0_fig} (Design 2a)
approximations \eqref{normal_alp_0}   and \eqref{eq:accurate_app2} are illustrated with a dotted blue and   dashed green lines respectively.
In Figures~\ref{without_pic1}--\ref{end_pic} (Design 2b) the dashed green line depicts approximation \eqref{eq:prod7a}.
From these figures, we can draw the following conclusions.
\begin{itemize}
\item Approximations \eqref{eq:accurate_app} and \eqref{eq:accurate_app2} are very accurate across all values of $\delta$ and $\alpha$. This is particularly evident for $d=20,50$.
\item Approximations \eqref{eq:prod5a} and \eqref{normal_alp_0} are  accurate only for very large values of $d$, like $d=50$.
\item Approximation \eqref{eq:prod6a} is generally accurate. For  $\delta$ close to one (for such values of $\delta$ the covering is very poor) and $n$ close to $2^d$ this approximation begins to worsen, see Figures~\ref{without_pic2} and  \ref{end_pic}.
\item A sensible choice of $\delta$ can dramatically increase the coverage proportion $C_d(\mathbb{Z}_n,r) $. This effect, which we call `$\delta$-effect',  is evident in all figures and is very important. It gets much stronger as $d$ increases.
\end{itemize}

\begin{figure}[h]
\centering
\begin{minipage}{.5\textwidth}
  \centering
  \includegraphics[width=1\linewidth]{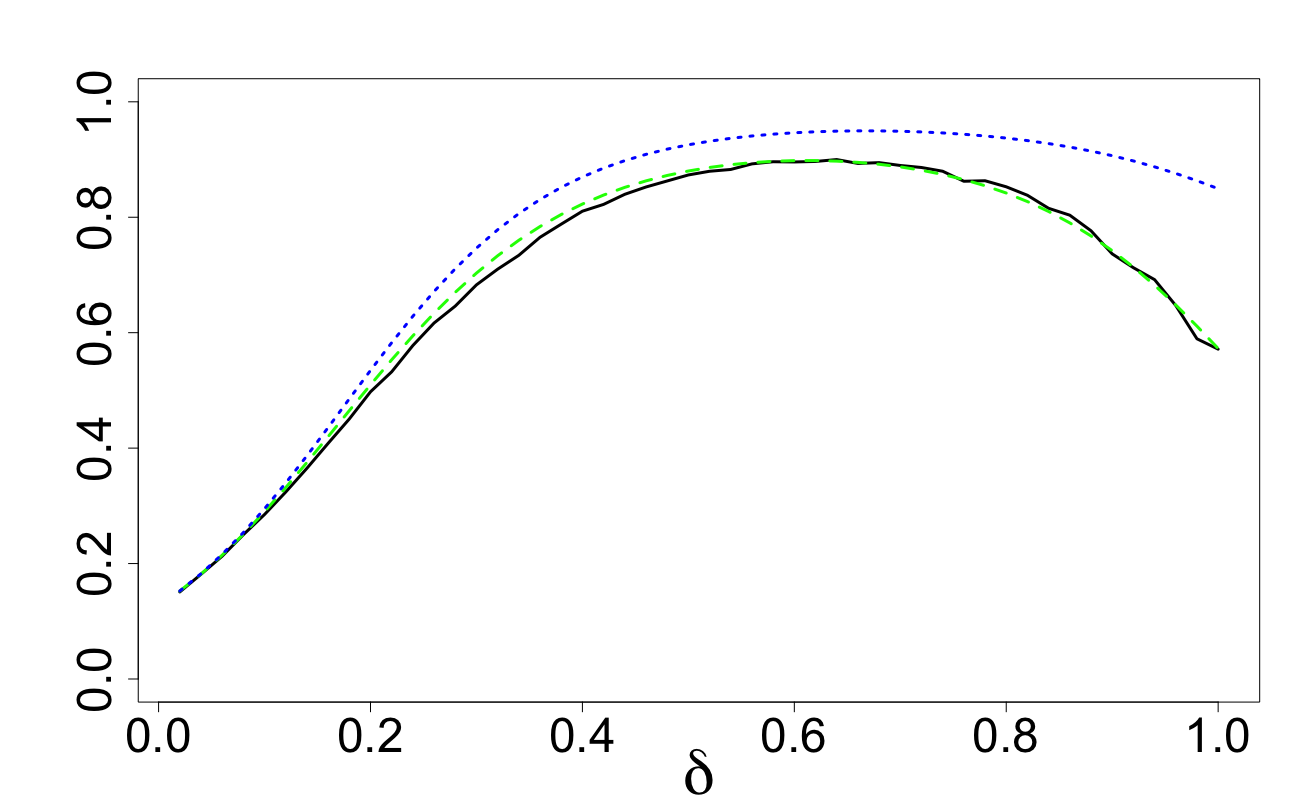}
  \caption{Design 1: $C_d(\mathbb{Z}_n,r) $ and approximations;\\ $d=10, \alpha=0.5, n=128$.}
  \label{pos_alpha1}
\end{minipage}%
\begin{minipage}{.5\textwidth}
  \centering
  \includegraphics[width=1\linewidth]{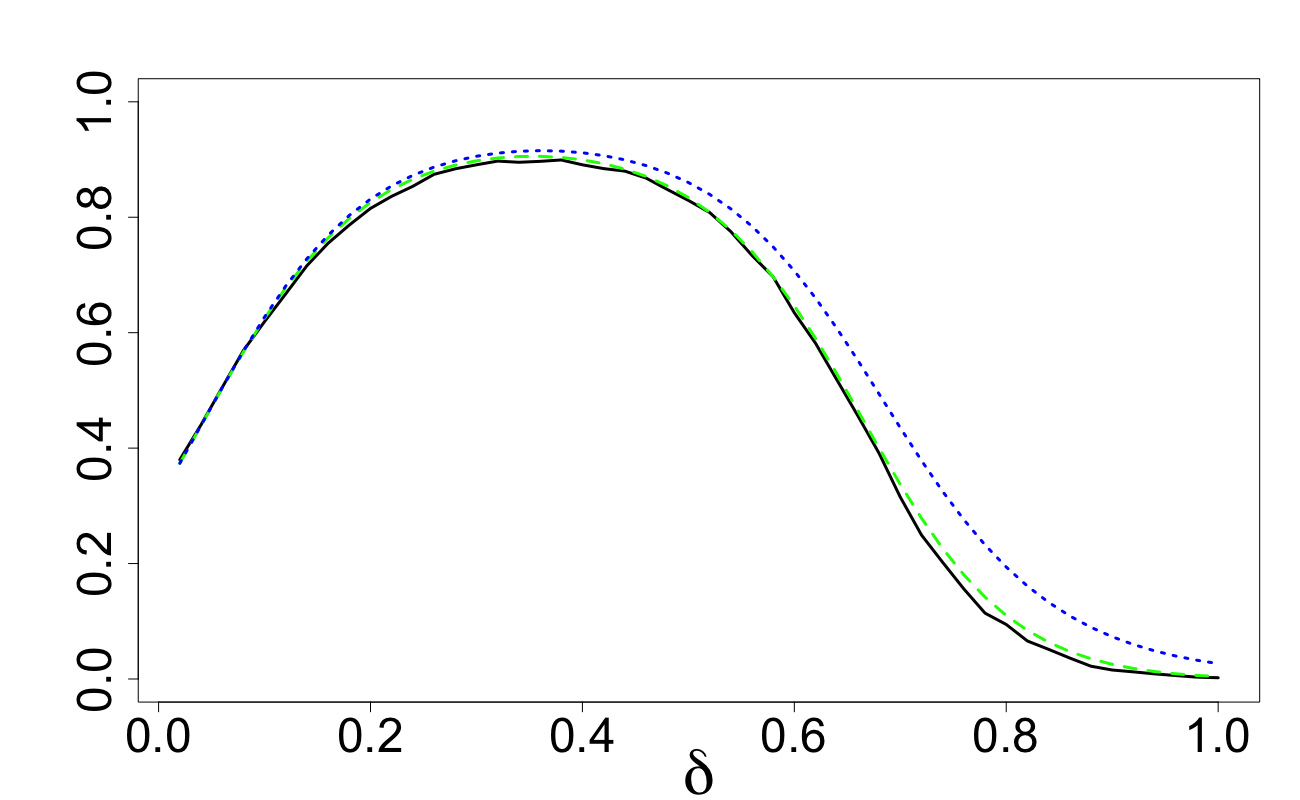}
  \caption{Design 1: $C_d(\mathbb{Z}_n,r) $ and approximations;\\ $d=20, \alpha=0.1, n=128$.}
\end{minipage}
\end{figure}

\begin{figure}[h]
\centering
\begin{minipage}{.5\textwidth}
  \centering
  \includegraphics[width=1\linewidth]{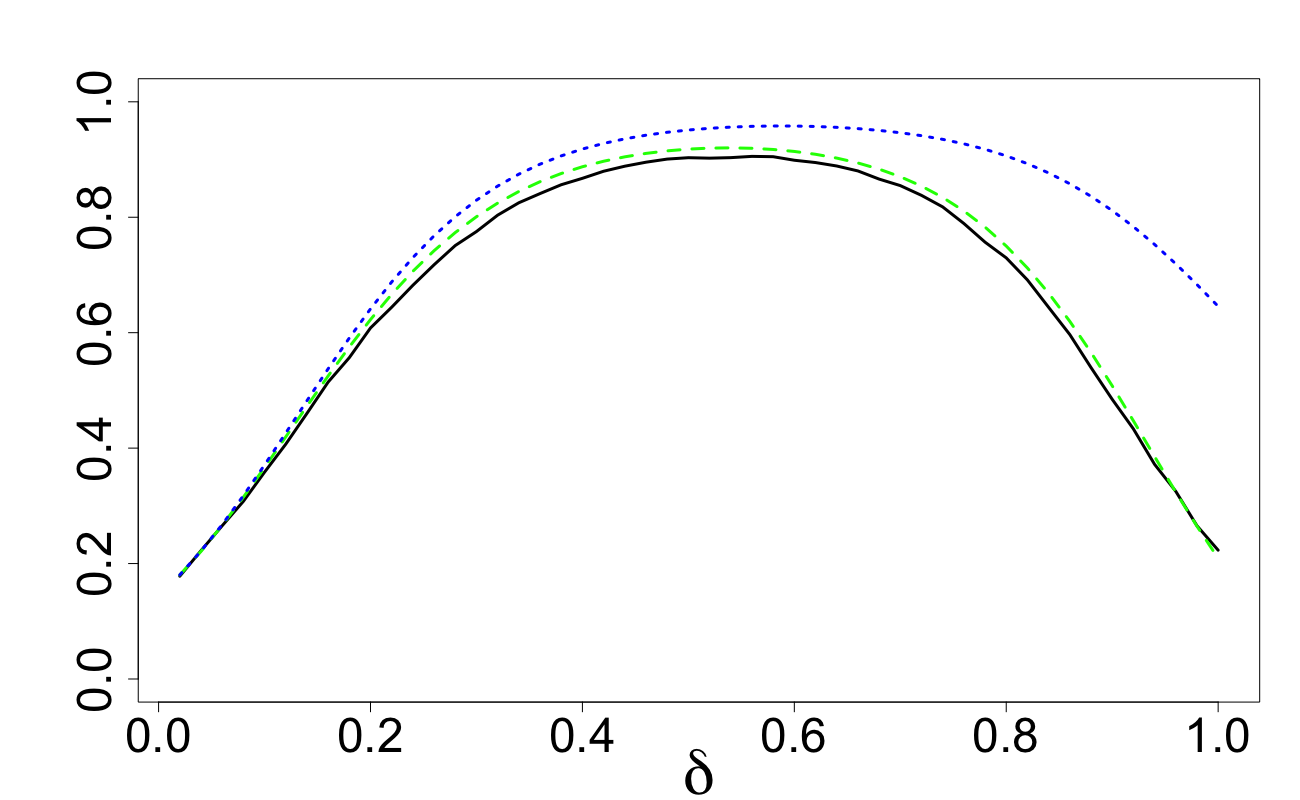}
  \caption{Design 1: $C_d(\mathbb{Z}_n,r) $ and approximations;\\ $d=20, \alpha=0.5, n=512$.}
\end{minipage}%
\begin{minipage}{.5\textwidth}
  \centering
  \includegraphics[width=1\linewidth]{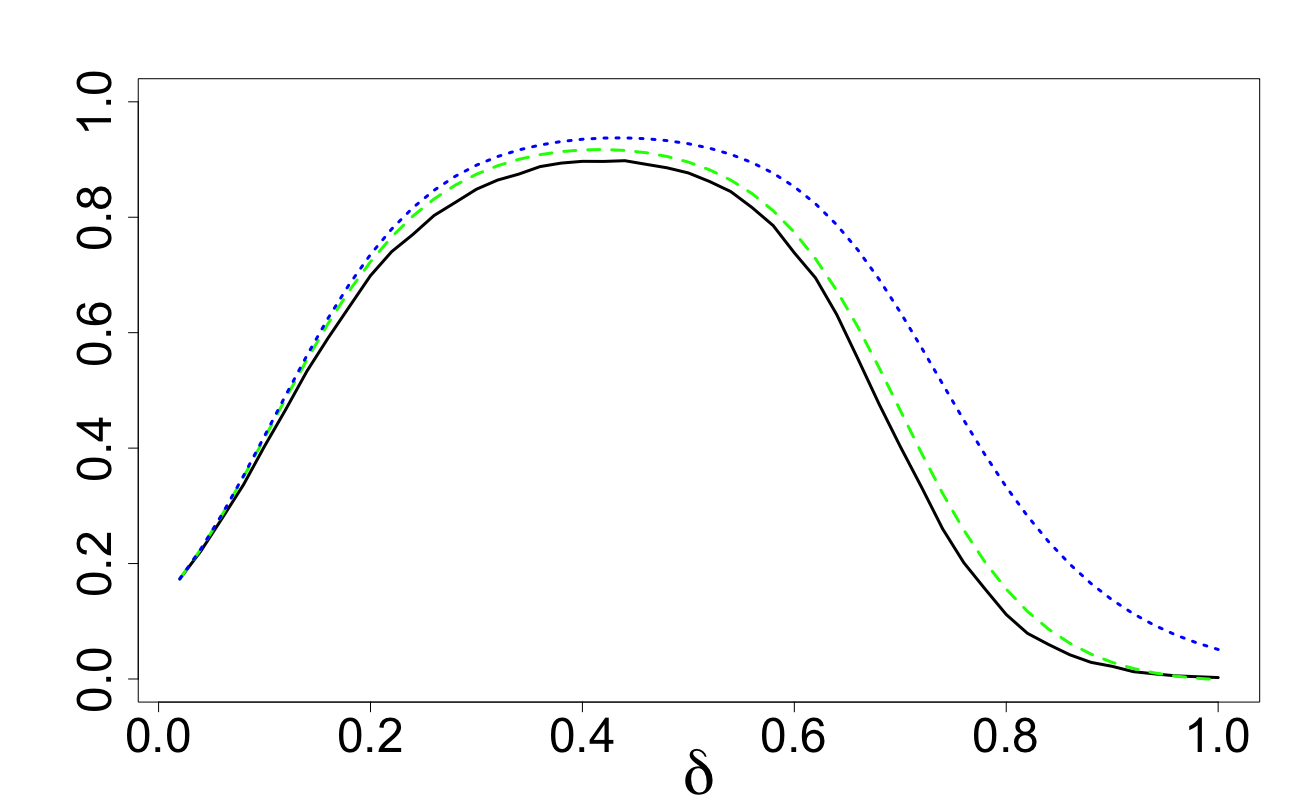}
  \caption{Design 1: $C_d(\mathbb{Z}_n,r) $ and approximations;\\ $d=20, \alpha=0.1, n=512$.}
\end{minipage}
\end{figure}

\begin{figure}[h]
\centering
\begin{minipage}{.5\textwidth}
  \centering
  \includegraphics[width=1\linewidth]{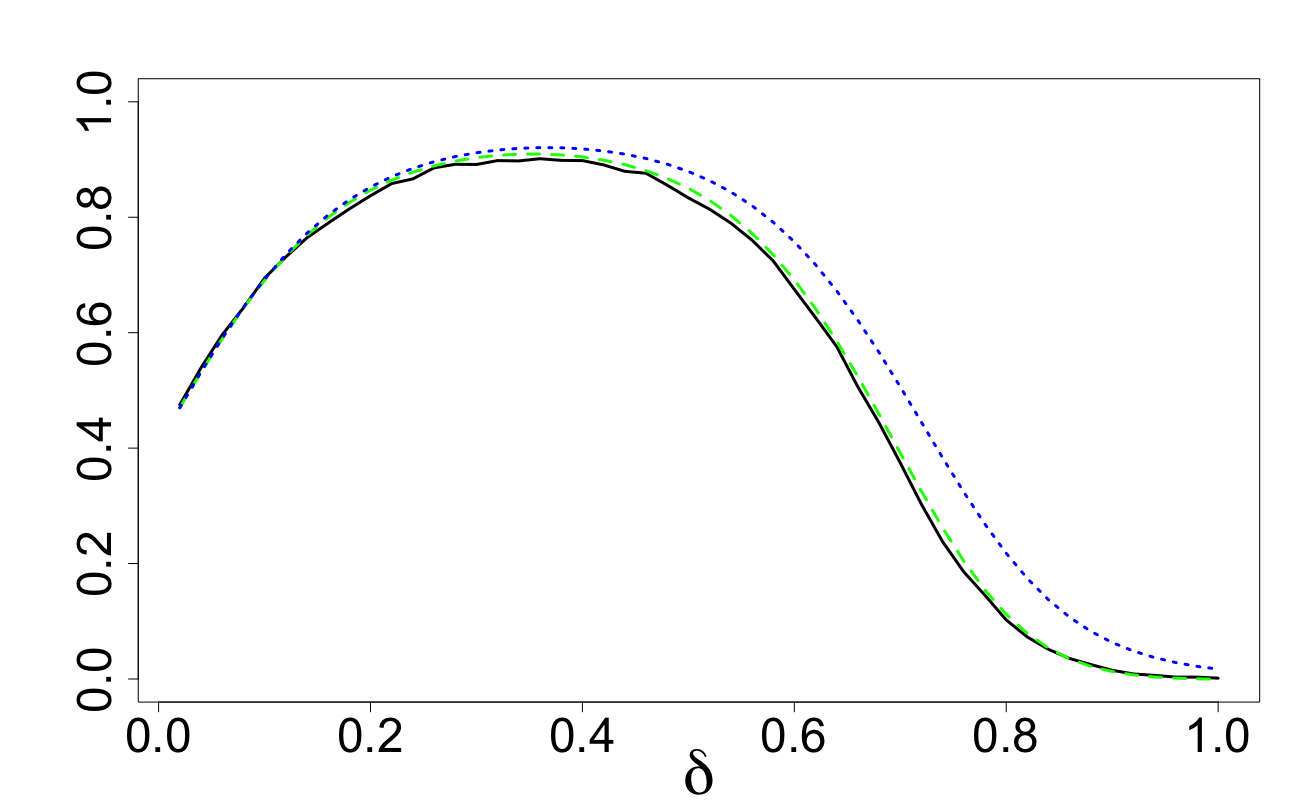}
  \caption{Design 1: $C_d(\mathbb{Z}_n,r) $ and approximations;\\ $d=50, \alpha=0.5, n=512$.}
\end{minipage}%
\begin{minipage}{.5\textwidth}
  \centering
  \includegraphics[width=1\linewidth]{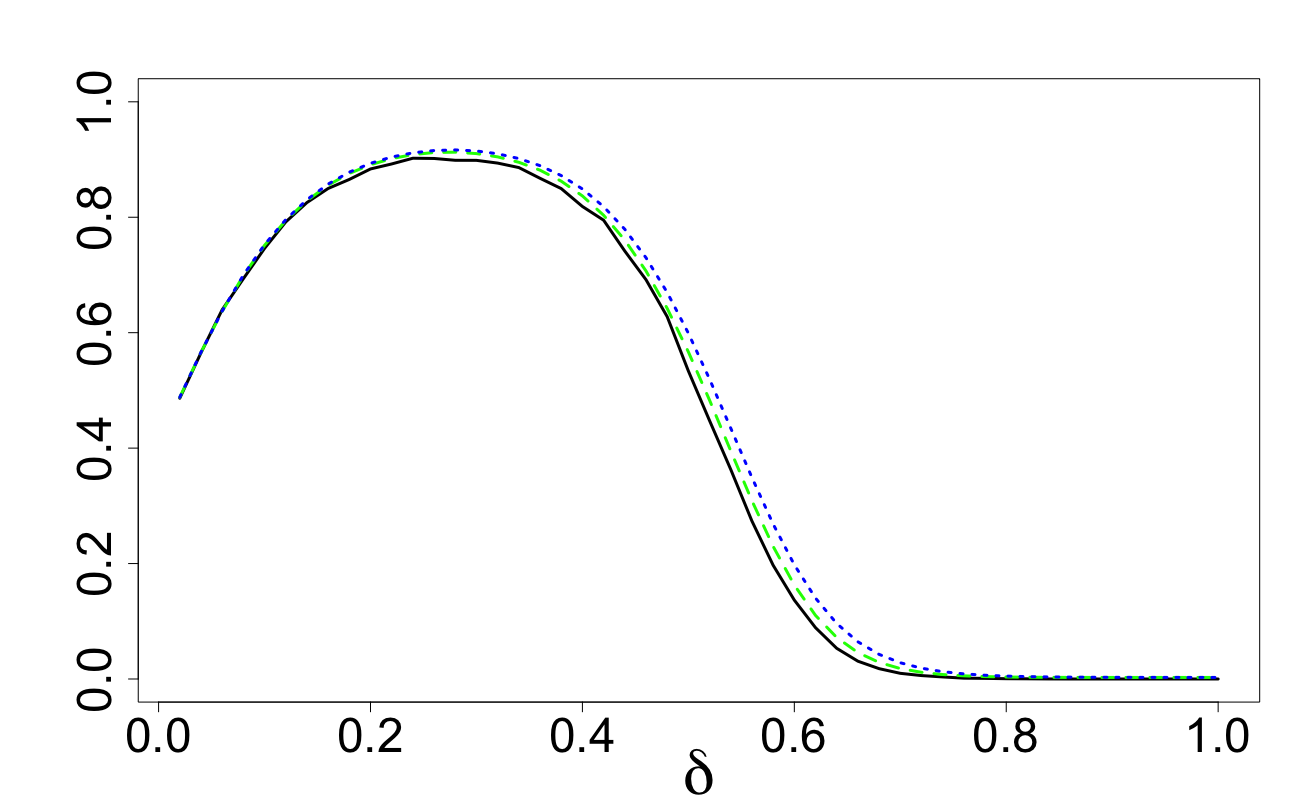}
  \caption{Design 1: $C_d(\mathbb{Z}_n,r) $ and approximations;\\ $d=50, \alpha=0.1, n=512$.}
\label{approx2}
\end{minipage}

\end{figure}

\begin{figure}[h]
\centering
\begin{minipage}{.5\textwidth}
  \centering
  \includegraphics[width=1\linewidth]{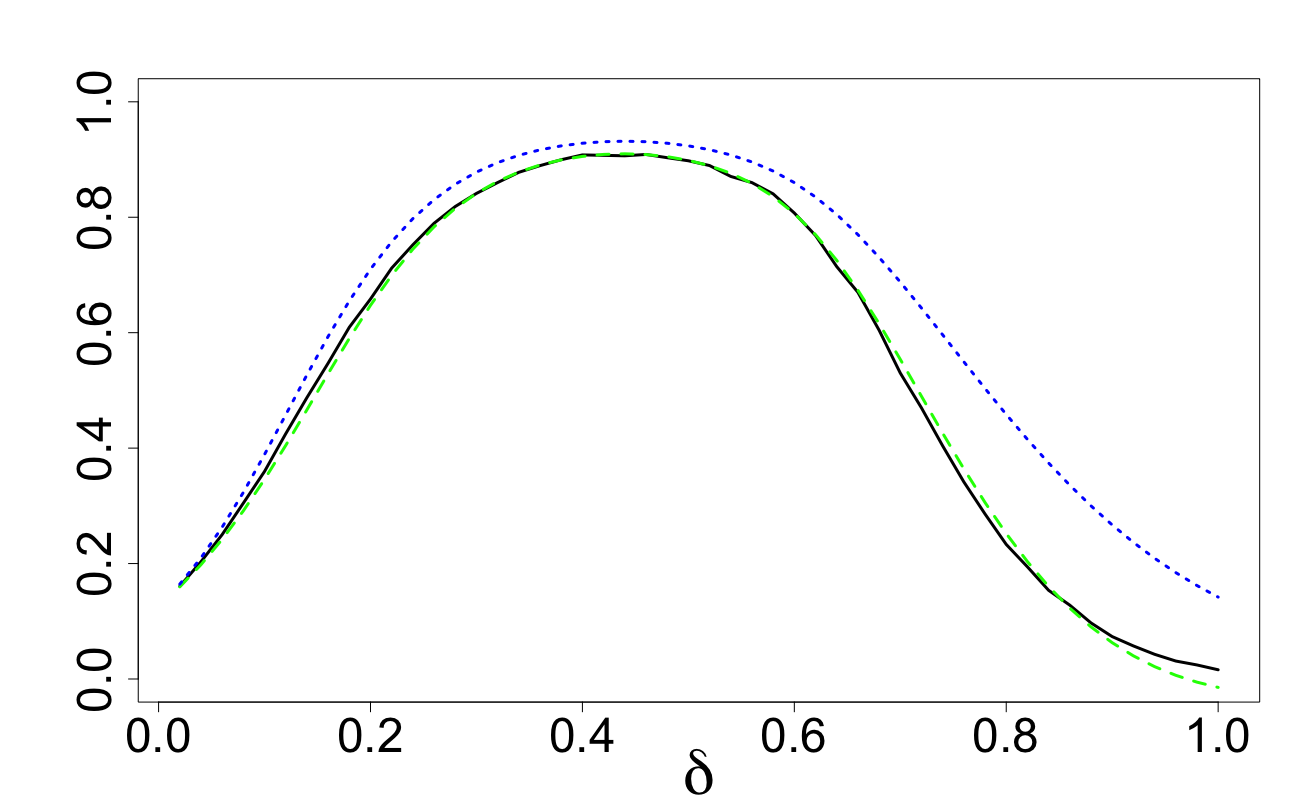}
  \caption{Design 2a: $C_d(\mathbb{Z}_n,r) $ and approximations;\\ $d=10, \alpha=0, n=128$.}
  \label{approx1}
\end{minipage}%
\begin{minipage}{.5\textwidth}
  \centering
  \includegraphics[width=1\linewidth]{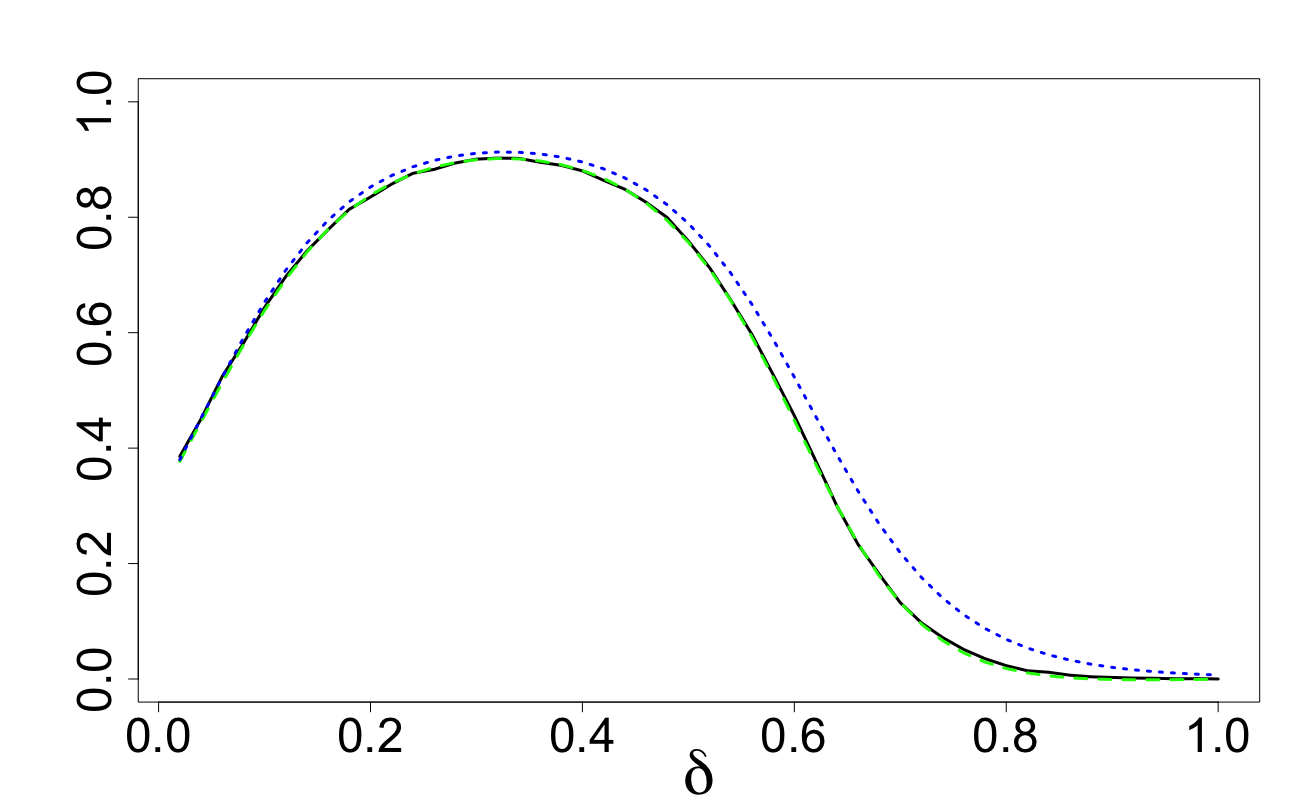}
  \caption{Design 2a: $C_d(\mathbb{Z}_n,r) $ and approximations; $d=20, \alpha=0, n=128$.}
\end{minipage}
\end{figure}

\begin{figure}[h]
\centering
\begin{minipage}{.5\textwidth}
  \centering
  \includegraphics[width=1\linewidth]{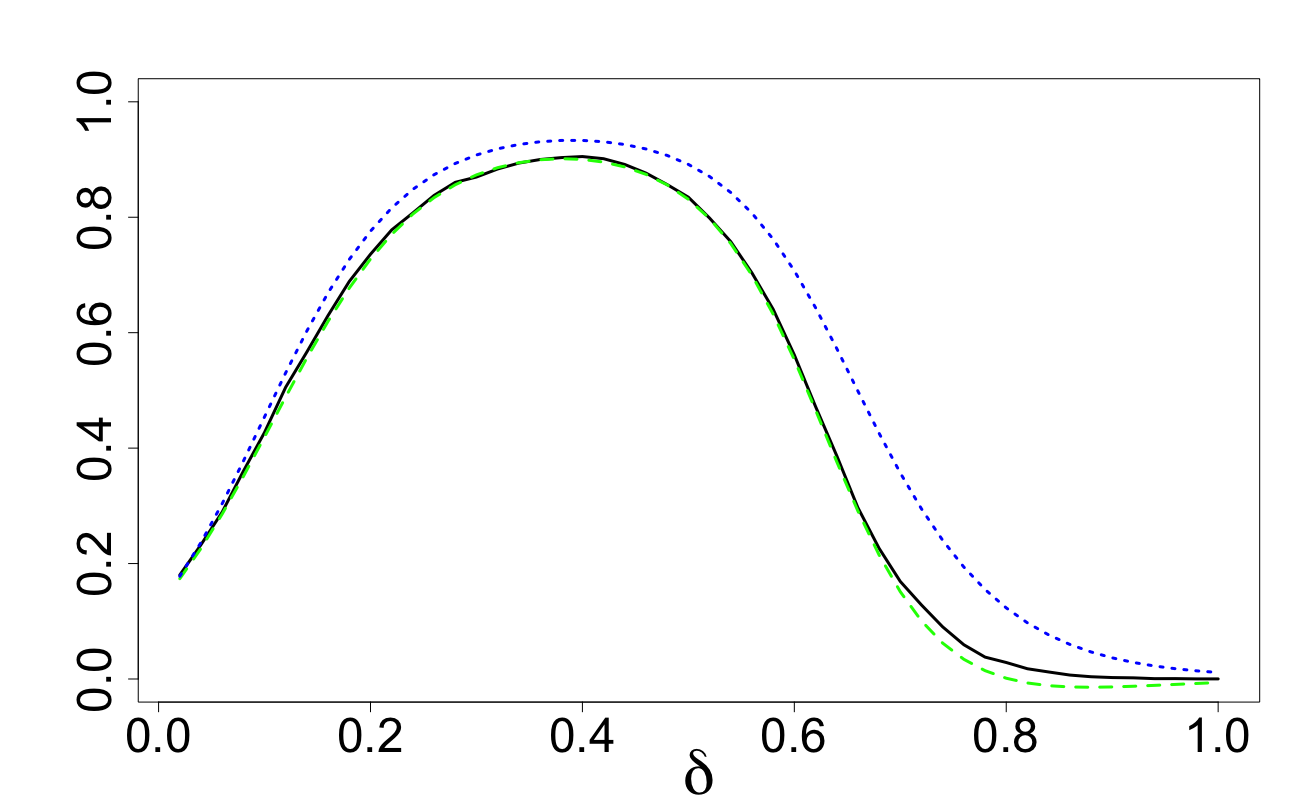}
  \caption{Design 2a: $C_d(\mathbb{Z}_n,r) $ and approximations;\\ $d=20, \alpha=0, n=512$.}
\end{minipage}%
\begin{minipage}{.5\textwidth}
  \centering
  \includegraphics[width=1\linewidth]{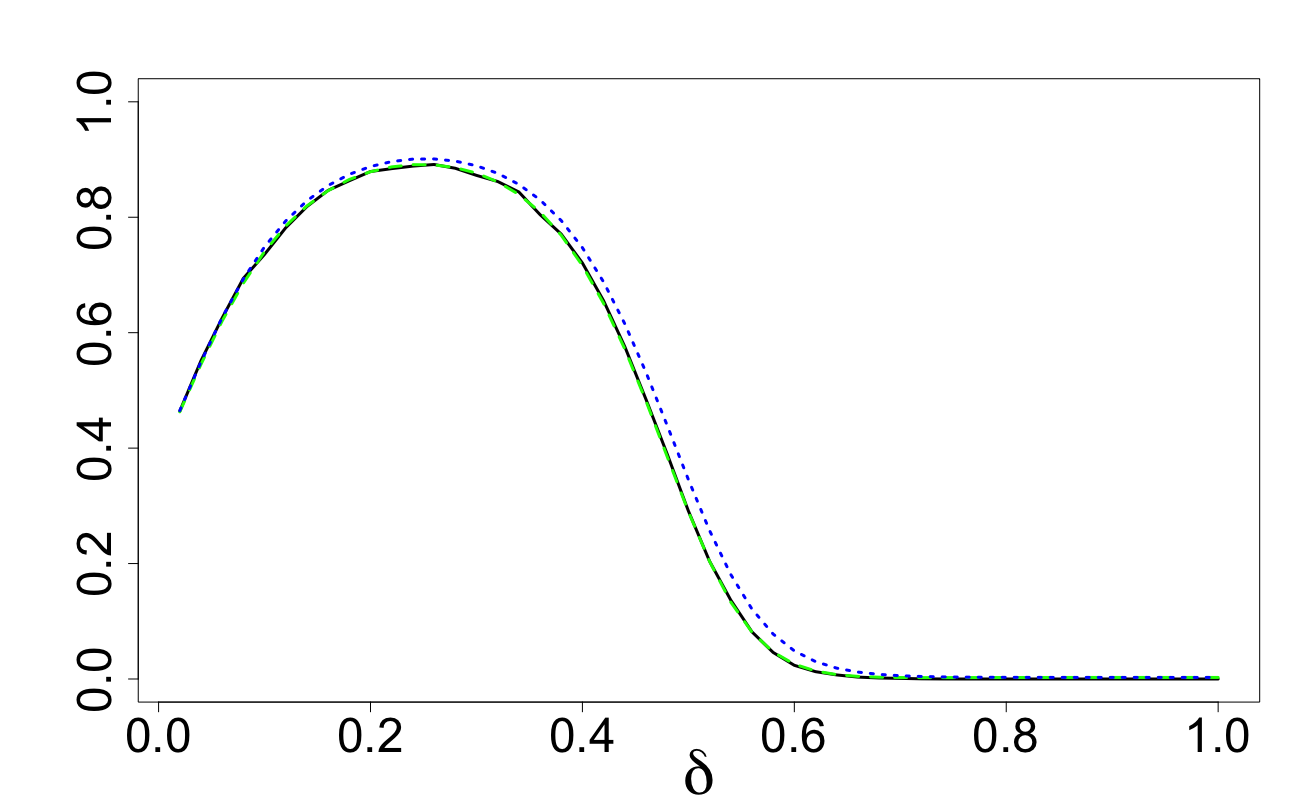}
  \caption{Design 2a: $C_d(\mathbb{Z}_n,r) $ and approximations;\\ $d=50, \alpha=0, n=512$.}
  \label{alpha_0_fig}
\end{minipage}
\end{figure}

\begin{figure}[h]
\centering
\begin{minipage}{.5\textwidth}
  \centering
  \includegraphics[width=1\linewidth]{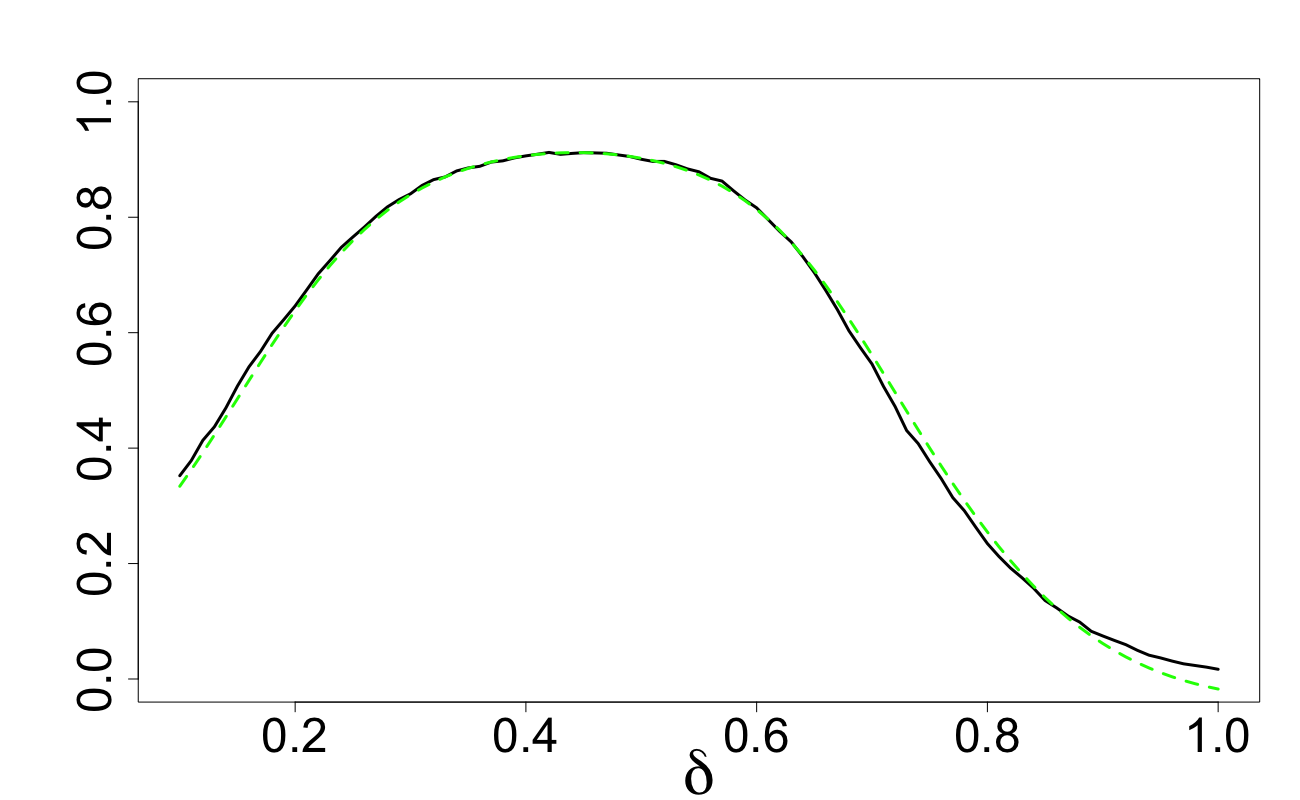}
  \caption{Design 2b: $C_d(\mathbb{Z}_n,r) $ and approxima-\\tion \eqref{eq:prod7a}; $d=10, n=128$.}
    \label{without_pic1}
\end{minipage}%
\begin{minipage}{.5\textwidth}
  \centering
  \includegraphics[width=1\linewidth]{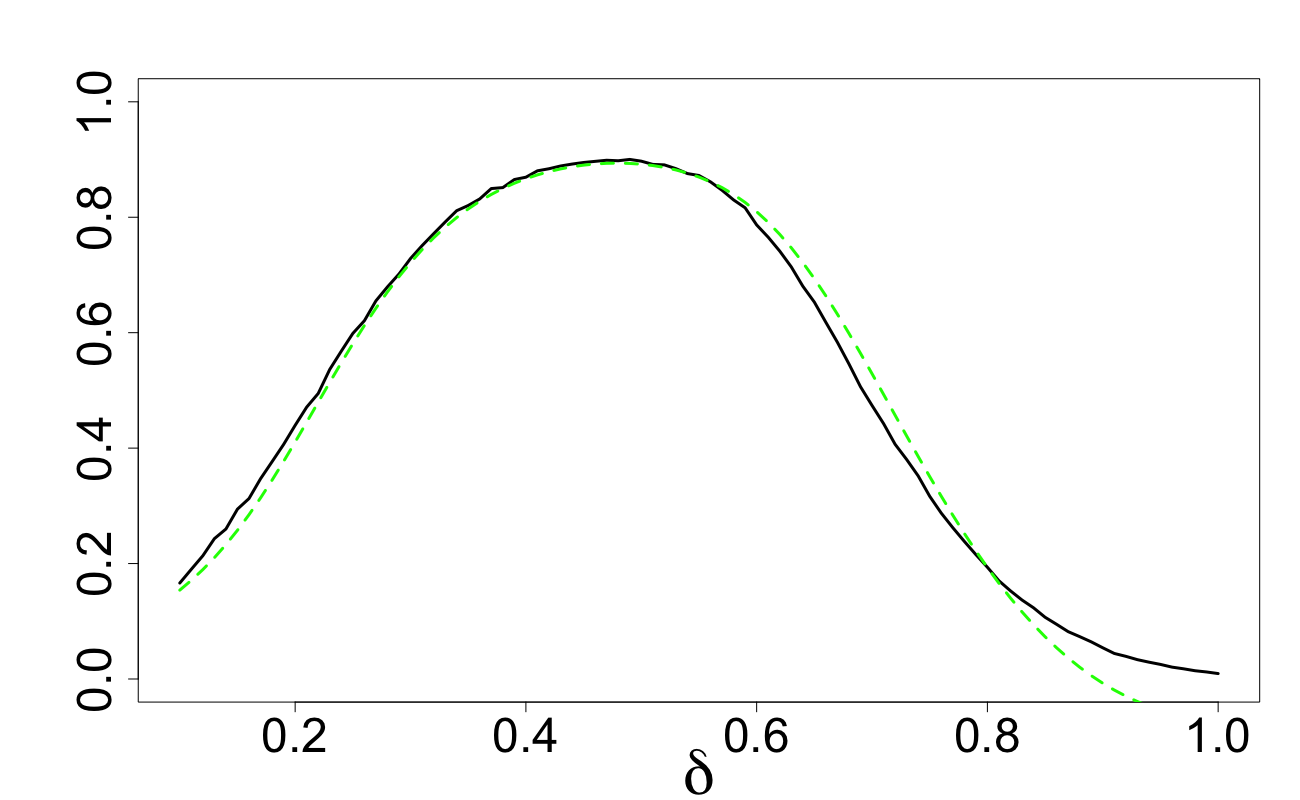}
  \caption{Design 2b: $C_d(\mathbb{Z}_n,r) $ and approximation\\ \eqref{eq:prod7a}; $d=10, n=256$.}
    \label{without_pic2}
\end{minipage}

\end{figure}

\begin{figure}[h]
\centering
\begin{minipage}{.5\textwidth}
  \centering
  \includegraphics[width=1\linewidth]{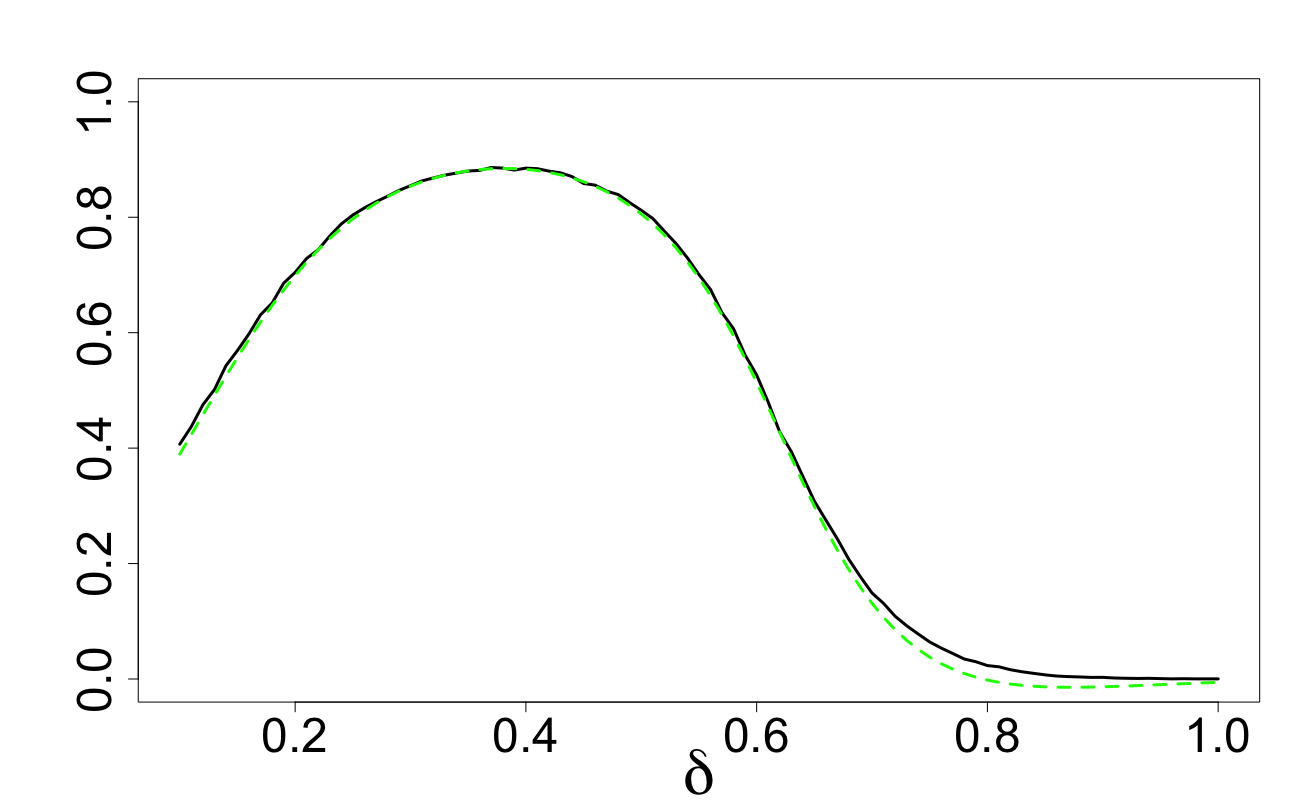}
  \caption{Design 2b: $C_d(\mathbb{Z}_n,r) $ and approxima-\\tion \eqref{eq:prod7a}; $d=20, n=512$.}
\end{minipage}%
\begin{minipage}{.5\textwidth}
  \centering
  \includegraphics[width=1\linewidth]{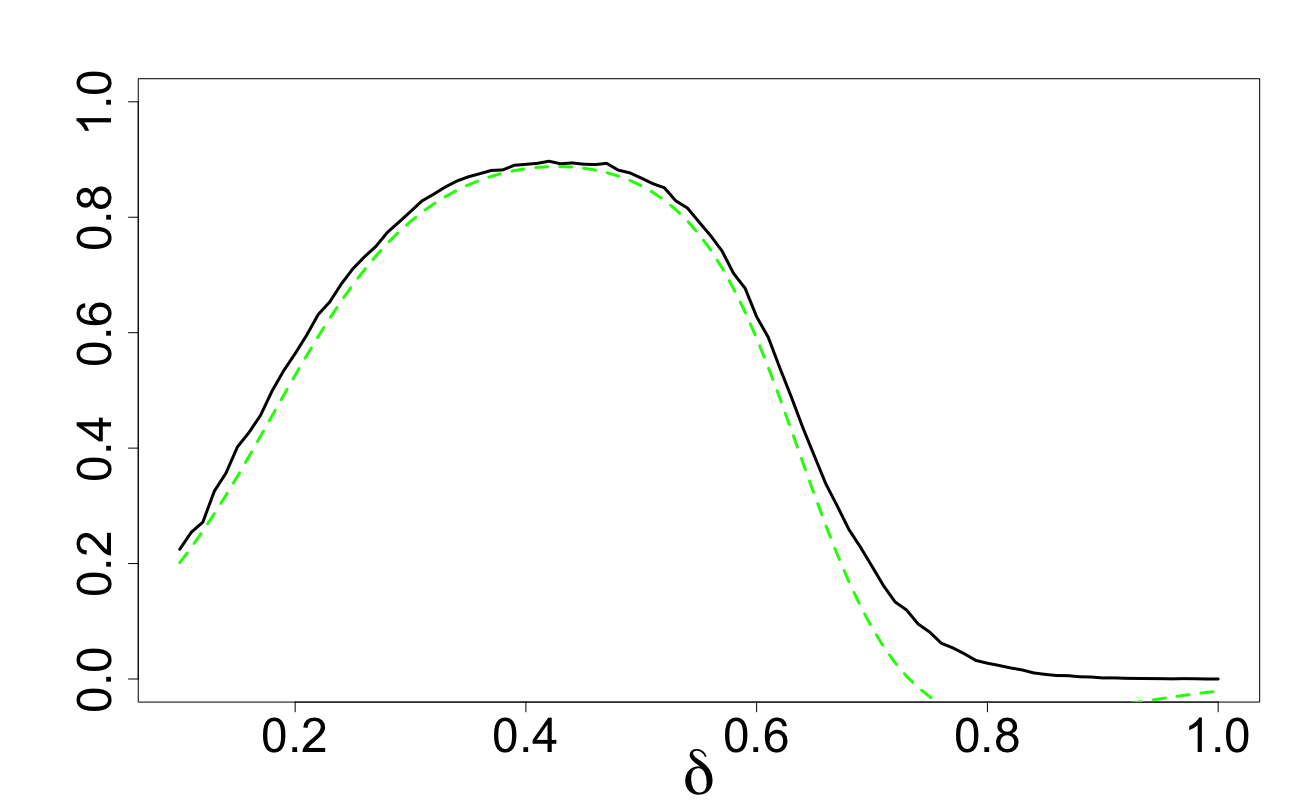}
  \caption{Design 2b: $C_d(\mathbb{Z}_n,r) $ and approximation\\ \eqref{eq:prod7a}; $d=20, n=2048$.}
  \label{end_pic}
\end{minipage}

\end{figure}

\subsection{Comparison across $\alpha$}

In Table~\ref{alpha_table}, for Design~2a and Design 1 with $\alpha=0.5,1,1.5$ we present the smallest values of $r$ required to achieve the 0.9-coverage on average. For these schemes, the value inside the brackets shows the average value of $\delta$ required to obtain this 0.9-coverage. Design 2b  is not used as $d$ is too small (for this design, we must have $n < 2^d$ and in these cases Design 2b provides better coverings than the other designs considered).

\begin{table}[h]
\centering
\begin{tabular}{ |p{2.5cm}||p{2cm}|p{2cm}|p{2cm}|p{2cm}|   }
 \hline
 \multicolumn{5}{|c|}{$d=5$} \\
 \hline
  & $n=25 $ &$n=50$& $n=100$ & $n=500$\\
 \hline
 Design 2a $(\alpha=0)$ & 1.051 (0.44)    &  0.885 (0.50)  & 0.812 (0.50)    & 0.798 (0.50)      \\
Design 1, $\alpha=0.5$& 1.072 (0.68)    & 0.905 (0.78)  & 0.770 (0.78)   &   0.540 (0.80)   \\
Design 1, $\alpha=1$&  1.072 (0.78)    & 0.931 (0.86)   &  0.798 (0.98) &  0.555  (1.00)   \\
Design 1, $\alpha=1.5$& 1.091 (0.92)     & 0.950 (0.96)   &0.820 (0.98)   & 0.589 (1.00)    \\
 \hline
\end{tabular}
\caption{Values of $r$ and $\delta$ (in brackets) to achieve 0.9 coverage for $d=5$.}
\label{alpha_table}
\end{table}
\begin{table}[h]
\centering
\begin{tabular}{ |p{2.5cm}||p{2cm}|p{2cm}|p{2cm}|p{2cm}|   }
 \hline
 \multicolumn{5}{|c|}{$d=10$} \\
 \hline
  & $n=500 $ &$n=1000$& $n=5000$ & $n=10000$\\
 \hline
 Design 2a $(\alpha=0)$ &  1.228 (0.50)    &  1.135 (0.50)   & 1.073   (0.50)    &  1.071 (0.50)       \\
Design 1, $\alpha=0.5$&  1.271 (0.69)     & 1.165 (0.73)  &  0.954 (0.76)   &  0.886 (0.78)     \\
Design 1, $\alpha=1$&   1.297 (0.87)   & 1.194 (0.90)    & 0.992 (0.93)   & 0.917 (0.95)     \\
Design 1, $\alpha=1.5$&  1.320 (1.00)   &  1.220  (1.00)   & 1.032 (1.00)   & 0.953 (1.00)    \\
 \hline
\end{tabular}
\caption{Values of $r$ and $\delta$ (in brackets) to achieve 0.9 coverage for $d=10$.}
\label{alpha_table2}
\end{table}

From Tables~\ref{alpha_table} and \ref{alpha_table2} we can make the following conclusions:
\begin{itemize}
\item For small $n$ ($n<2^d$ or $n\simeq 2^d$), Design 2a provides a more efficient covering than other three other schemes and hence   smaller values of $\alpha$ are better.
\item For $n>2^d$, Design 2a begins to become impractical since a large proportion of points  duplicate. This is reflected in Table~\ref{alpha_table} by comparing $n=100$ and $n=500$ for Design 2a; there is only a small reduction in $r$ despite a large increase in $n$. Moreover, for values of $n>>2^d$,  Design 2a provides a very inefficient covering.
\item For $n>>2^d$, from looking at Design 1 with $\alpha=0.5$ and $n=500$, it would appear beneficial to choose  $\alpha \in (0,1)$  rather than $\alpha>1$ or $\alpha=0$.
\end{itemize}

Using approximations  \eqref{eq:accurate_app} and  \eqref{eq:accurate_app2}, in Figures~\ref{alpha_comparison1}--\ref{alpha_comparison2} we depict $C_d(\mathbb{Z}_n,r)$ across $\delta$ for different choices of~$\alpha$. In Figures~\ref{alpha_comparison1}--\ref{alpha_comparison2}, the red line, green line, blue line and cyan line depict approximation \eqref{eq:accurate_app2} ($\alpha=0$) and approximation \eqref{eq:accurate_app} with $\alpha=0.5$, $\alpha=1$ and $\alpha=1.5$ respectively.  These figures demonstrate the clear benefit of choosing a smaller $\alpha$, at least for these values of $n$ and $d$.

\begin{figure}[h]
\centering
\begin{minipage}{.5\textwidth}
  \centering
  \includegraphics[width=1\linewidth]{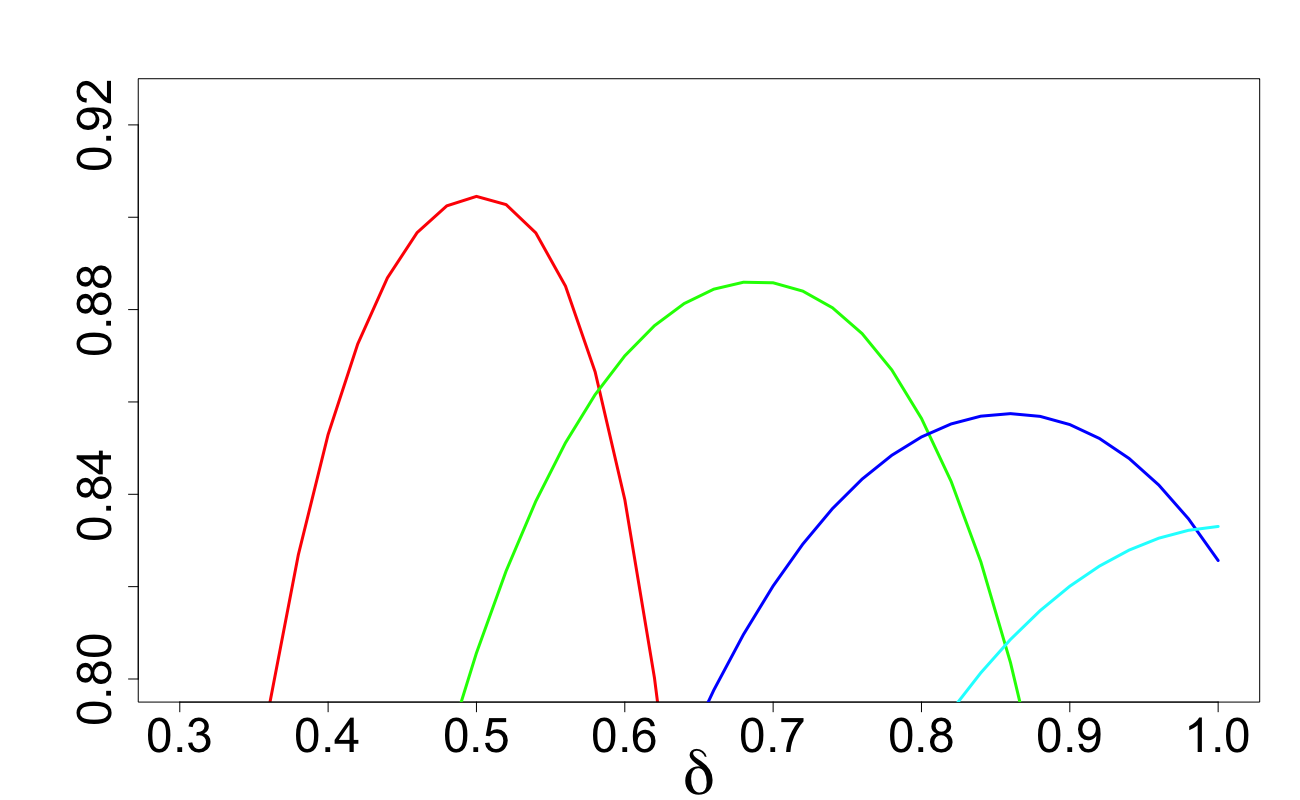}
  \caption{$d=10$, $n=512$, $r=1.228$ }
  \label{alpha_comparison1}
\end{minipage}%
\begin{minipage}{.5\textwidth}
  \centering
  \includegraphics[width=1\linewidth]{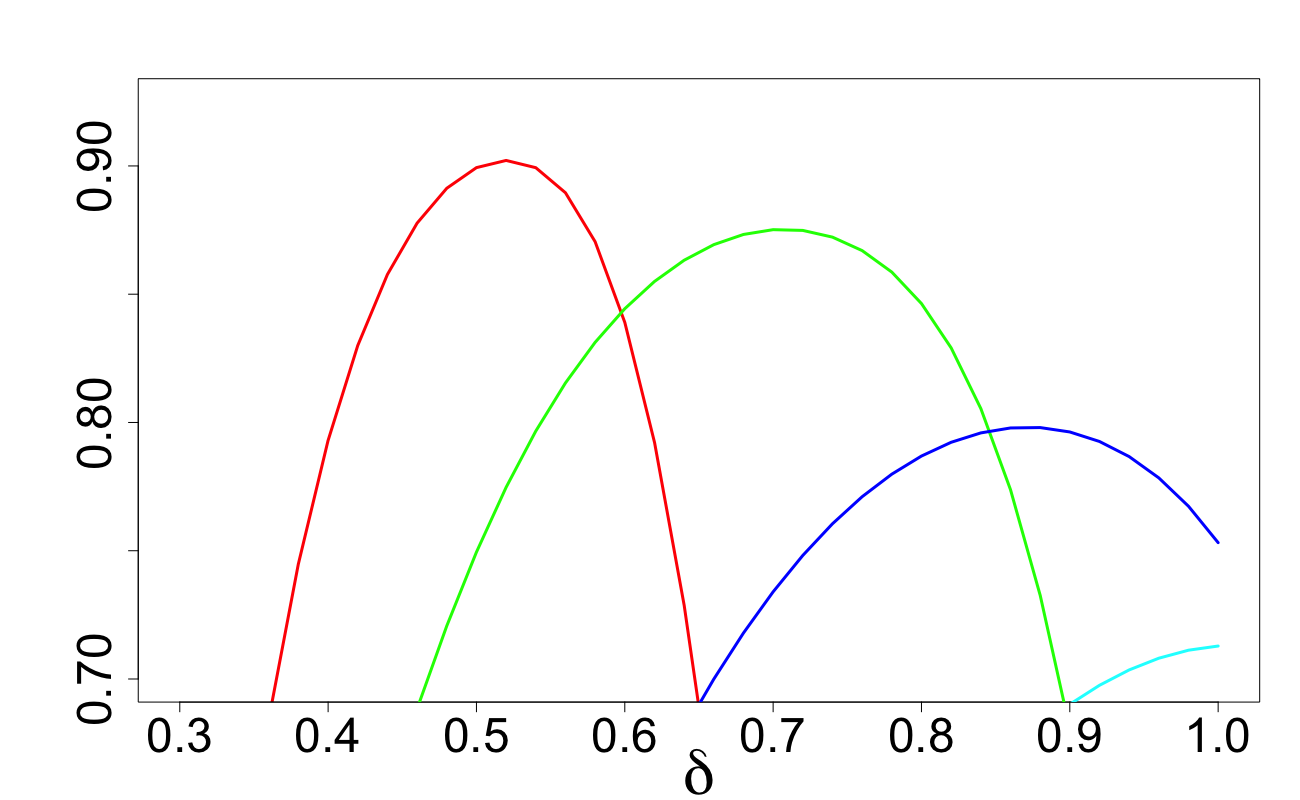}
  \caption{$d=10$, $n=1024$, $r=1.13$ }
  \label{alpha_comparison2}
\end{minipage}

\end{figure}

\section{Quantization in a cube }
\label{sec:quantization}

\subsection{Quantization error and its relation to weak covering}

In this section, we will study the following characteristic of a design $\mathbb{Z}_n $.\\

 \begin{bf}{Quantization error.}\end{bf} Let  $U=(u_1, \ldots, u_d)$ be uniform random vector on $[-1,1]^d$.
The mean squared quantization error for a design $\mathbb{Z}_n=\{Z_1, \ldots, Z_n\} \subset \mathbb{R}^d$ is defined by
\be
\label{eq:rho_s}
\theta(\mathbb{Z}_n)=\mathbb{E}_U \varrho^2(U,\mathbb{Z}_n)\,, \;\;{\rm where}\;\;
\varrho(U,\mathbb{Z}_n)= \min_{Z_i \in \mathbb{Z}_n} \|U-Z_i\|\, .
\ee
If the design $\mathbb{Z}_n$ is randomized then we consider the expected value $\mathbb{E}_{\mathbb{Z}_n} \theta(\mathbb{Z}_n)$
of $\theta(\mathbb{Z}_n)$ as the main characteristic without stressing this.

 The mean squared quantization error $\theta(\mathbb{Z}_n)$ is related to our main quantity $C_d(\mathbb{Z}_n,r)$ defined in \eqref{eq:cover1}: indeed, $C_d(\mathbb{Z}_n,r)$, as a function of $r\geq 0$, is the c.d.f. of the r.v. $\varrho(U,\mathbb{Z}_n)$ while $\theta(\mathbb{Z}_n)$ is the second moment of the distribution with this c.d.f.:
\be\label{eq:relation}
\theta(\mathbb{Z}_n)= \int_{r\geq 0} r^2 d  C_d(\mathbb{Z}_n,r)\, .
\ee
This relation will allow us to use the approximations derived above for $C_d(\mathbb{Z}_n,r)$ in order to construct approximations for the quantization error $\theta(\mathbb{Z}_n)$.

%

\subsection{Quantization error for Design 1 }

Using approximation \eqref{eq:accurate_app} for the quantity $C_d(\mathbb{Z}_n,r)$, we obtain
\be\label{quant_density}
&&\frac{d}{dr}(C_d(\mathbb{Z}_n,r)) \cong f_{\alpha,\delta}(r):= \frac{n \cdot r}{\delta} \int_{-\infty}^{\infty} \varphi(s)\varphi(c_s) \psi_{2,\alpha}(s) \times \nonumber \\
&&\times\left [\frac{\sqrt{2\alpha+1}}{\sqrt{s' + k}}+c_{d,\alpha} \frac{\alpha  \left( s'+\frac{d\delta^2(2\alpha-1)}{3(2\alpha+5)(2\alpha+1)} \right)}{ (2\alpha+3)\left (s'+k\right)^{2} } \left\{{\delta(c_s^3-c_s)} - \frac{\sqrt{2\alpha+1}(r^2-\frac{d\delta^2}{2\alpha+1}-s')}{ \sqrt{s'+k}}   \right\}  \right] ds \,.
\ee
By then using relation \eqref{eq:relation} we obtain the following  approximation for the mean squared quantization error with Design 1:
\be\label{final_quant}
\theta(\mathbb{Z}_n)\cong  \int_{ 0}^{\infty} r^2 f_{\alpha,\delta}(r) dr \, .
\ee

By taking $\alpha=1$ in \eqref{quant_density} we obtain:
\bea
 f_{1,\delta}(r) \!\!\!&:=&\!\! \frac{n \cdot r}{\delta} \int_{-\infty}^{\infty} \varphi(s)\varphi(c_s) \psi_{2,1}(s) \! \left [\frac{\sqrt{3}}{\sqrt{s' + k}}+c_{d,1} \frac{  \left( s'+\frac{d\delta^2}{63} \right)}{ 5\left (s'+k\right)^{2} } \left\{{\delta(c_s^3-c_s)} - \frac{\sqrt{3}(r^2-\frac{d\delta^2}{3}-s')}{ \sqrt{s'+k}}   \right\}  \right] ds \,.
\eea
with $\psi_{2,1}$ defined in \eqref{eq:psi2-1}.
The resulting approximation
\bea
\theta(\mathbb{Z}_n)\cong  \int_{ 0}^{\infty} r^2 f_{1,\delta}(r) dr \, .
\eea
coincides with \cite[formula 31]{us}.

\subsection{Quantization error for Design 2a }

Using approximation \eqref{eq:accurate_app2} for the quantity $C_d(\mathbb{Z}_n,r)$, we have:

\be\frac{d}{dr}(C_d(\mathbb{Z}_n,r)) \cong f_{0,\delta;n}(r):=\;\;\;\;\;\;\;\;\;\;\;\;\;
\;\;\;\;\;\;\;\;\;\;\;\;\;\;\;\;\;\;\;\;\;\;\;\;\;\;\;\;\;\;\;\;\;\;\;\;\;\;\;\;\;\;\;\;\;
\;\;\;\;\;\;\;\;\;\;\;\;\;\;\;\;\;\;\;\;\;\;\;\;\;\;\;\;\;\;\nonumber \\
 \frac{n \cdot r}{\delta} \int_{-\infty}^{\infty} \frac{\varphi(s)\varphi(c_s) \psi_{4,n}(s)}{\sqrt{s'}}
\left [1+\left(1+\frac3d \right)\frac{(2(s'-d/3)/\sqrt{5}+d/5)(6c_s^2-c_s^4-3)}{12(s')^2}   \right] ds \, ,\label{eq:des2a_quant}
\ee

where $\psi_{4,n}(\cdot)$ is defined in \eqref{eq:psi4}.
From \eqref{eq:relation} we then obtain the following  approximation for  the mean squared quantization error with Design 2a:
\be\label{final_quant_alpha_zero}
\theta(\mathbb{Z}_n)\cong  \int_{ 0}^{\infty} r^2 f_{0,\delta;n}(r) dr \, .
\ee

\subsection{Quantization error for Design 2b }

Similarly to \eqref{final_quant_alpha_zero}, for Design 2b, we use the approximation
\be\label{final_quant_alpha_zero2b}
\theta(\mathbb{Z}_n)\cong  \int_{ 0}^{\infty} r^2 f_{0,\delta;m}(r) dr \, .
\ee
 where $f_{0,\delta,m}(r)$ is defined by
\eqref{eq:des2a_quant}
and $m=m_{n,d}$ is defined in \eqref{eq:m_n}.

\subsection{Accuracy of  approximations for quantization error and the $\delta$-effect}

\label{sec:q_eff}

In this section, we assess the accuracy of approximations \eqref{final_quant}, \eqref{final_quant_alpha_zero} and \eqref{final_quant_alpha_zero2b}. Using a black line we depict $\mathbb{E}_{\mathbb{Z}_n} \theta(\mathbb{Z}_n)$ obtained via Monte Carlo simulations. Depending on the value of $\alpha$, in Figures~\ref{worsen_figure111}--\ref{quant_fig_2} approximation \eqref{final_quant} or \eqref{final_quant_alpha_zero} is shown using a red line. In Figures~\ref{quant_without1}--\ref{quant_without2}, approximation \eqref{final_quant_alpha_zero2b} is  depicted with a red line. From the figures below we can see that all approximations are generally very accurate. Approximation \eqref{final_quant_alpha_zero} is much more accurate  than approximation \eqref{final_quant} across all choices of $\delta$ and $n$ and this can be explained by the additional term taken in the general expansion; see Section~\ref{sec:General_expansion}. This high accuracy is also seen with approximation \eqref{final_quant_alpha_zero2b}. The accuracy of approximation \eqref{final_quant} seems to worsen for large $\delta$, $n$ and $d$ not too large like $d=20$, see Figures~\ref{worsen_figure1}--\ref{worsen_figure2}. For $d=50$, all approximations are extremely accurate for all choices of $\delta$ and $n$. Figures~\ref{worsen_figure111}--\ref{quant_fig_2} very clearly demonstrate the $\delta$-effect implying that  a sensible choice of $\delta$ is crucial for good quantization.

\begin{figure}[h]
\centering
\begin{minipage}{.5\textwidth}
  \includegraphics[width=1\linewidth]{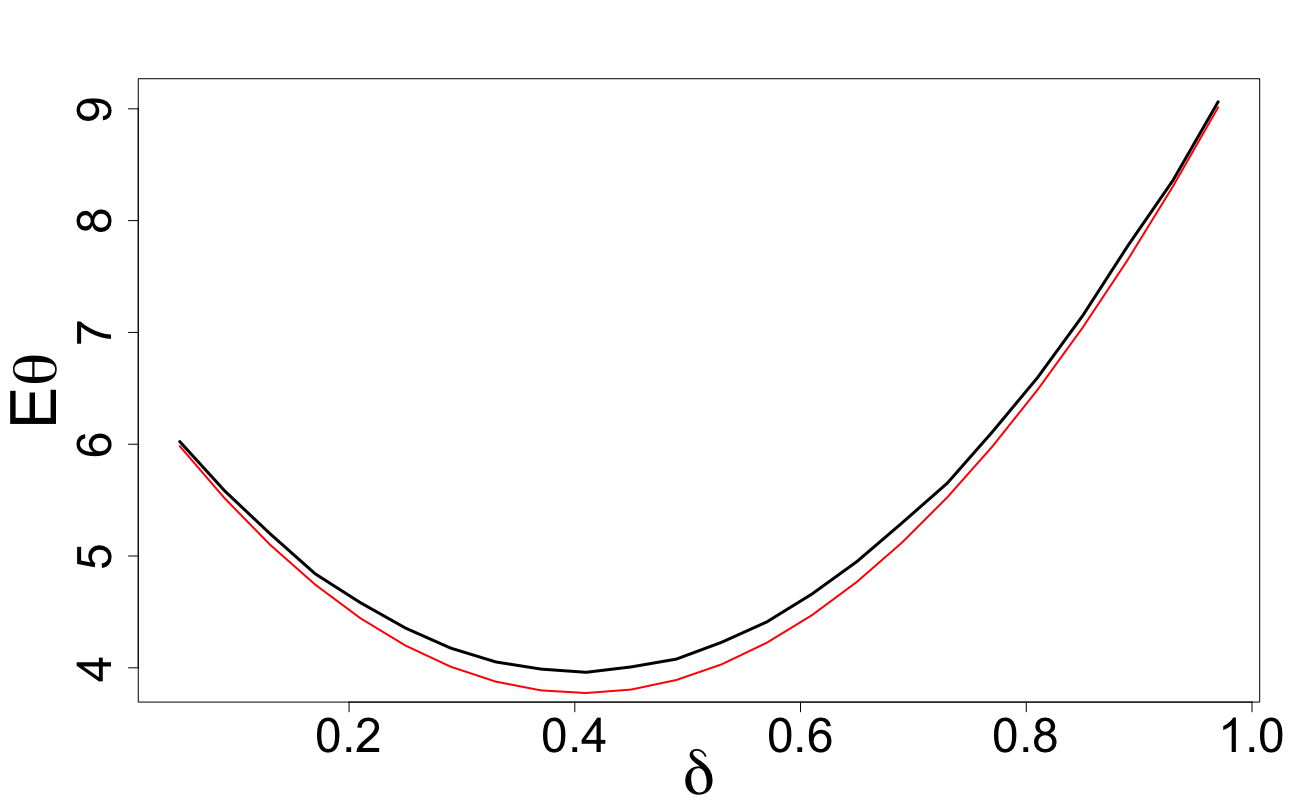}
  \caption{$\mathbb{E} \theta(\mathbb{Z}_n)$ and approximation \eqref{final_quant}$: d=20$, \\$\alpha=1$, $n=500$ }
\label{worsen_figure111}
\end{minipage}%
\begin{minipage}{.5\textwidth}
  \centering
  \includegraphics[width=1\linewidth]{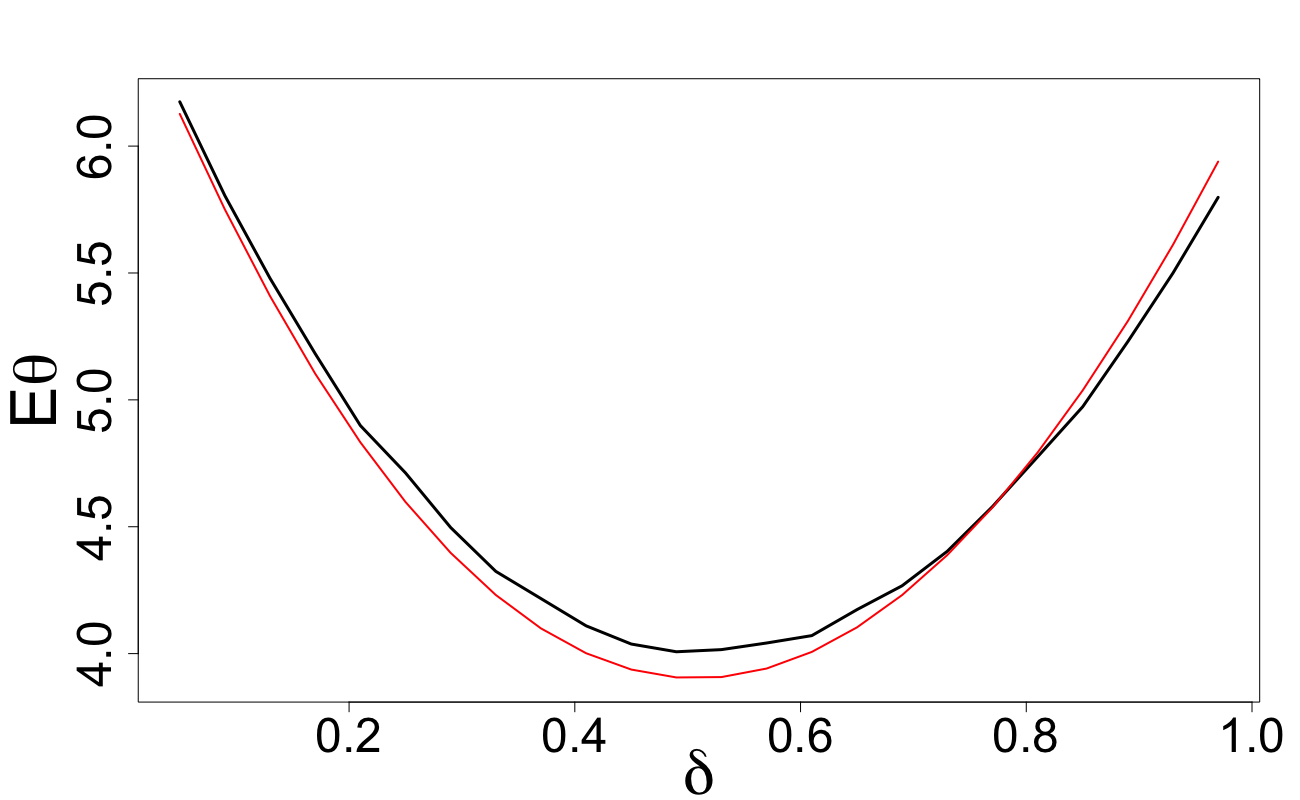}
  \caption{$\mathbb{E} \theta(\mathbb{Z}_n)$ and approximation \eqref{final_quant}$: d=20$, \\$\alpha=0.5$, $n=500$ }
  \label{worsen_figure112}
\end{minipage}
\end{figure}

\begin{figure}[!h]
\centering
\begin{minipage}{.5\textwidth}
  \includegraphics[width=1\linewidth]{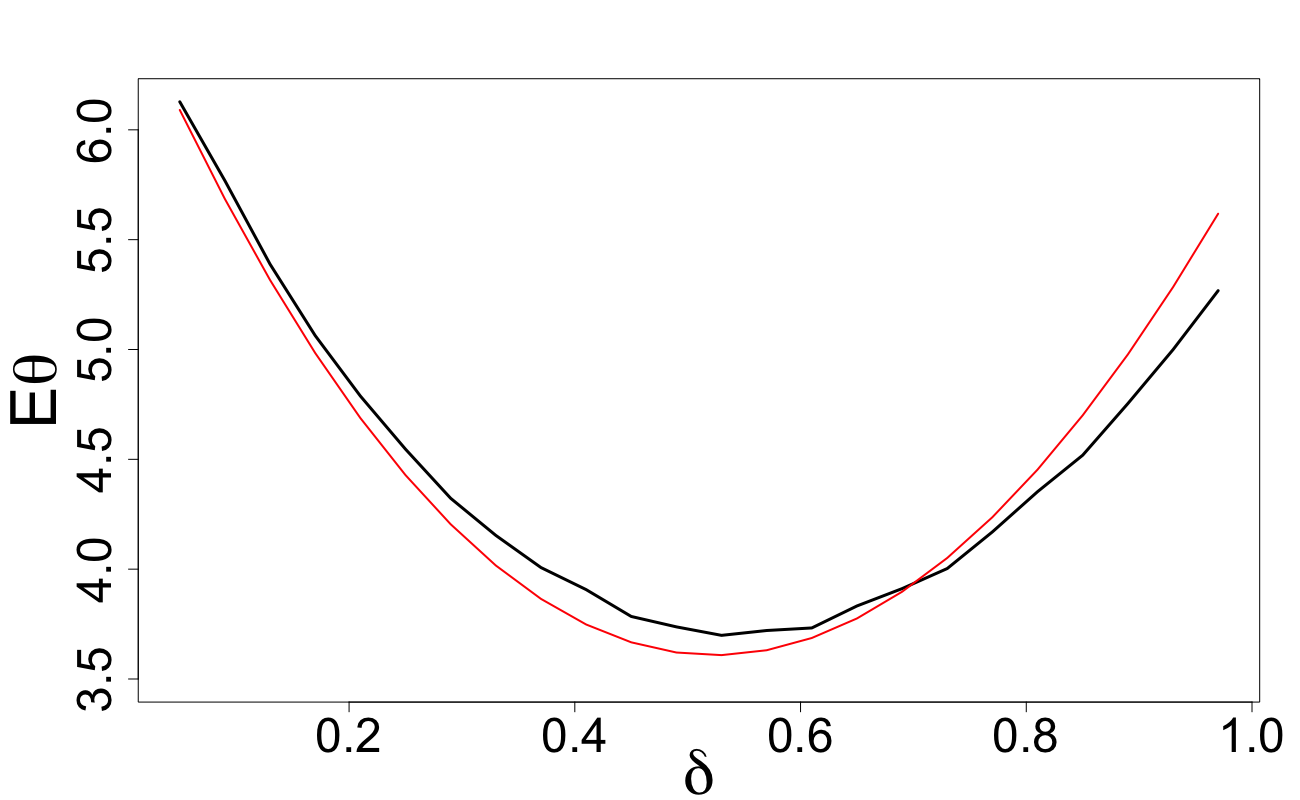}
  \caption{$\mathbb{E} \theta(\mathbb{Z}_n)$ and approximation \eqref{final_quant}: $d=20$,\\ $\alpha=0.5$, $n=1000$ }
  \label{worsen_figure1}
\end{minipage}%
\begin{minipage}{.5\textwidth}
  \centering
  \includegraphics[width=1\linewidth]{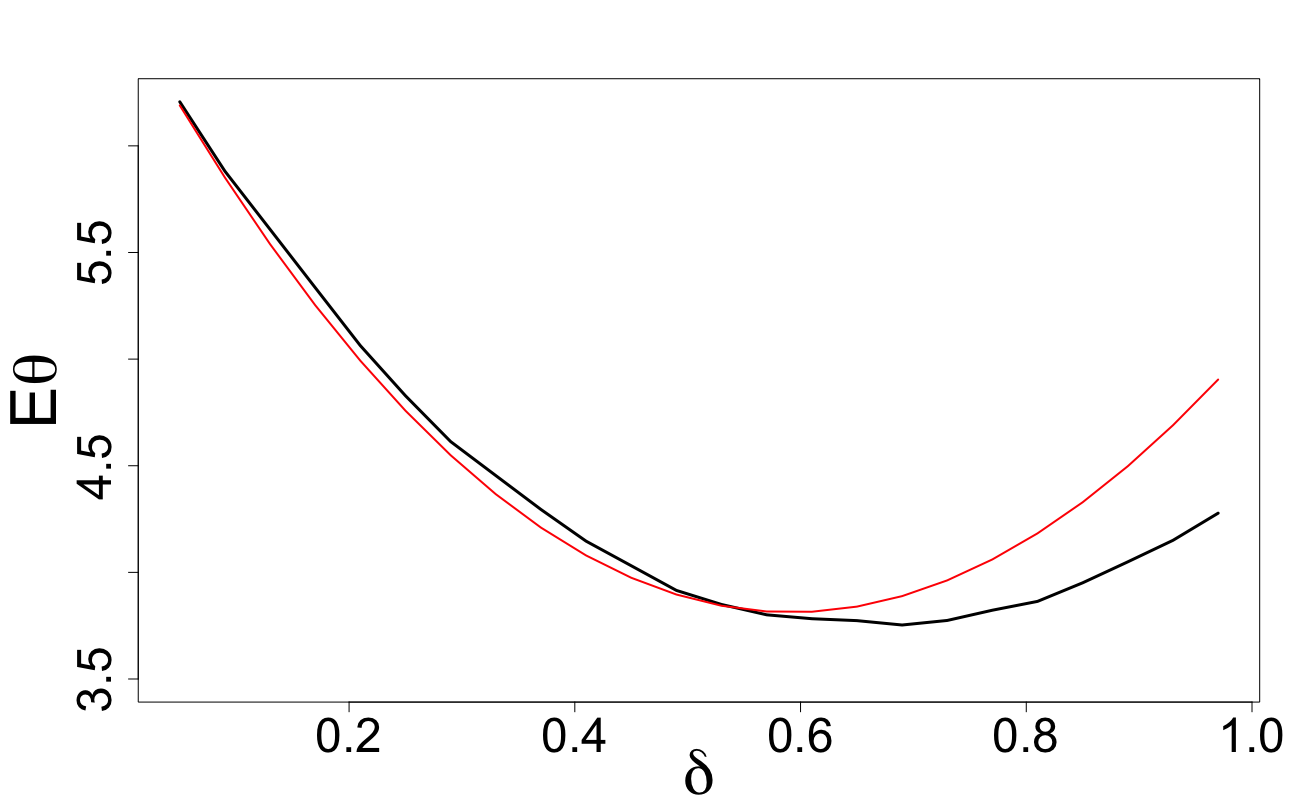}
  \caption{$\mathbb{E} \theta(\mathbb{Z}_n)$ and approximation \eqref{final_quant}: $d=20$,\\ $\alpha=1$, $n=1000$ }
    \label{worsen_figure2}
\end{minipage}
\end{figure}

\begin{figure}[!h]
\centering
\begin{minipage}{.5\textwidth}
  \includegraphics[width=1\linewidth]{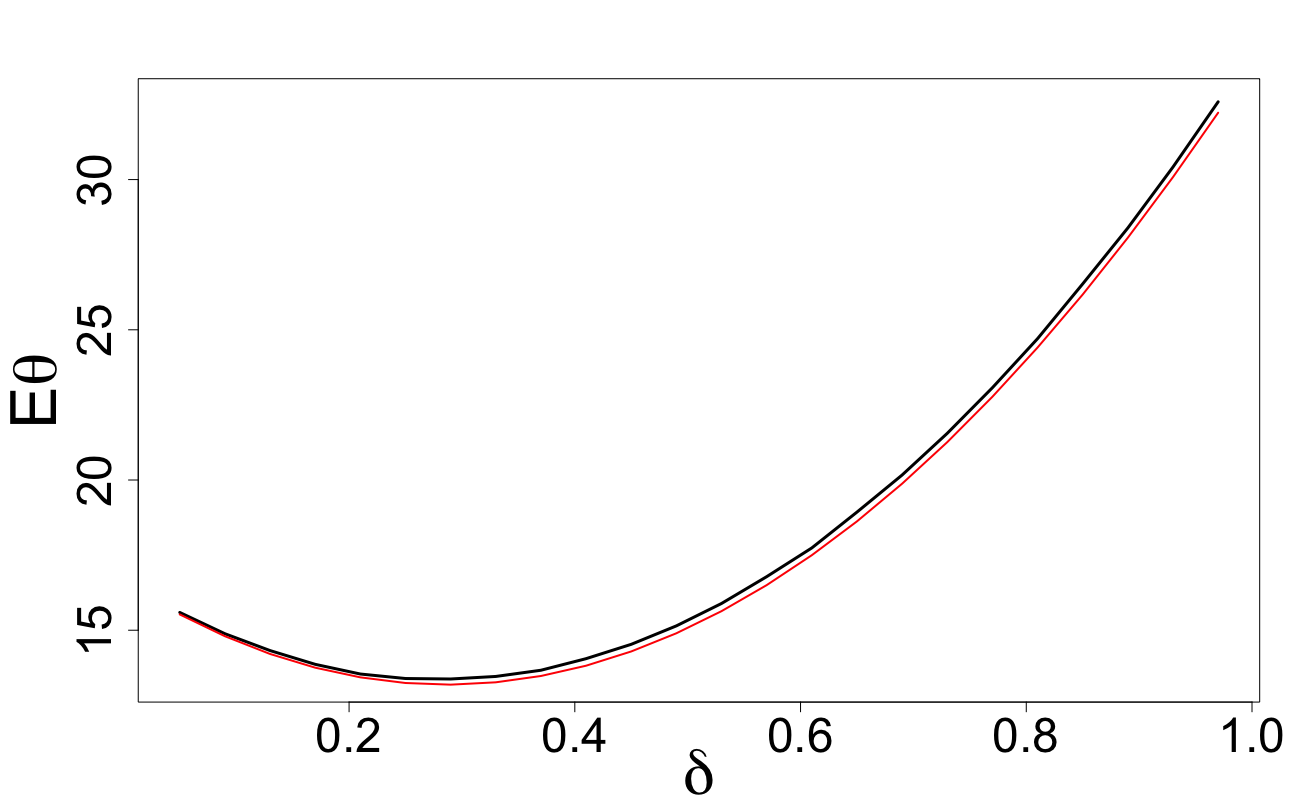}
  \caption{$\mathbb{E} \theta(\mathbb{Z}_n)$ and approximation \eqref{final_quant}: $d=50$, \\ $\alpha=0.1$, $n=1000$ }
\end{minipage}%
\begin{minipage}{.5\textwidth}
  \centering
  \includegraphics[width=1\linewidth]{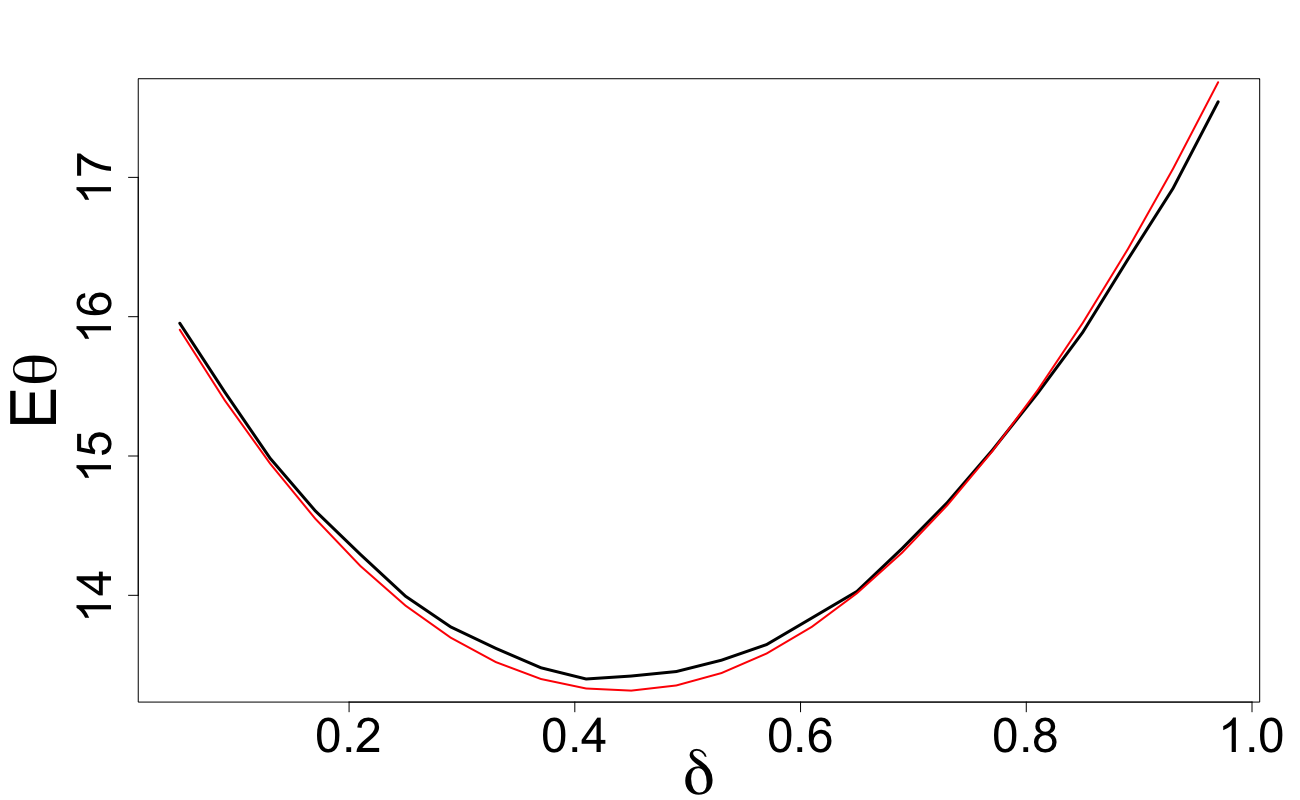}
  \caption{$\mathbb{E} \theta(\mathbb{Z}_n)$ and approximation \eqref{final_quant}: $d=50$, \\ $\alpha=1$, $n=1000$ }
    \label{quant_fig_2}
\end{minipage}
\end{figure}

\begin{figure}[!h]
\centering
\begin{minipage}{.5\textwidth}
  \includegraphics[width=1\linewidth]{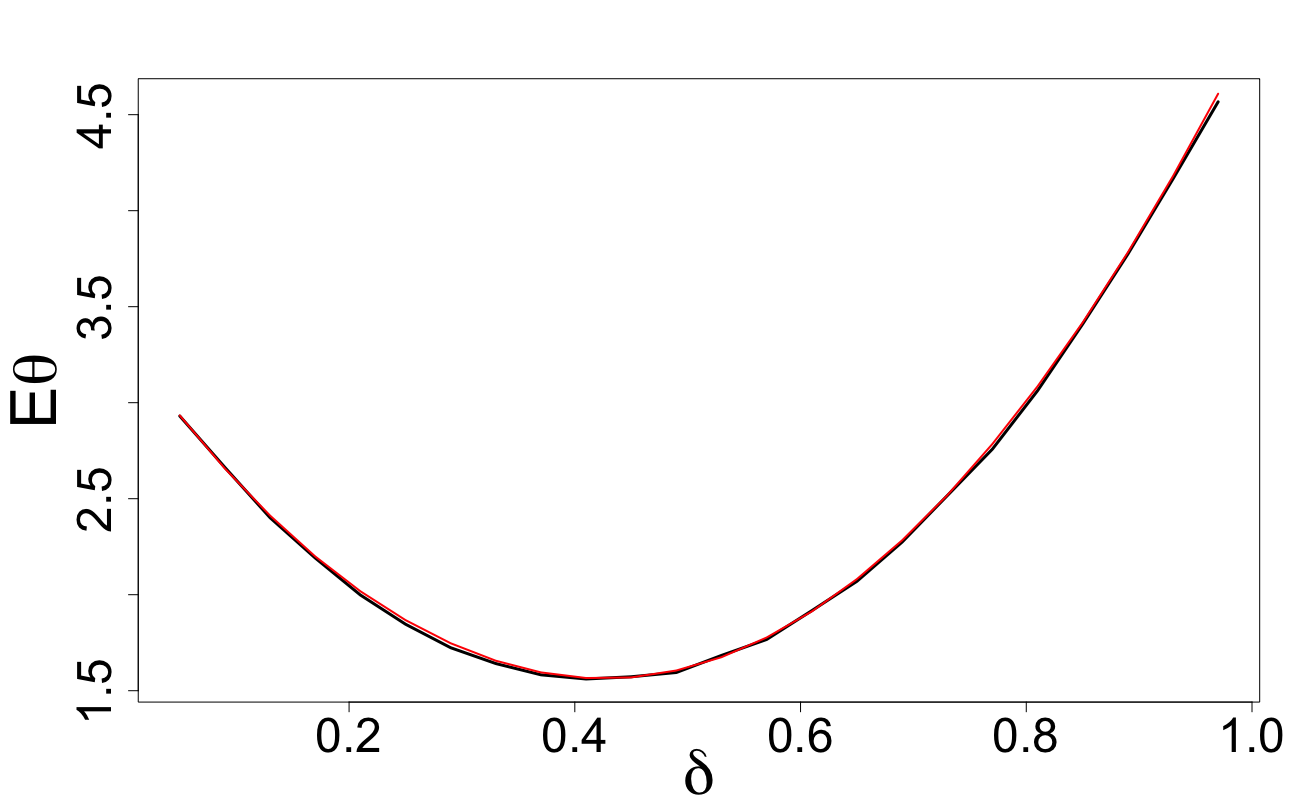}
  \caption{$\mathbb{E} \theta(\mathbb{Z}_n)$ and approximation \eqref{final_quant_alpha_zero}: $d=10$,\\$\alpha=0$, $n=100$ }
  \label{quant_fig_1}
\end{minipage}%
\begin{minipage}{.5\textwidth}
  \centering
  \includegraphics[width=1\linewidth]{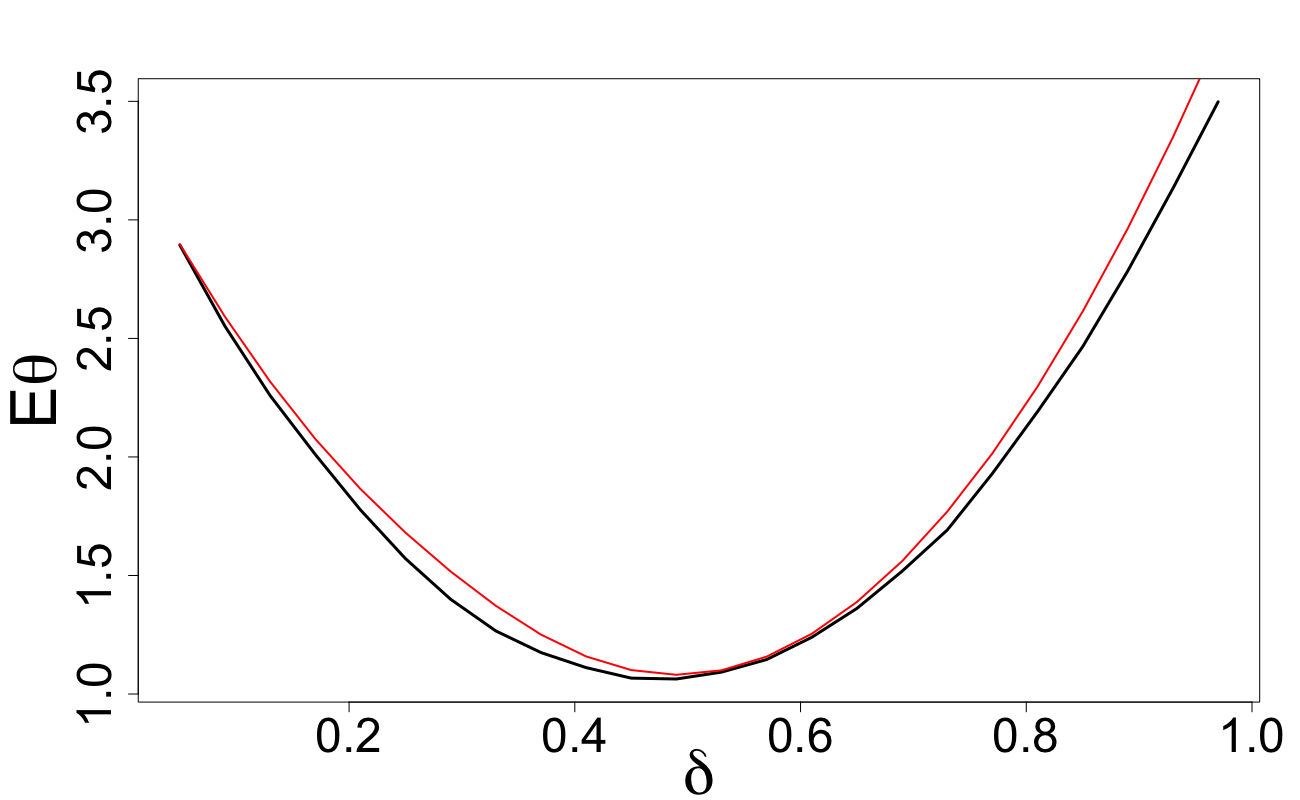}
  \caption{$\mathbb{E} \theta(\mathbb{Z}_n)$ and approximation \eqref{final_quant_alpha_zero}: $d=10$, \\$\alpha=0$, $n=500$ }
\end{minipage}
\end{figure}

\begin{figure}[h]
\centering
\begin{minipage}{.5\textwidth}
  \includegraphics[width=1\linewidth]{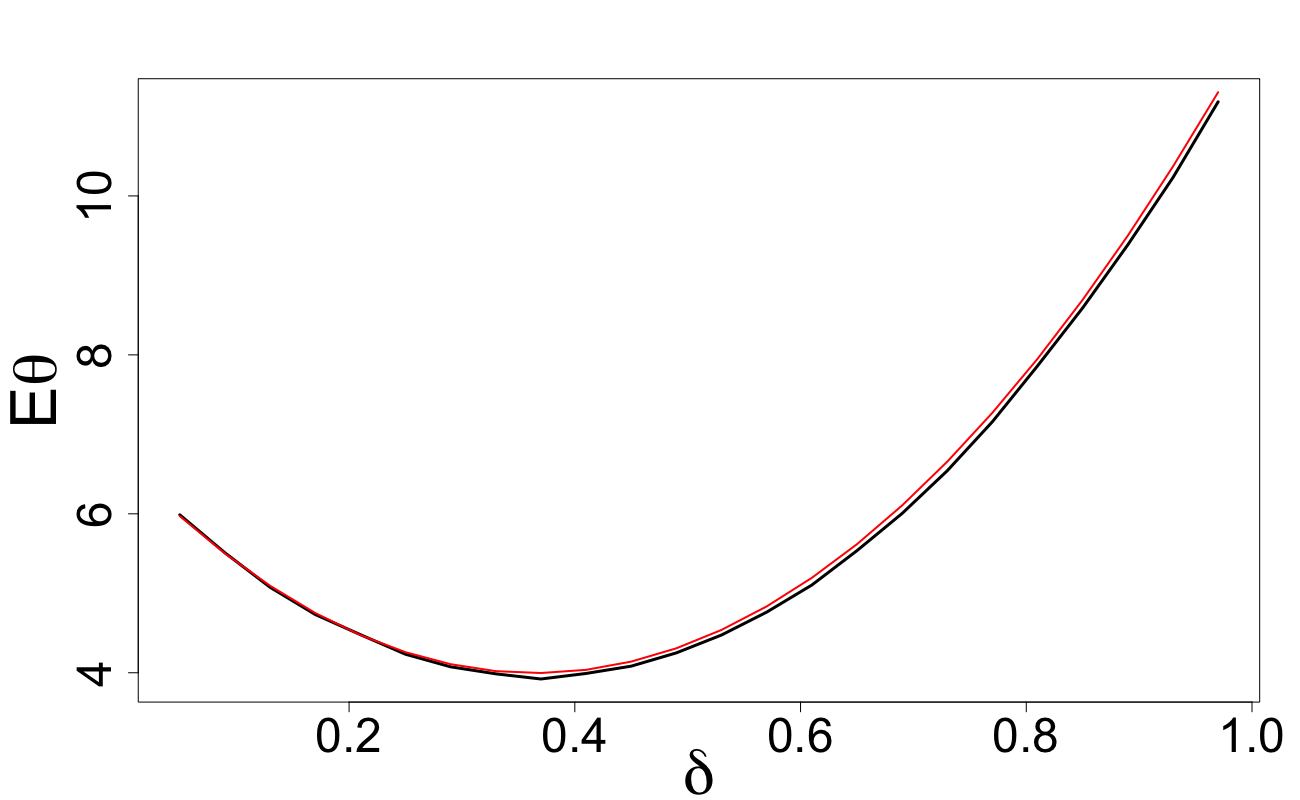}
  \caption{$\mathbb{E} \theta(\mathbb{Z}_n)$ and approximation \eqref{final_quant_alpha_zero}: $d=20$,\\ $\alpha=0$, $n=500$ }
\end{minipage}%
\begin{minipage}{.5\textwidth}
  \centering
  \includegraphics[width=1\linewidth]{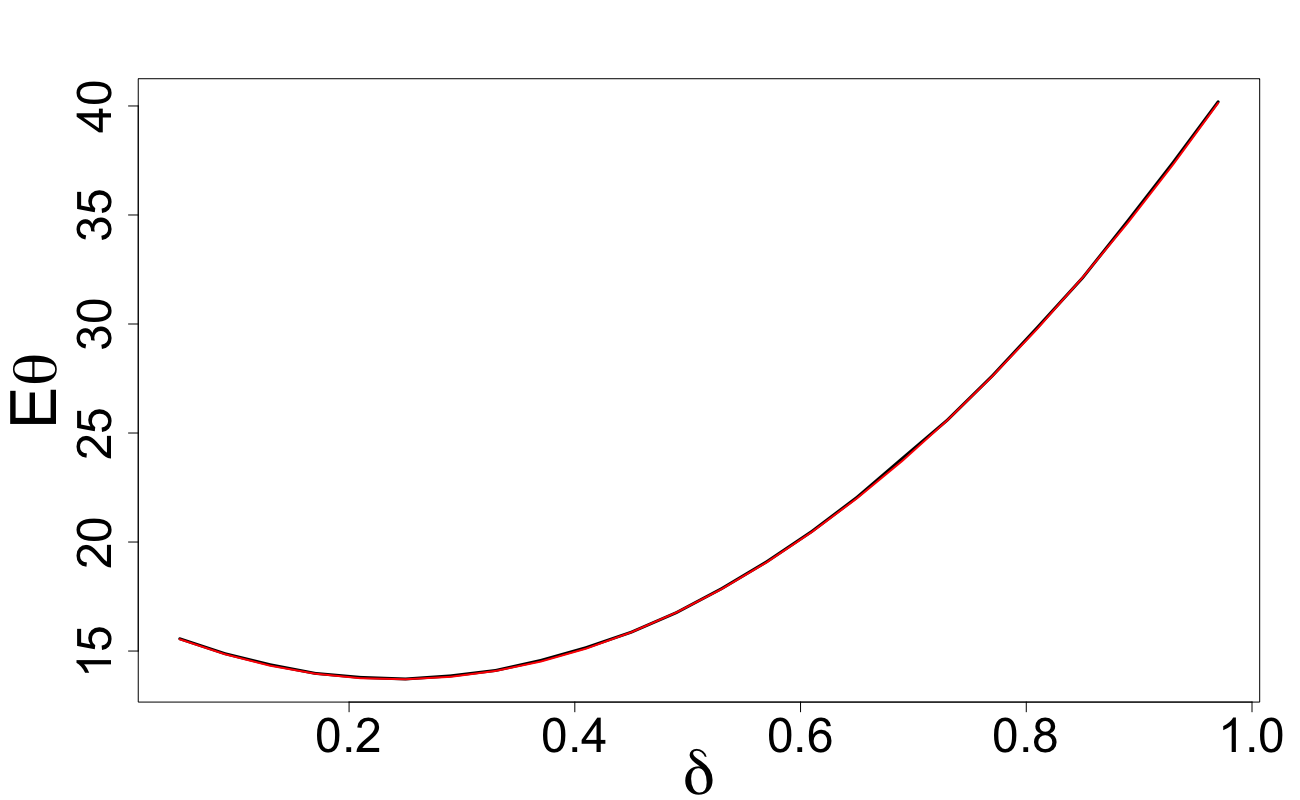}
  \caption{$\mathbb{E} \theta(\mathbb{Z}_n)$ and approximation \eqref{final_quant_alpha_zero}: $d=50$, \\$\alpha=0$, $n=500$ }
\end{minipage}
\end{figure}

\begin{figure}[!h]
\centering
\begin{minipage}{.5\textwidth}
  \includegraphics[width=1\linewidth]{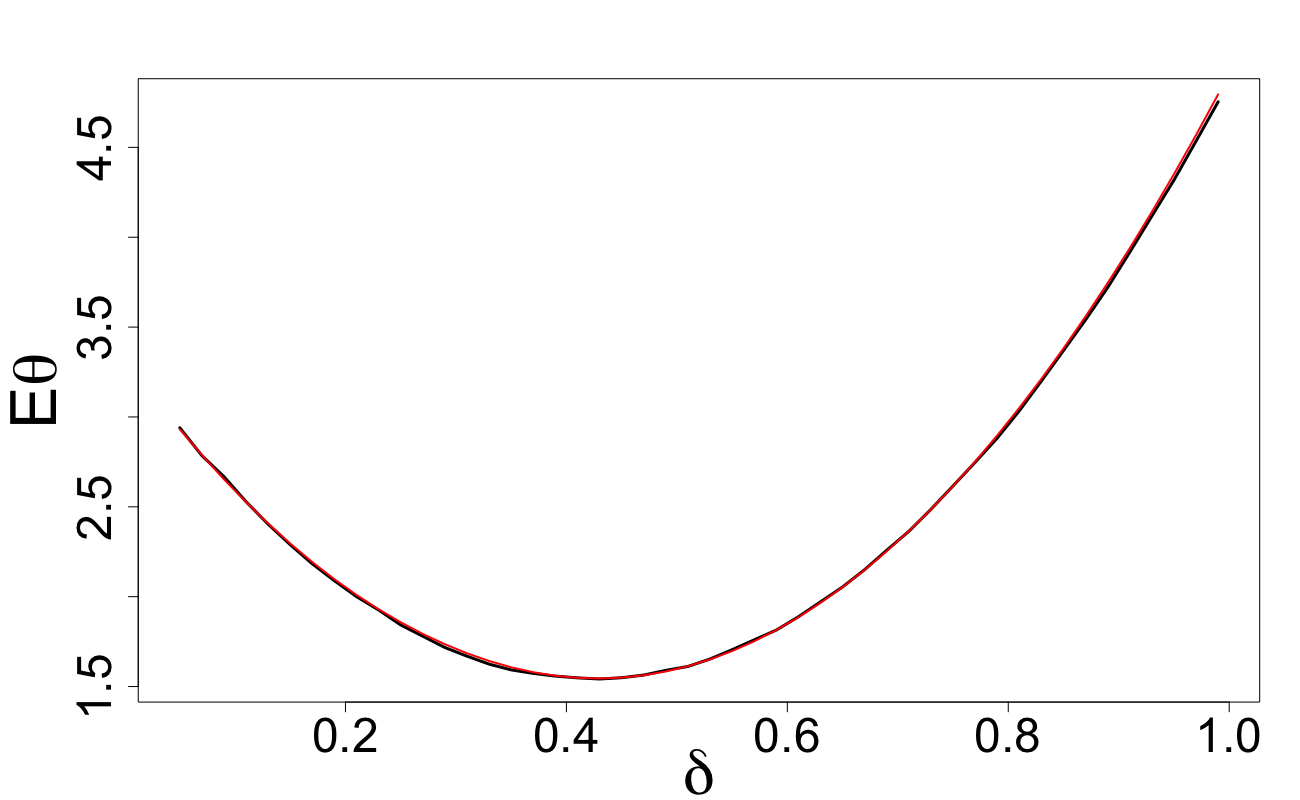}
  \caption{$\mathbb{E} \theta(\mathbb{Z}_n)$ and approximation \eqref{final_quant_alpha_zero2b}: $d=10$, \\  $n=100$ }
  \label{quant_without1}
\end{minipage}%
\begin{minipage}{.5\textwidth}
  \centering
  \includegraphics[width=1\linewidth]{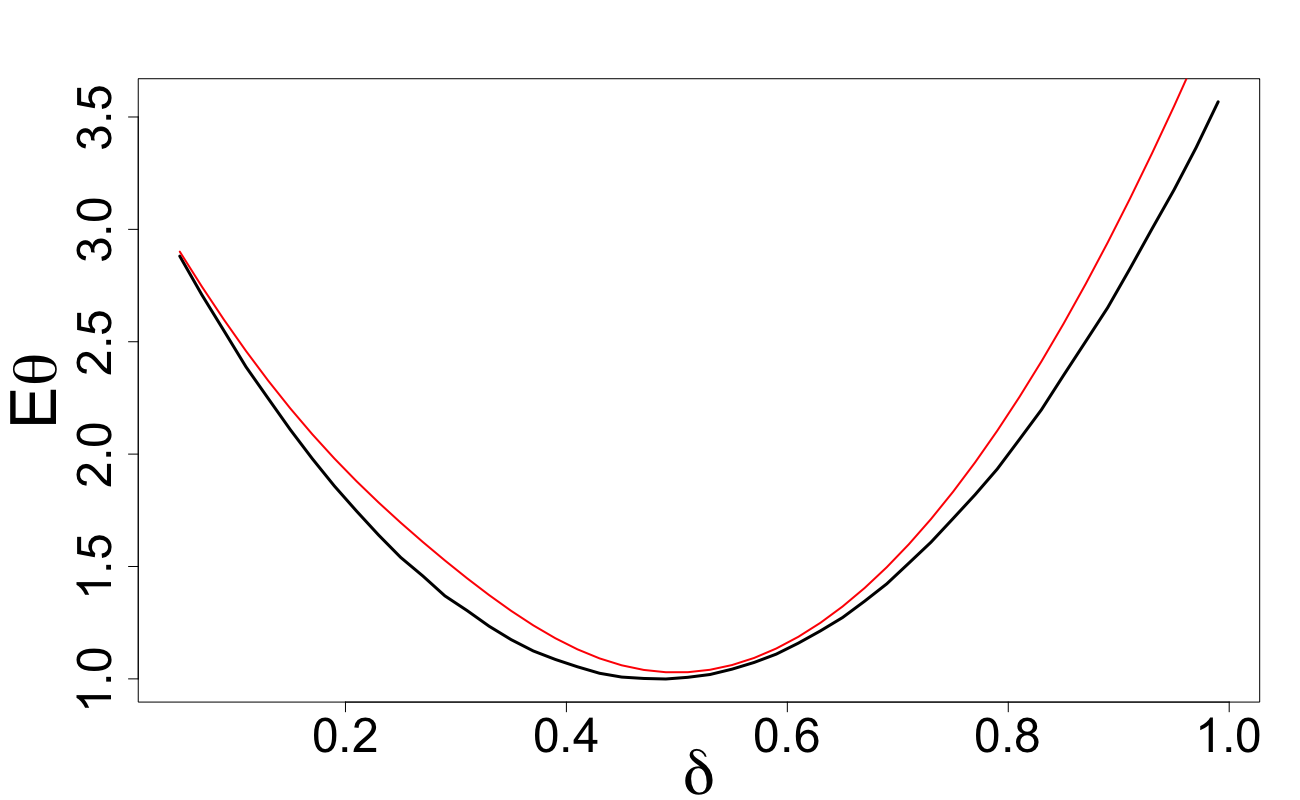}
  \caption{$\mathbb{E} \theta(\mathbb{Z}_n)$ and approximation \eqref{final_quant_alpha_zero2b}: $d=10$, \\  $n=500$ }
\end{minipage}
\end{figure}
\clearpage

\begin{figure}[!h]
\centering
\begin{minipage}{.5\textwidth}
  \includegraphics[width=1\linewidth]{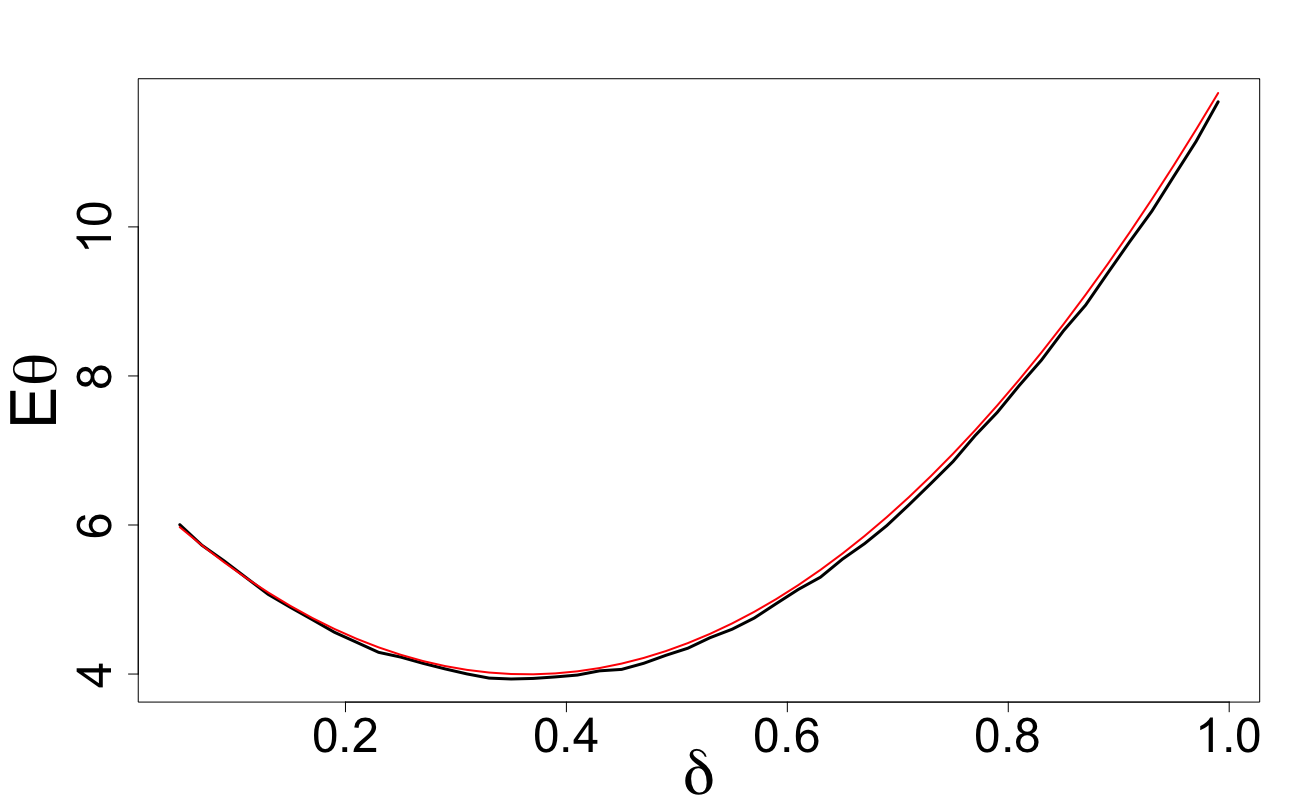}
  \caption{$\mathbb{E} \theta(\mathbb{Z}_n)$ and approximation \eqref{final_quant_alpha_zero2b}: $d=20$, \\ $n=500$ }
\end{minipage}%
\begin{minipage}{.5\textwidth}
  \centering
  \includegraphics[width=1\linewidth]{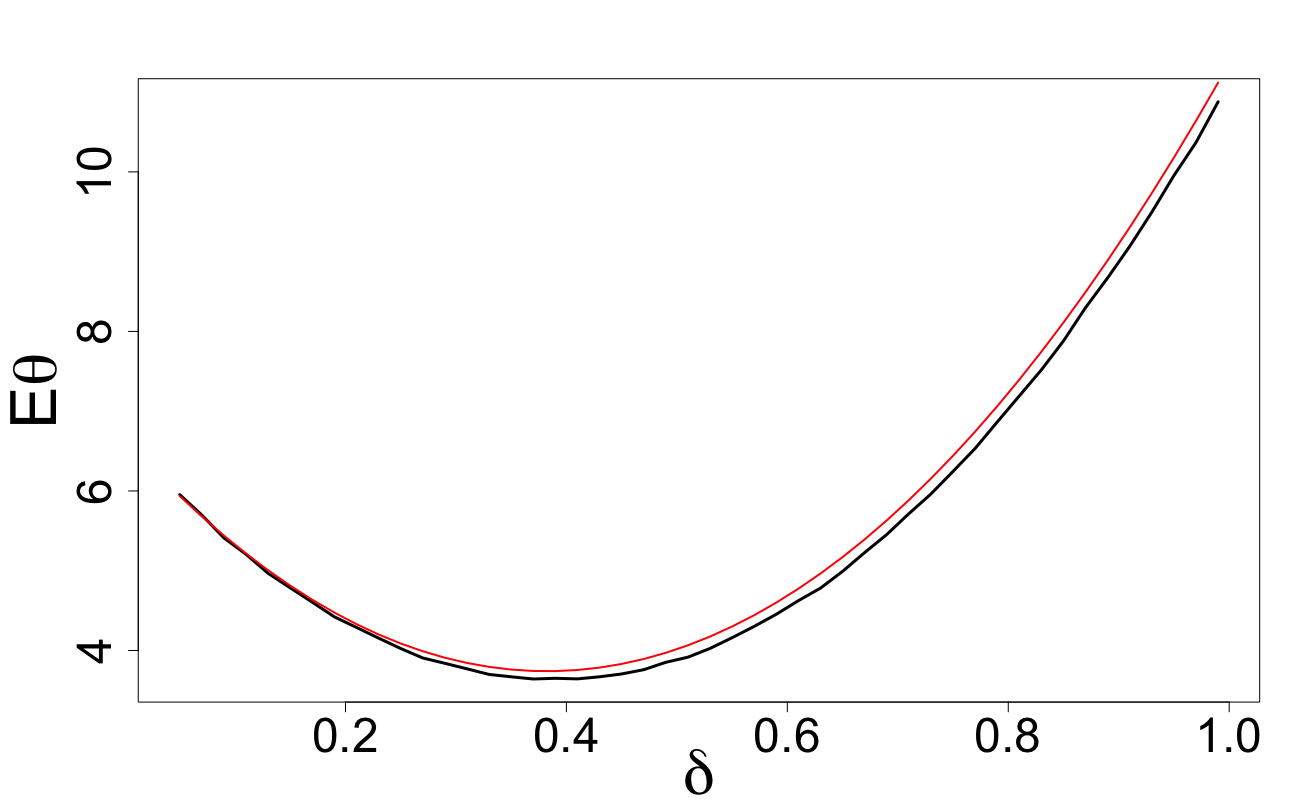}
  \caption{$\mathbb{E} \theta(\mathbb{Z}_n)$ and approximation \eqref{final_quant_alpha_zero2b}: $d=20$, \\  $n=1000$ }
  \label{quant_without2}
\end{minipage}
\end{figure}

\section{Comparative numerical studies of covering properties for several designs}
\label{sec:numeric2}

Let us extend the range of designs considered above by additing the following two designs.\\

{\bf Design 3.} {\it
$Z_1, \ldots, Z_n$ are taken from a low-discrepancy Sobol's sequence on  the cube  $[-\delta,\delta]^d$.
}\\

{\bf Design 4.} {\it
$Z_1, \ldots, Z_n$ are taken from the minimum-aberration $2^{d-k}$ fractional factorial design on the vertices of the cube  $[-\delta,\delta]^d$.
}\\

Unlike Designs 1, 2a, 2b and 3, Design 4 is non-adaptive and defined only for a particular $n$ of the form $n=2^{d-k}$ with some $k\geq 0$.
We have included this design into the list of all designs as "the golden standard". In view of the numerical study in \cite{us} and theoretical arguments in \cite{us1}, Design 4 with $k=1$ and optimal $\delta$ provides the best  quantization we were able to find; moreover, we have conjectured
in \cite{us1} that  Design 4 with $k=1$ and  optimal $\delta$ provides minimal normalized mean squared quantization error for all designs with $n \leq 2^d$. We repeat, Design 4 is defined for one particular value of $n$ only.

\subsection{Covering comparisons}
In Tables~\ref{table_d_10}--\ref{Table_d_20},  we present results of Monte Carlo simulations where we have computed  the smallest values of $r$ required to achieve the 0.9-coverage on average
(on average, for Designs~1, 2a, 2b). The value inside the brackets shows the  value of $\delta$ required to obtain the 0.9-coverage.

\begin{table}[h]
\centering
\begin{tabular}{ |p{2.5cm}||p{2cm}|p{2cm}|p{2cm}|p{2cm}|   }
 \hline
 \multicolumn{5}{|c|}{$d=10$} \\
 \hline
  & $n=64$ &$n=128$& $n=512$ & $n=1024$\\
 \hline
Design 1, $\alpha=0.5$&   1.629 (0.58) &  1.505 (0.65)  & 1.270 (0.72)& 1.165 (0.75)  \\
Design 1, $\alpha=1.5$&  1.635 (0.80)  & 1.525 (0.88)  &1.310 (1.00) & 1.210 (1.00) \\
Design 2a &  1.610 (0.38)    &  1.490 (0.46)   & 1.228 (0.50)    &1.132 (0.50)     \\
Design 2b &  1.609 (0.41)    & 1.475  (0.43)     & 1.178 (0.49)   & 1.075 (0.50)   \\
Design 3 & 1.595 (0.72)  & 1.485 (0.80) & 1.280 (0.85) & 1.170 (0.88)   \\
{   Design 3, $\delta=1$}& 1.678 (1.00)   &  1.534 (1.00) & 1.305  (1.00) & 1.187  (1.00) \\
Design 4  & 1.530 (0.44) & 1.395 (0.48)  & 1.115 (0.50)   & 1.075 (0.50)  \\
 \hline
\end{tabular}
\caption{Values of $r$ and $\delta$ (in brackets) to achieve 0.9 coverage for $d=10$.}
\label{table_d_10}
\end{table}
\begin{table}[!h]
\centering
\begin{tabular}{ |p{2.5cm}||p{2cm}|p{2cm}|p{2cm}|p{2cm}|   }
 \hline
 \multicolumn{5}{|c|}{$d=20$} \\
 \hline
  & $n=64$ &$n=128$& $n=512$ & $n=1024$\\
 \hline
Design 1, $\alpha=0.5$& 2.540 (0.44)  & 2.455 (0.48)  &  2.285 (0.55) &2.220 (0.60)  \\
Design 1, $\alpha=1.5$& 2.545 (0.60) & 2.460 (0.65)  & 2.290 (0.76) &2.215 (0.84)  \\
Design 2a &   2.538 (0.28)   &   2.445 (0.30)  &  2.270 (0.36)  & 2.180  (0.42)   \\
Design 2b &  2.538 (0.29)    &    2.445 (0.30)  & 2.253 (0.37)   &  2.173 (0.42)  \\
Design 3  &  2.520 (0.50)   & 2.445 (0.60)  &  2.285 (0.68) & 2.196 (0.72)  \\
{   Design 3, $\delta=1$}&  2.750 (1.00)   &  2.656 (1.00) &  2.435 (1.00) &  2.325 (1.00) \\
Design 4 & 2.490 (0.32) & 2.410 (0.35) &  2.220 (0.40) &2.125 (0.44)   \\

 \hline
\end{tabular}
\caption{Values of $r$ and $\delta$ (in brackets) to achieve 0.9 coverage for $d=20$.}
\label{Table_d_20}
\end{table}

From Tables~\ref{table_d_10}--\ref{Table_d_20}  we draw the following conclusions:
\begin{itemize}
\item Designs 2a and especially 2b provide very high quality coverage (on average) whilst being online  procedures (that is, nested designs);
\item Design 2b has significant benefits over Design 2a for values of $n$ close to $2^d$;
\item properly $\delta$-tuned deterministic non-nested Design 4 provides superior  covering;
\item coverage properties of   $\delta$-tuned  low-discrepancy sequences are much better than of the original low-discrepancy sequences;
\item coverage  of an unadjusted low-discrepancy sequence is  poor.
\end{itemize}

In Figures~\ref{function_of_r1}--\ref{function_of_r2}, after fixing $n$ and $\delta$, we plot $C_d(\mathbb{Z}_n,r) $ as a function of $r$
for the following designs: Design~1 with $\alpha=1$ (red line), Design 2a (blue line), Design 2b  (green line) and Design 3 with  $\delta=1$ (black line).  For Design 1 with $\alpha=1$, Design 2a and Design 2b, we have used approximations \eqref{eq:psi2-1}, \eqref{eq:psi4} and \eqref{eq:prod7a} respectively to depict $C_d(\mathbb{Z}_n,r) $ whereas for Design 3, we have used Monte Carlo simulations.
For the first three designs, depending of the choice of $n$, the value of $\delta$ has been fixed based on the optimal value for quantization; these are the values inside the brackets in Tables~\ref{d_10_quantization}--\ref{d_20_quantization}.

From Figure~\ref{function_of_r1}, we see that Design 2b is superior and uniformly dominates all other designs for this choice of $d$ and $n$ (at least when the level of coverage is greater than 1/2). In Figure~\ref{function_of_r2}, since $n<<2^d$, the values of $C_d(\mathbb{Z}_n,r) $ for Designs 2a and 2b practically coincide and the green line hides under the blue. In both figures we see that Design 3 with an unadjusted $\delta$ provides a very inefficient covering.

\begin{figure}[!h]
\centering
\begin{minipage}{.5\textwidth}
  \includegraphics[width=1\linewidth]{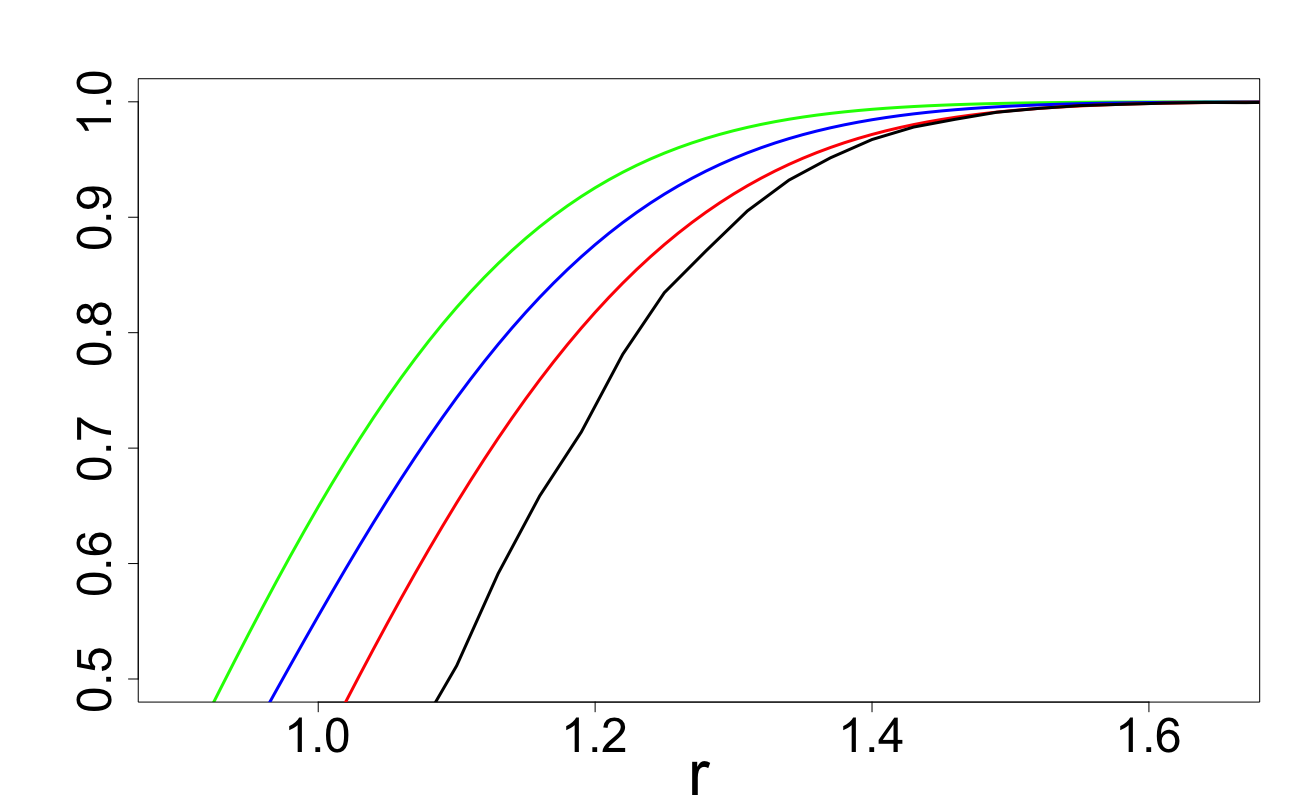}
  \caption{$C_d(\mathbb{Z}_n,r) $ as a function of $r$ for several \\designs: $d=10$, $n=512$ }
  \label{function_of_r1}
\end{minipage}%
\begin{minipage}{.5\textwidth}
  \centering
  \includegraphics[width=1\linewidth]{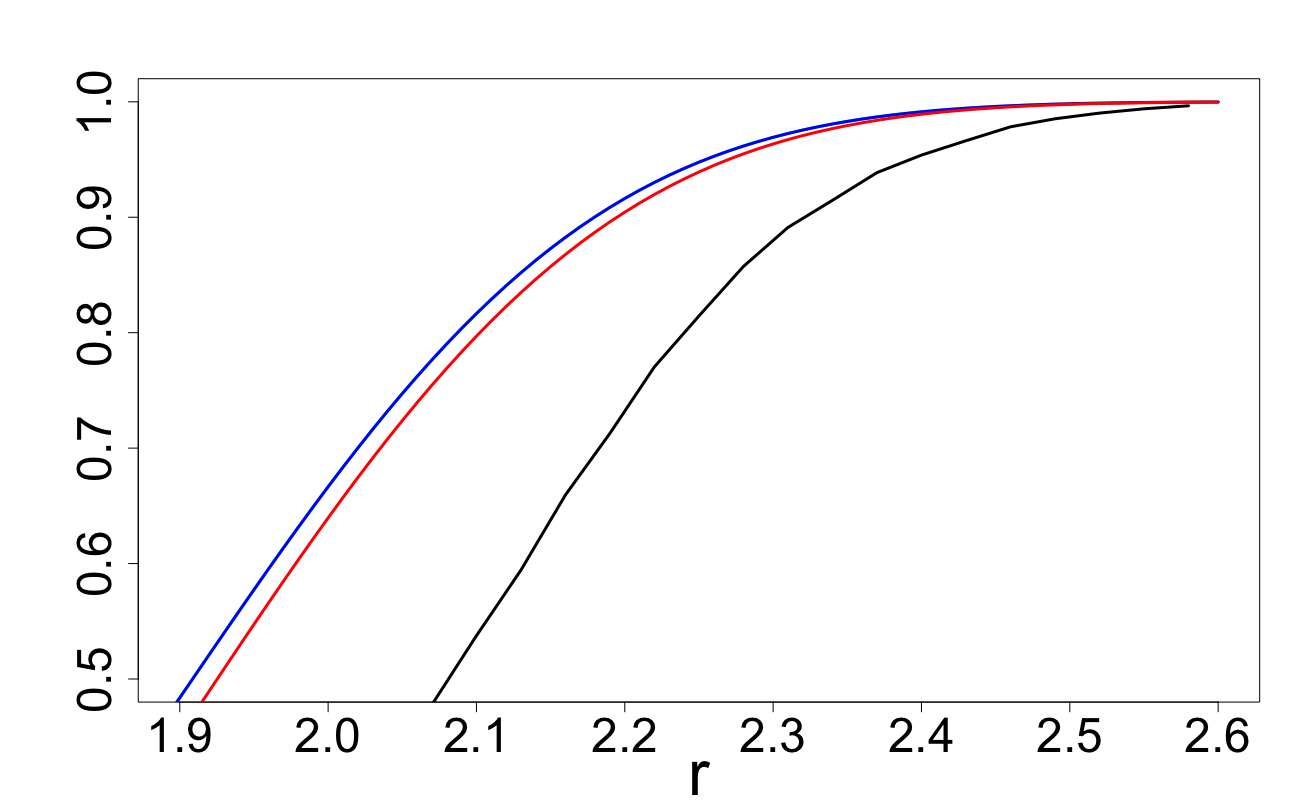}
  \caption{$C_d(\mathbb{Z}_n,r) $ as a function of $r$ for several\\ designs: $d=20$, $n=1024$ }
    \label{function_of_r2}
\end{minipage}
\end{figure}

\subsection{Quantization comparisons}

As follows from results of \cite[Ch.6]{niederreiter1992random}, for efficient covering schemes the order of convergence of the covering radius to 0 as $n \to \infty$ is $n^{-1/d}$. Therefore, for the mean squared distance (which is the quantization error) we should expect the order $n^{-2/d}$ as $n \to \infty$. Therefore, for sake of comparison of quantization errors $\theta_n$ across $n$ we renormalize this error from $\mathbb{E} \theta_n$ to $n^{2/d} \mathbb{E} \theta_n  $.

In Figure~\ref{dom1}--\ref{d_20_quantization}, we present the minimum value of $n^{2/d} \mathbb{E} \theta_n  $ for a selection of designs. In these tables, the value within the brackets corresponds to the value of $\delta$ where the minimum of $n^{2/d} \mathbb{E} \theta_n$ was obtained.

\begin{table}[h]
\centering
\begin{tabular}{ |p{2.5cm}||p{2cm}|p{2cm}|p{2cm}|p{2cm}|   }
 \hline
 \multicolumn{5}{|c|}{$d=10$} \\
 \hline
  & $n=64$ &$n=128$& $n=512$ & $n=1024$\\
 \hline
Design 1, $\alpha=0.5$& 4.072 (0.56)  &  4.013 (0.60) & 3.839 (0.68) & 3.770   (0.69) \\
Design 1, $\alpha=1$ & 4.153 (0.68)   & 4.105 (0.72) &  3.992 (0.80) &3.925 (0.84)  \\
Design 1, $\alpha=1.5$ &  4.164 (0.80)   & 4.137 (0.86)  &  4.069 (0.96) & 4.026 (0.98)    \\
Design 2a &  3.971 (0.38)  & 3.866 (0.44)  & 3.670 (0.48)  & 3.704  (0.50)   \\
Design 2b &  3.955 (0.40)  & 3.798 (0.44)  & 3.453 (0.48) &  3.348  (0.50)   \\
Design 3&   3.998 (0.68)  & 3.973 (0.76) & 3.936 (0.80) &   3.834 (0.82) \\
Design 3, $\delta=1$& 4.569   (1.00)  & 4.425  (1.00) & 4.239 (1.00) & 4.094  (1.00) \\

Design 4 & 3.663 (0.40) &3.548 (0.44)   & 3.221 (0.48)  & 3.348  (0.50)   \\
 \hline
\end{tabular}
\caption{Minimum value of $n^{2/d} \mathbb{E} \theta_n$ and $\delta$ (in brackets) across selected designs;  $d=10$.  }
\label{d_10_quantization}
\end{table}

\begin{table}[h]
\centering
\begin{tabular}{ |p{2.5cm}||p{2cm}|p{2cm}|p{2cm}|p{2cm}|   }
 \hline
 \multicolumn{5}{|c|}{$d=20$} \\
 \hline
  & $n=64$ &$n=128$& $n=512$ & $n=1024$\\
 \hline
 Design 1,  $\alpha=0.5$& 7.541 (0.40) &  7.515 (0.44) & 7.457 (0.52)  &  7.421 (0.54) \\
 Design 1,  $\alpha=1$ & 7.552 (0.52)   & 7.563 (0.56)  &  7.528 (0.64)  & 7.484 (0.68)  \\
 Design 1,  $\alpha=1.5$&7.561 (0.60)   & 7.571 (0.64)   & 7.556 (0.74)   & 7.527 (0.78)  \\
 Design 2a & 7.488 (0.30)  & 7.461 (0.33)   & 7.346  (0.35) & 7.248 (0.39)   \\
 Design 2b & 7.487 (0.29)  &  7.458 (0.34)  & 7.345 (0.36)  & 7.234 (0.40)  \\
 Design 3&  7.445 (0.48)   &7.464 (0.56)  &7.487 (0.64)  &  7.453 (0.66)  \\
Design 3, $\delta=1$& 9.089   (1.00)  &  9.133 (1.00) &  8.871 (1.00) & 8.681  (1.00) \\

 Design 4 &  7.298 (0.32) &  7.270  (0.33)&  7.133 (0.36) & 7.016 (0.40)  \\
 \hline
\end{tabular}
\caption{Minimum value of $n^{2/d} \mathbb{E} \theta_n$ and $\delta$ (in brackets) across selected designs;  $d=20$.  }
\label{d_20_quantization}
\end{table}

In Figure~\ref{dom1}, we depict the c.d.f.'s for the distance $ \varrho(X,\mathbb{Z}_n)$  for  Design 2a with $\delta=0.5$ (in red) and Design 3 with $\delta=0.8$ (in black). We can see that for $d=10$ and $n=512$, Design 2a stochastically dominates Design 3. The style of Figure~\ref{dom2} is the same as figure Figure~\ref{dom1}, however we set $n=1024$ and Design 2a is replaced with Design 2b with $\delta=0.5$ (we also set $\delta=0.82$ for Design 3). Here we see a very clear stochastic dominance of the Design 2b over Design 4. All findings are consistent with Tables~\ref{d_10_quantization} and~\ref{d_20_quantization}. In Figures~\ref{dom1} and~\ref{dom2},  values of the parameter $\delta$ for all designs 
are chosen as numerically optimal, in accordance with 
Table~\ref{d_10_quantization}.


\begin{figure}[!h]
\centering
\begin{minipage}{.5\textwidth}
  \includegraphics[width=1\linewidth]{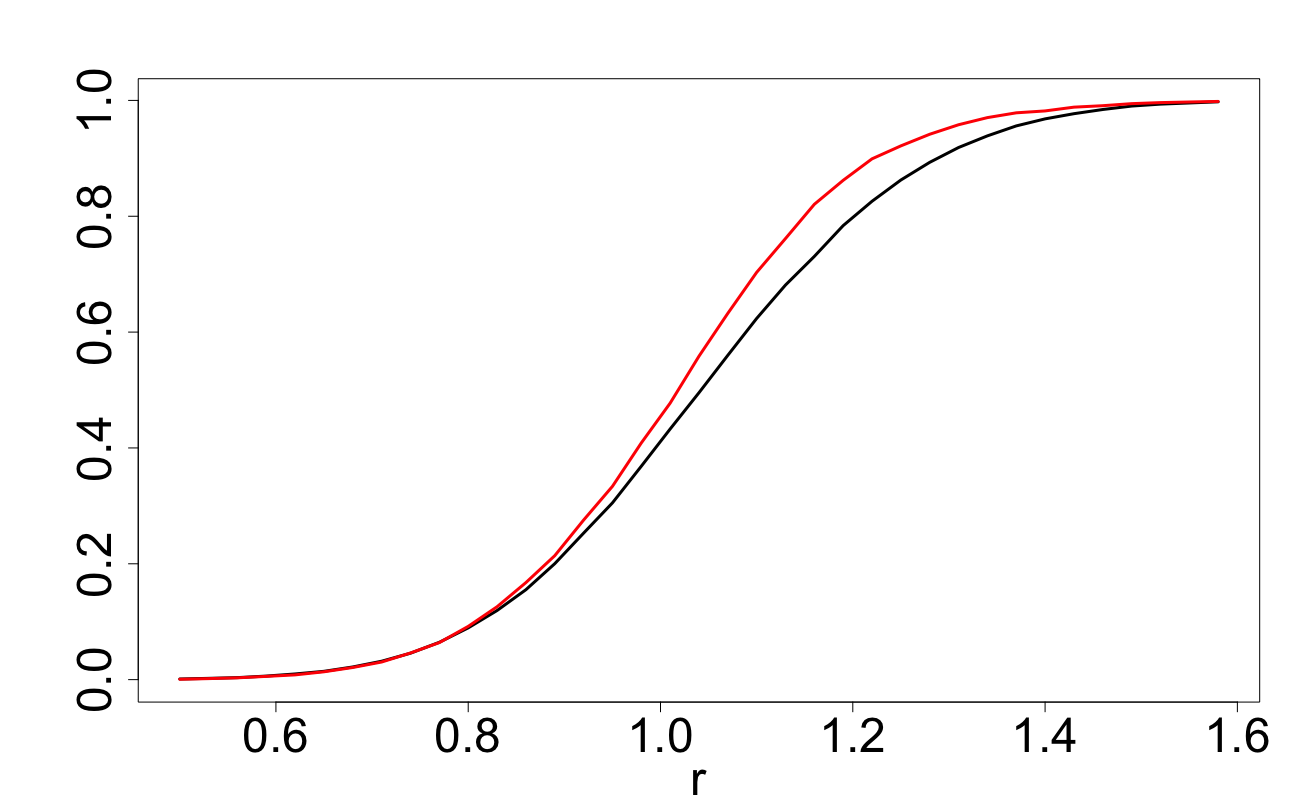}
  \caption{$d=10, n=512$: Design 2a with \\$\delta=0.5$ stochastically dominates  Design 3 \\ with $\delta=0.8$. }
  \label{dom1}
\end{minipage}%
\begin{minipage}{.5\textwidth}
  \centering
  \includegraphics[width=1\linewidth]{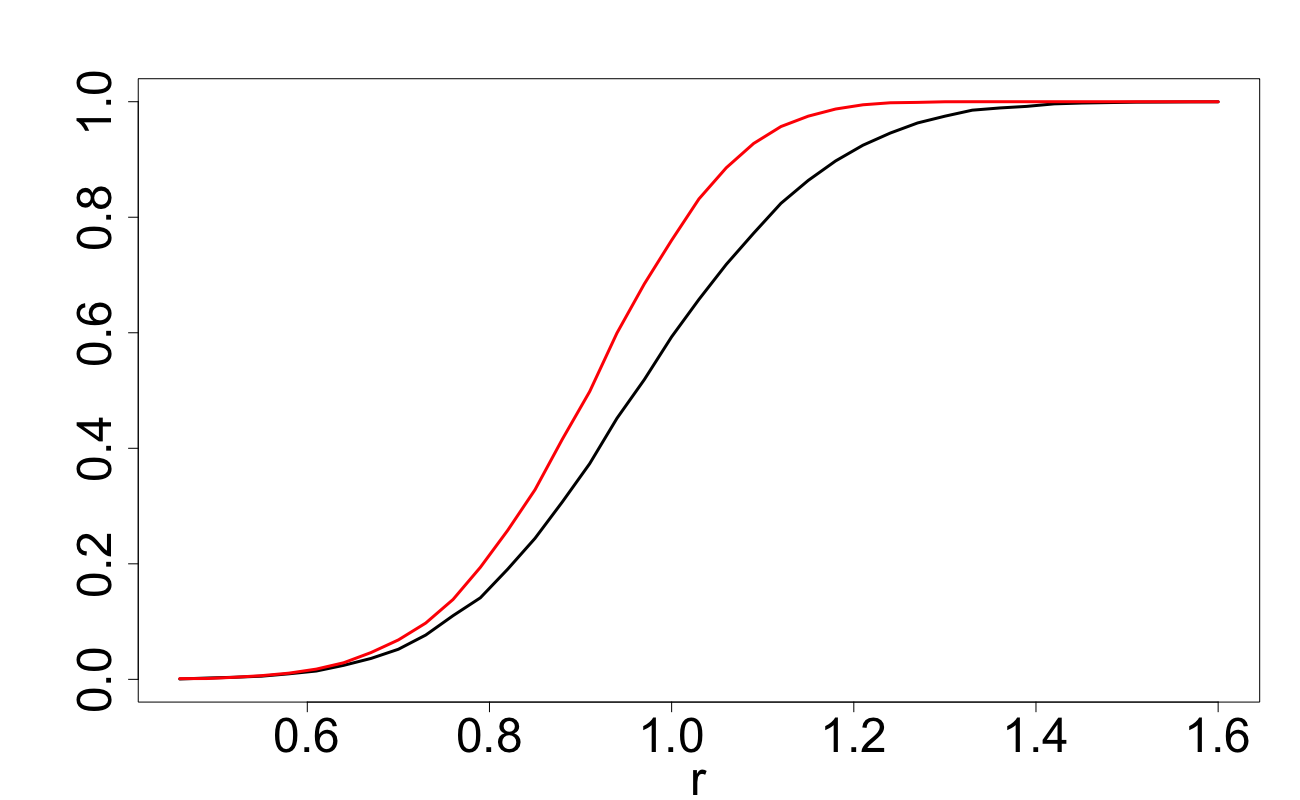}
  \caption{$d=10, n=1024$: Design 2b with \\$\delta=0.5$ stochastically dominates  Design 3 \\ with $\delta=0.82$. }
    \label{dom2}
\end{minipage}
\end{figure}

We make the following conclusions from analyzing results of this section:

\begin{itemize}
  \item Designs 2a and 2b provide very good quantization per point. As expected, Design 2b is superior over Design 2a when $n$ is close to $ 2^d$; see Table~\ref{d_10_quantization}.
  \item Properly $\delta$-tuned non-nested Design 4 is provides the best quantization per point of all designs considered.
  \item Properly $\delta$-tuned Design 3 is comparable in performance to Design 1 but it is not as efficient as Designs 2a, 2b and 4.
  \end{itemize}

\section{Covering and quantization in the $d$-simplex}

\label{sec:simplex}

\subsection{ Characteristics of interest}
Consider the standard orthogonal $d$-simplex
\bea
{\cal S}_d := \left\{(u_1,u_2,\ldots, u_d)\in \mathbb{R}^d \, \big |\,  \sum_{i=1}^{d}u_i \leq 1 \text{ and } u_i\ge0 \text{ for all } i  \right  \} \,
\eea
with   $\text{vol}({\cal S}_d) = {1}/{d!}$.
For a design $\mathbb{Z}_n=\{Z_1, \ldots, Z_n\}$, consider the following two characteristics:
\begin{enumerate}
\item[(a)] the proportion of the simplex ${\cal S}_{d}$ covered by ${\cal B}_d(\mathbb{Z}_n,r)$:
\be
\mbox{$C_d(\mathbb{Z}_n,r):= d!\,$vol$({\cal S}_{d} \cap {\cal B}_d(\mathbb{Z}_n,r))$}\,,
\ee

\item[(b)] $\theta(\mathbb{Z}_n)=\mathbb{E}_U\min_{i=1, \ldots, n} \|U-Z_i\|^2$, 
the mean squared quantization error for $\mathbb{Z}_n$, where  $U=(u_1, \ldots, u_d)$ is a random vector uniformly distributed  in ${\cal S}_d$.
\end{enumerate}
In this section, we investigate whether the $\delta$-effect seen in Sections  \ref{sec:numeric}, \ref{sec:q_eff} and \ref{sec:numeric2}  for the cube is present for the simplex ${\cal S}_d$. We will consider two possible ways of scaling points in ${\cal S}_d$.
Define the two  $\delta$-simplices ${\cal S}^{(\delta)}_{d,1}$ and ${\cal S}^{(\delta)}_{d,2} $ as follows:
\bea
{\cal S}^{(\delta)}_{d,1}:=\delta \cdot {\cal S}_{d} \,,
\eea
\vspace{-0.5cm}
\bea
{\cal S}^{(\delta)}_{d,2} := \left\{(u_1,u_2,\ldots, u_d)\in \mathbb{R}^d \, \big |\,  \sum_{i=1}^{d}u_i \leq \frac{d+\delta}{d+1} \text{ and } u_i\ge \frac{1-\delta}{d+1} \text{ for all } i  \right  \} \, .
\eea

By construction, the value of $\delta$ in ${\cal S}^{(\delta)}_{d,2}$ scales the simplex around its centroid ${ S}^\ast_{d} = \left(\frac{1}{d+1},\frac{1}{d+1},\ldots,\frac{1}{d+1} \right)$, where for $\delta=1$, we have ${\cal S}^{(\delta)}_{d,2}={\cal S}_d$. Simple depictions of ${\cal S}^{(\delta)}_{d,1}$ and ${\cal S}^{(\delta)}_{d,2}$ are given in Figures~\ref{simplex1}--\ref{simplex2}.

\begin{figure}[!h]
\centering
\begin{minipage}{.5\textwidth}
 \begin{center}
\begin{tikzpicture}[scale=5.5]
\draw (0,0) node[anchor=north]{$(0,0)$}
  -- (1,0) node[anchor=north]{$(1,0)$}
  -- (0,1) node[anchor=south]{$(0,1)$}
  -- cycle;

  \draw (0,0) node[anchor=north]{$(0,0)$}
  -- (0.5,0) node[anchor=north]{$(1/2,0)$}
  -- (0,0.5) node[anchor=east]{$(0,1/2)$}
  -- cycle;
  \end{tikzpicture}
\end{center}
\caption{${\cal S}_d$ and ${\cal S}^{(\delta)}_{d,1}$ with $d=2$ and $\delta=0.5$}
\label{simplex1}

 \end{minipage}%
\begin{minipage}{.5\textwidth}
\begin{center}
\begin{tikzpicture}[scale=5.5]
\draw (0,0) node[anchor=north]{$(0,0)$}
  -- (1,0) node[anchor=north]{$(1,0)$}
  -- (0,1) node[anchor=south]{$(0,1)$}
  -- cycle;

  \draw (0.1666,0.1666) node[anchor=north]{$(1/6,\!1/6)$}
  -- (0.666,0.1666) node[anchor=north]{$(2/3,1/6)$}
  -- (0.1666,0.666) node[anchor=south]{$(1/6,\!2/3)$}
  -- cycle;
\draw (1/3,1/3) circle[radius=0.3pt] node[anchor=north]{${ S}^\ast_{d}$};
\end{tikzpicture}
\end{center}
\caption{${\cal S}_d$ and ${\cal S}^{(\delta)}_{d,2}$ with $d=2$ and $\delta=0.5$}
\label{simplex2}

  \end{minipage}
\end{figure}

We will numerically assess  covering and quantization characteristics for the following two designs.\\

{\bf Design S1.} {\it
$Z_1, \ldots, Z_n$ are i.i.d. random vectors uniformly distributed  in the $\delta$-scaled simplex ${\cal S}^{(\delta)}_{d,1} $,
where $\delta \in [0,1] $ is a parameter.\\
}

{\bf Design S2.} {\it
$Z_1, \ldots, Z_n$ are i.i.d. random vectors uniformly distributed  in the $\delta$-scaled simplex ${\cal S}^{(\delta)}_{d,2} $,
where $\delta \in [0,1] $ is a parameter.\\
}

To simulate points $Y$ uniformly distributed in the simplex ${\cal S}_{d} $, we can simply generate
$d$ i.i.d. uniformly distributed points in $[0,1]$, add 0 and 1 to the collection of points  and take the first $d$ spacings (out of the total number $d+1$ of these spacings). Points  $Y' = \delta Y $  and $Y'' = \delta\cdot (Y-{ S}^\ast_{d})+{ S}^\ast_{d} $  are then uniform
in ${\cal S}^{(\delta)}_{d,1}$ and ${\cal S}^{(\delta)}_{d,2}$ respectively. This procedure can be easily performed in R using the package `uniformly'.

\subsection{Numerical investigation of the $\delta$-effect for $d$-simplex}
Using the above procedure, we numerically study characteristics of Designs S1 and  S2.  In Figures~\ref{simplex_1_pic}--\ref{simplex_1_pic_end} we plot ${ C}_d(\mathbb{Z}_n,r) $ as a functions of $\delta \in [0,1]$ across $n, r$ and $d$ for Design S1. The corresponding results for Design S2 are given in Figures~\ref{simplex_2_pic}--\ref{simplex_2_pic_end}. In Figures~\ref{quant_s1}--\ref{quant_s1_end} and Figures~\ref{quant_s2}--\ref{quant_s2_end}, we depict $\mathbb{E}\theta(\mathbb{Z}_n)$ for Designs S1 and S2 respectively for different $n$ and $d$. In each figure we plot values of $\mathbb{E}\theta(\mathbb{Z}_n)$ for different values of $r$; a step in $r$ increase gives the next curve up.

\begin{figure}[!h]
\centering
\begin{minipage}{.5\textwidth}
  \includegraphics[width=1\linewidth]{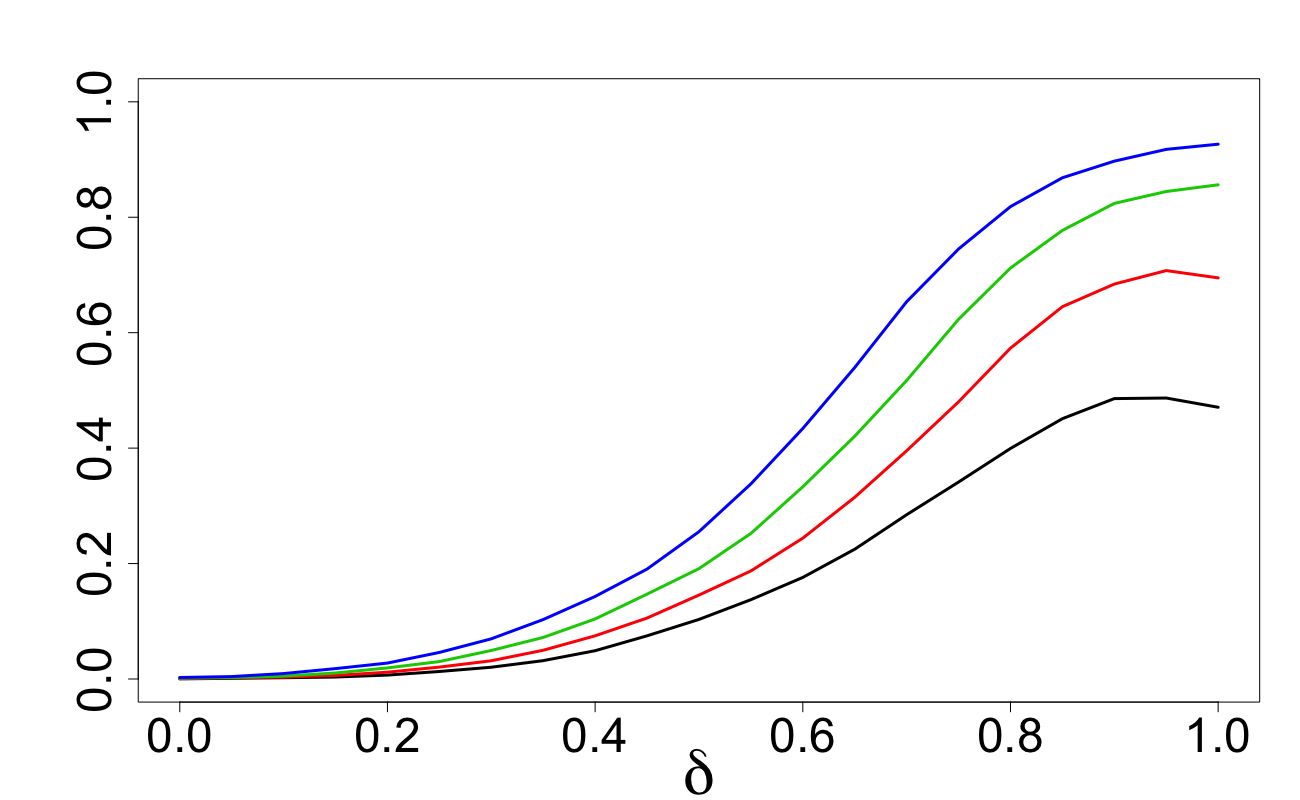}
  \caption{ ${ C}_d(\mathbb{Z}_n,r) $ for Design S1: $d=5$, $n=128$,\\ $r$ from $0.11$ to $0.17$ increasing by $0.02$. }
  \label{simplex_1_pic}
\end{minipage}%
\begin{minipage}{.5\textwidth}
  \centering
  \includegraphics[width=1\linewidth]{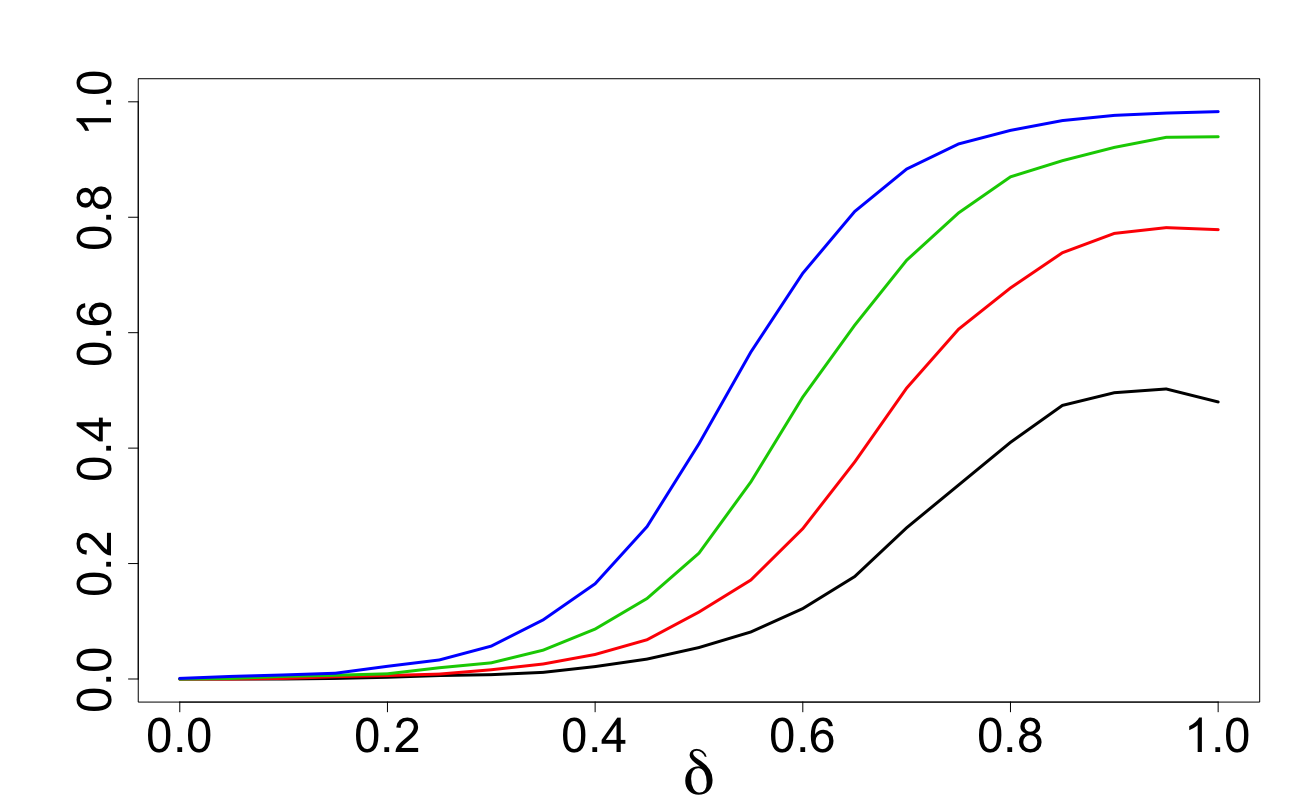}
  \caption{${ C}_d(\mathbb{Z}_n,r) $ for Design S1: $d=10$, $n=512$,\\ $r$ from $0.13$ to $0.19$ increasing by $0.02$. }
\end{minipage}
\end{figure}

\begin{figure}[!h]
\centering
\begin{minipage}{.5\textwidth}
  \includegraphics[width=1\linewidth]{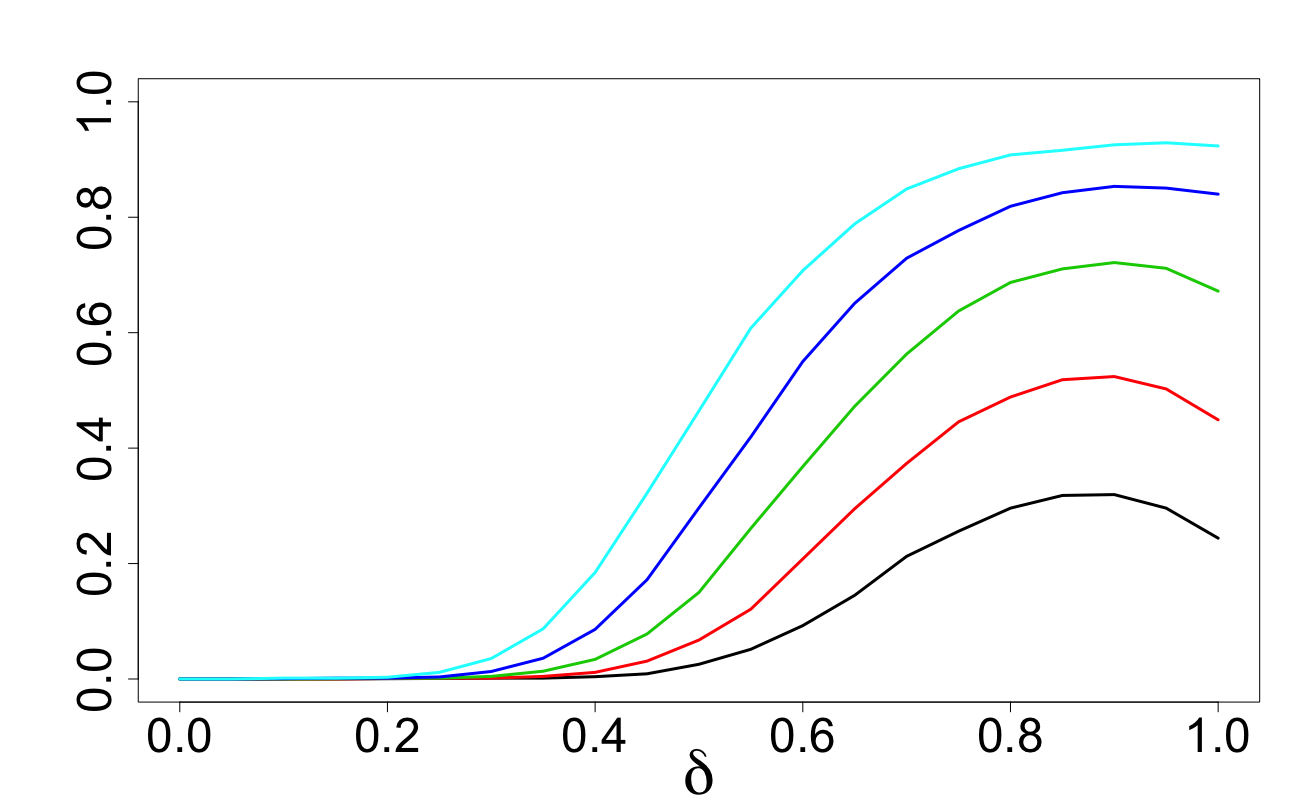}
  \caption{${ C}_d(\mathbb{Z}_n,r) $ for Design S1: $d=20$, \\$n=1024$, $r$ from $0.13$ to $0.17$ increasing by $0.01$. }
\end{minipage}%
\begin{minipage}{.5\textwidth}
  \centering
  \includegraphics[width=1\linewidth]{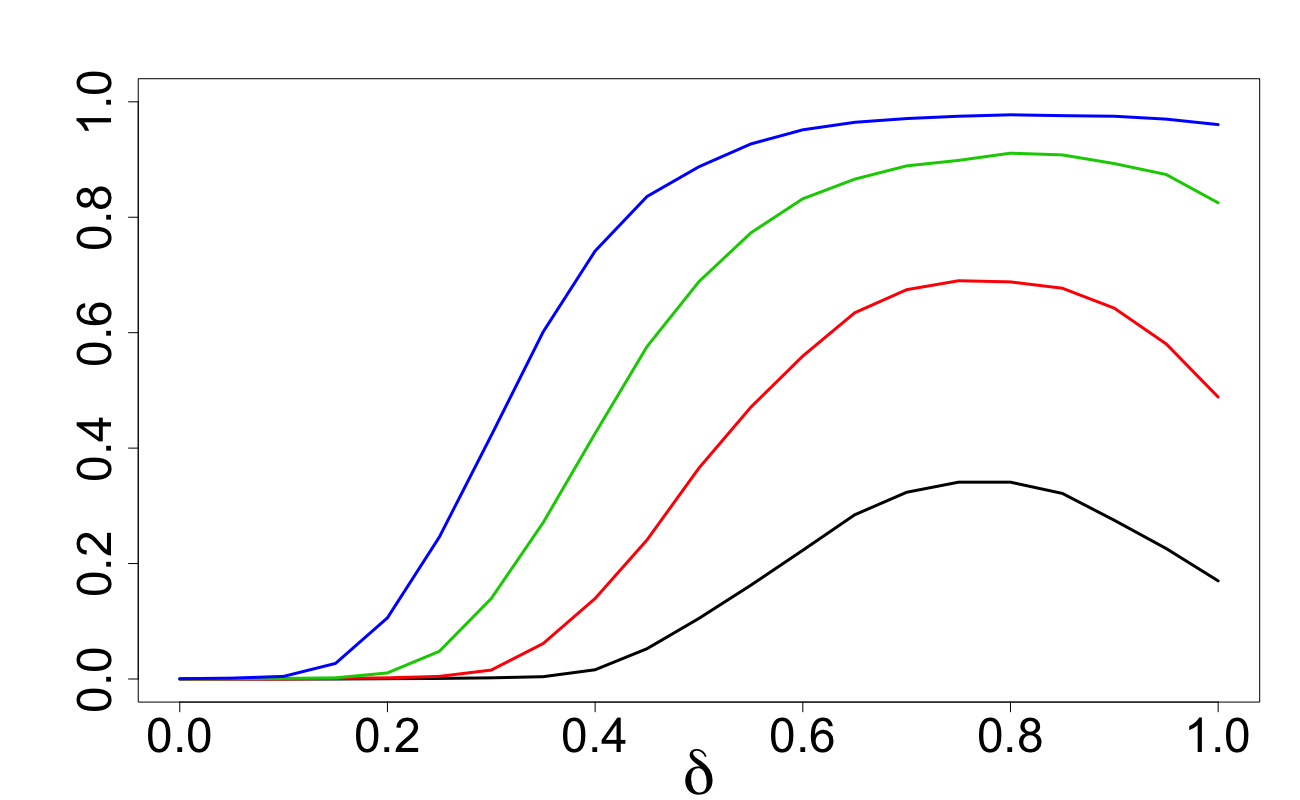}
  \caption{${ C}_d(\mathbb{Z}_n,r) $ for Design S1: $d=50$, $n=1024$,\\ $r$ from $0.12$ to $0.15$ increasing by $0.01$. }
    \label{simplex_1_pic_end}
\end{minipage}
\end{figure}

\begin{figure}[!h]
\centering
\begin{minipage}{.5\textwidth}
  \includegraphics[width=1\linewidth]{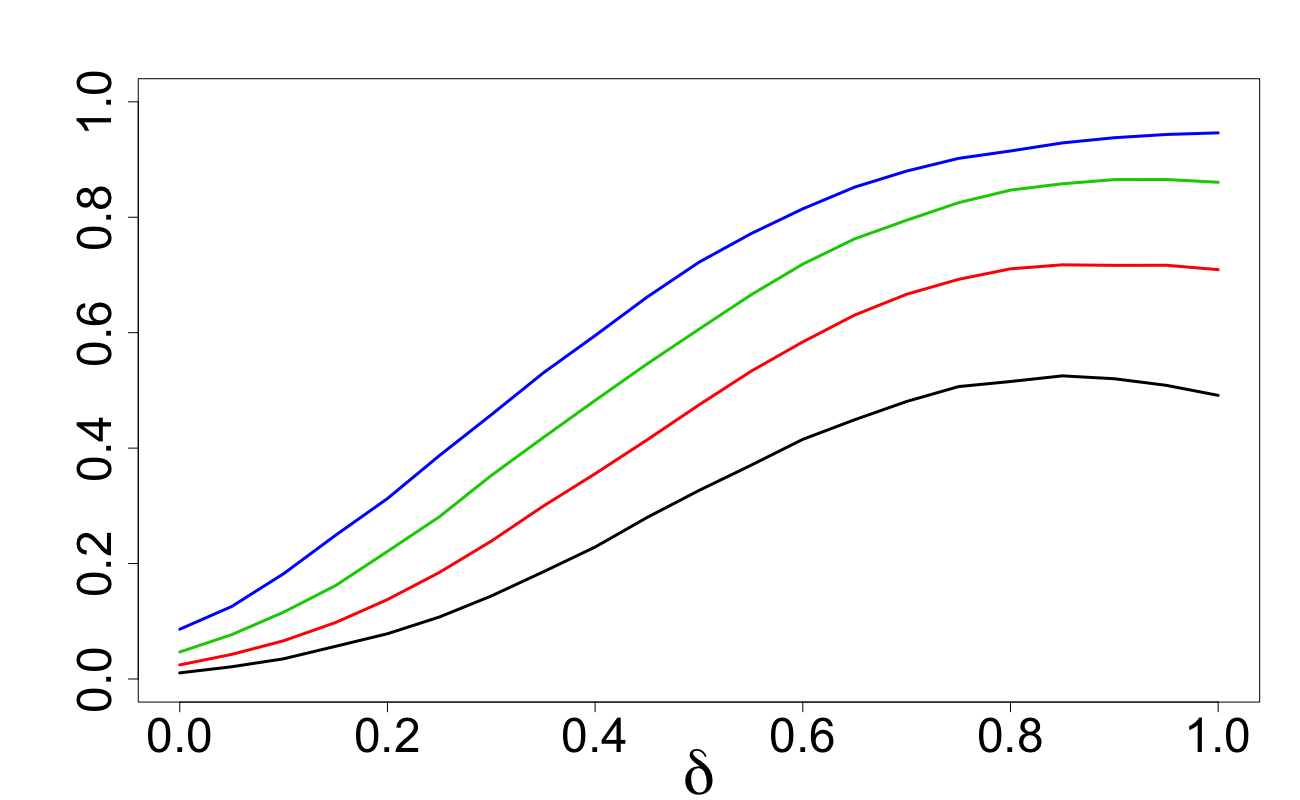}
  \caption{${ C}_d(\mathbb{Z}_n,r) $ for Design S2: $d=5$,\\ $n=128$, $r$ from $0.11$ to $0.17$ increasing by $0.02$. }
  \label{simplex_2_pic}
\end{minipage}%
\begin{minipage}{.5\textwidth}
  \centering
  \includegraphics[width=1\linewidth]{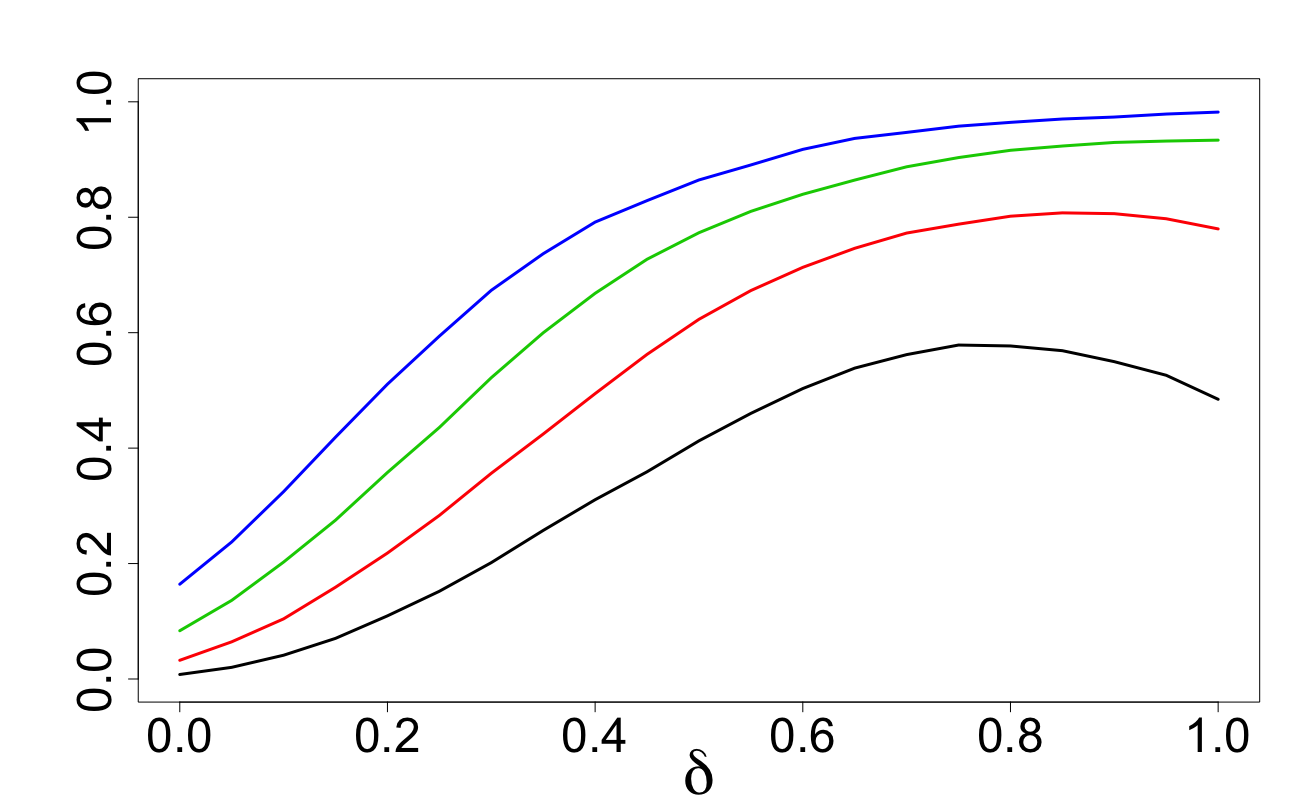}
  \caption{${ C}_d(\mathbb{Z}_n,r) $ for Design S2: $d=10$, $n=512$,\\ $r$ from $0.13$ to $0.19$ increasing by $0.02$.}
\end{minipage}
\end{figure}

\begin{figure}[!]
\centering
\begin{minipage}{.5\textwidth}
  \includegraphics[width=1\linewidth]{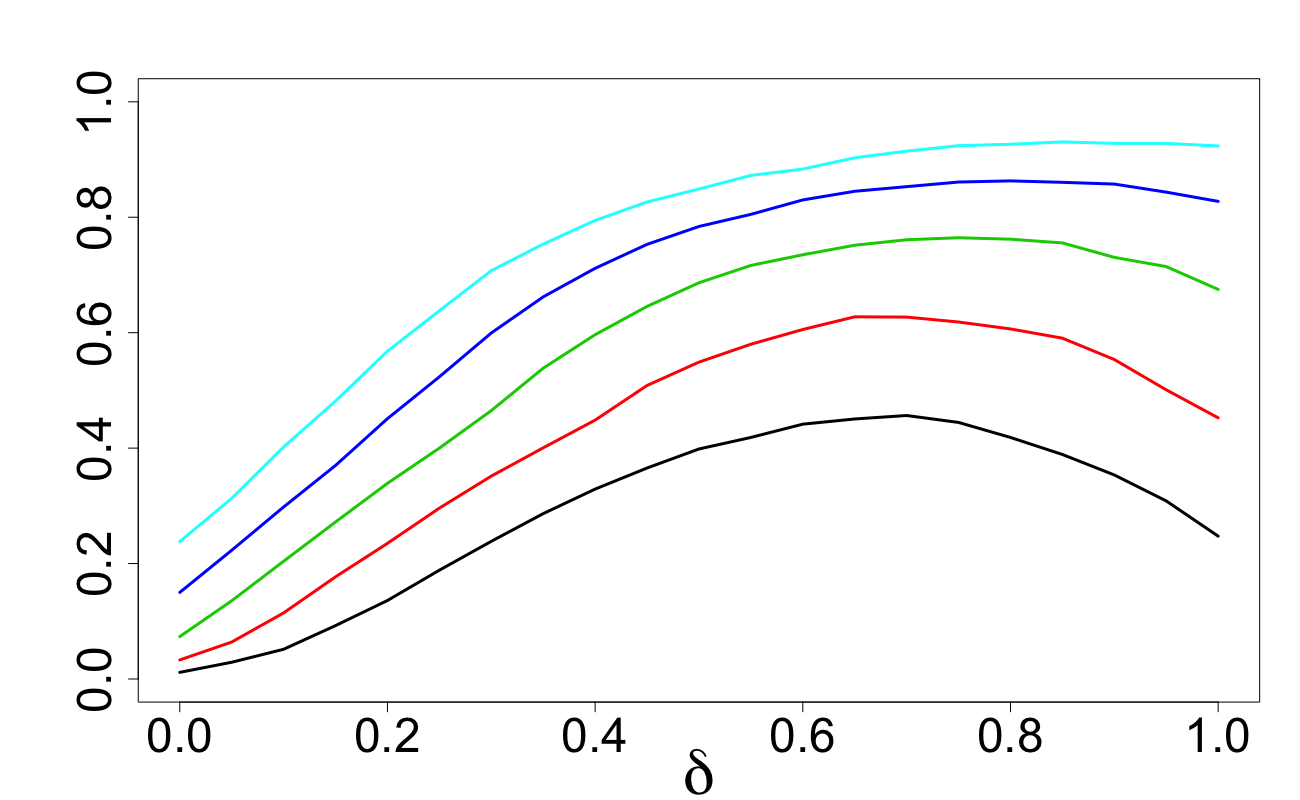}
  \caption{${ C}_d(\mathbb{Z}_n,r) $ for Design S2: $d=20$, \\$n=1024$, $r$ from $0.13$ to $0.17$ increasing by $0.01$.}
\end{minipage}%
\begin{minipage}{.5\textwidth}
  \centering
  \includegraphics[width=1\linewidth]{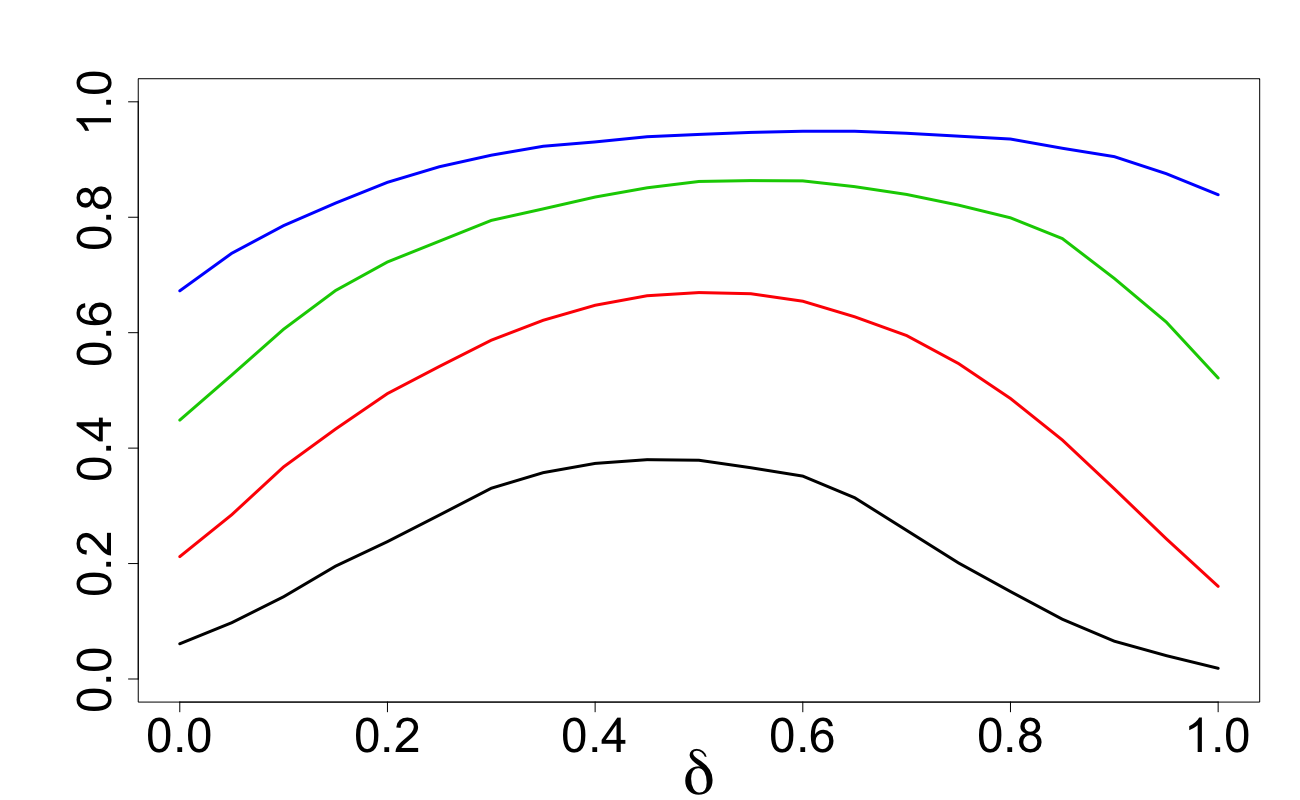}
  \caption{${ C}_d(\mathbb{Z}_n,r) $ for Design S2: $d=50$, $n=1024$,\\ $r$ from $0.11$ to $0.14$ increasing by $0.01$.}
  \label{simplex_2_pic_end}
\end{minipage}
\end{figure}

\begin{figure}[!h]
\centering
\begin{minipage}{.5\textwidth}
  \includegraphics[width=1\linewidth]{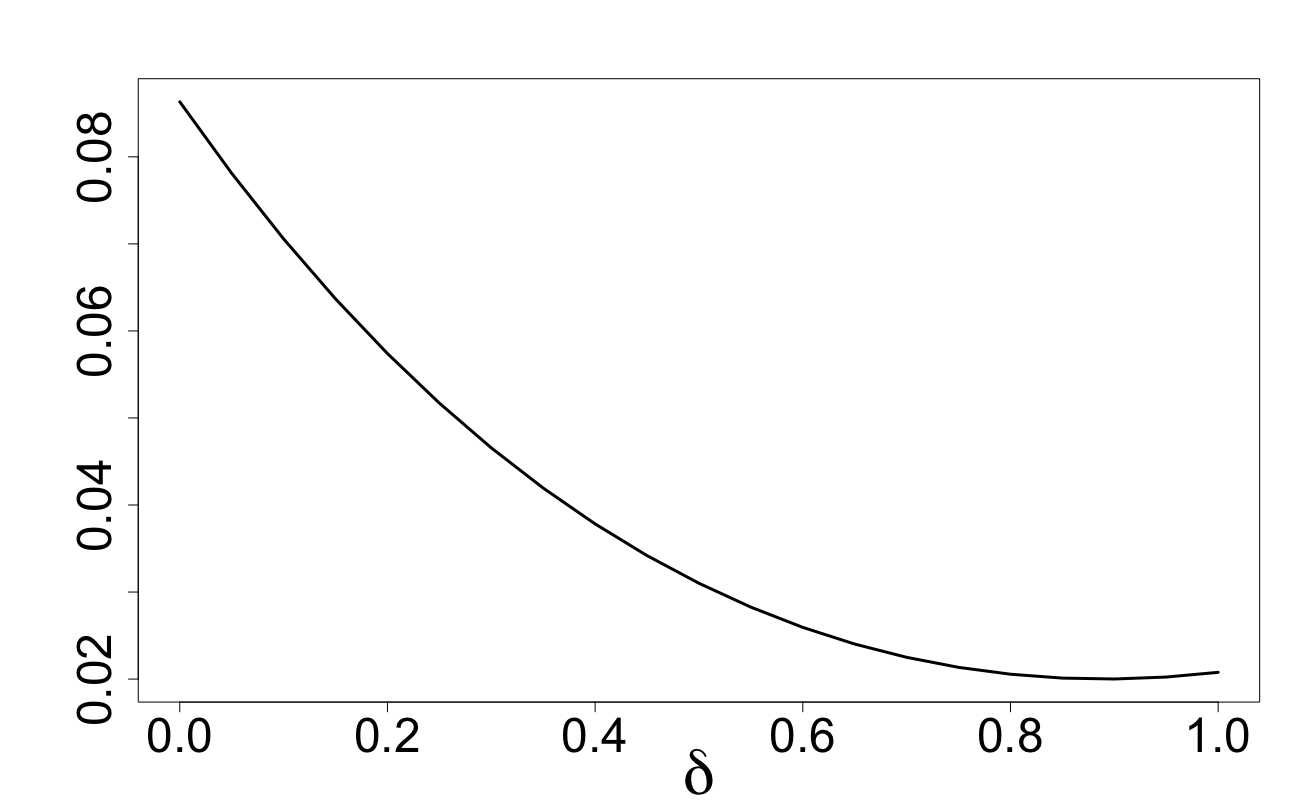}
  \caption{ $\mathbb{E}\theta(\mathbb{Z}_n)$ for Design S1: $d=20$, $n=1024$. }
  \label{quant_s1}
\end{minipage}%
\begin{minipage}{.5\textwidth}
  \centering
  \includegraphics[width=1\linewidth]{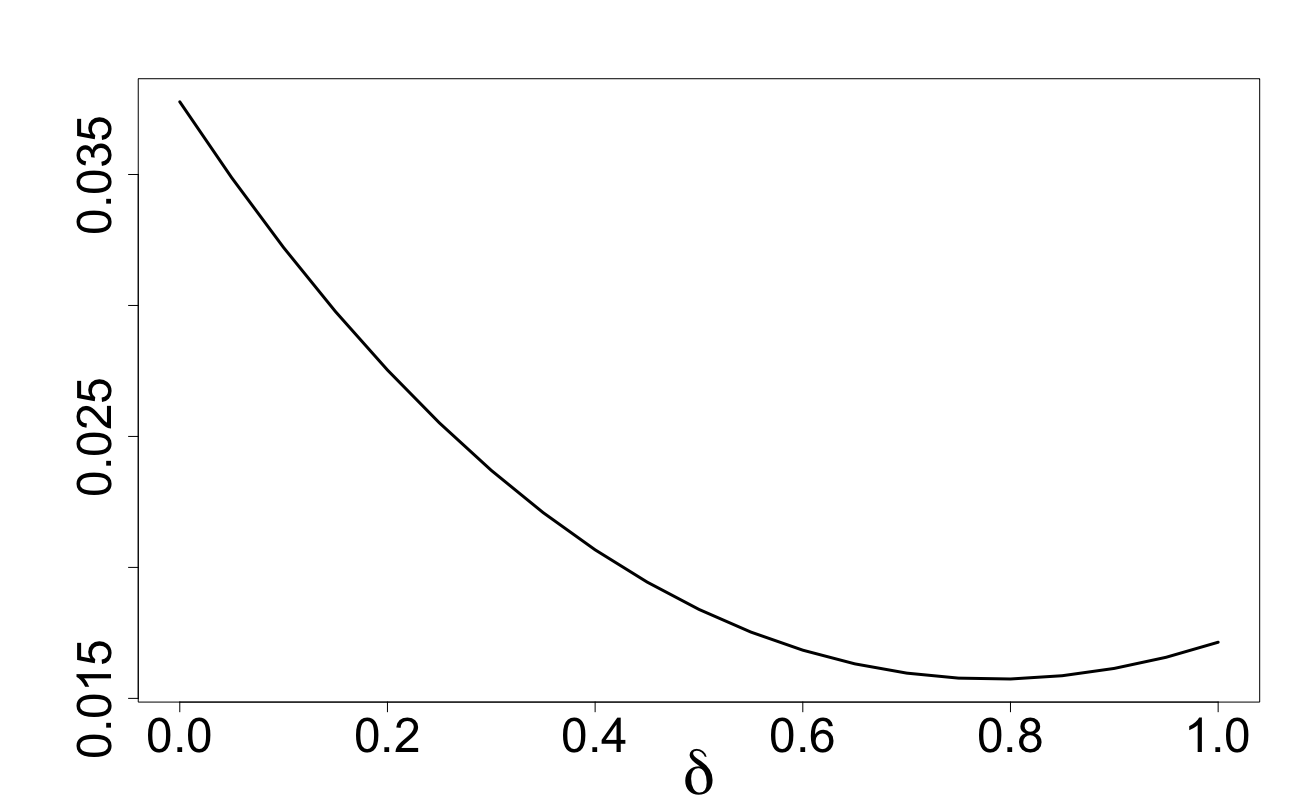}
  \caption{ $\mathbb{E}\theta(\mathbb{Z}_n)$ for Design S1: $d=50$, $n=1024$. }
    \label{quant_s1_end}
\end{minipage}
\end{figure}

\begin{figure}[!h]
\centering
\begin{minipage}{.5\textwidth}
  \includegraphics[width=1\linewidth]{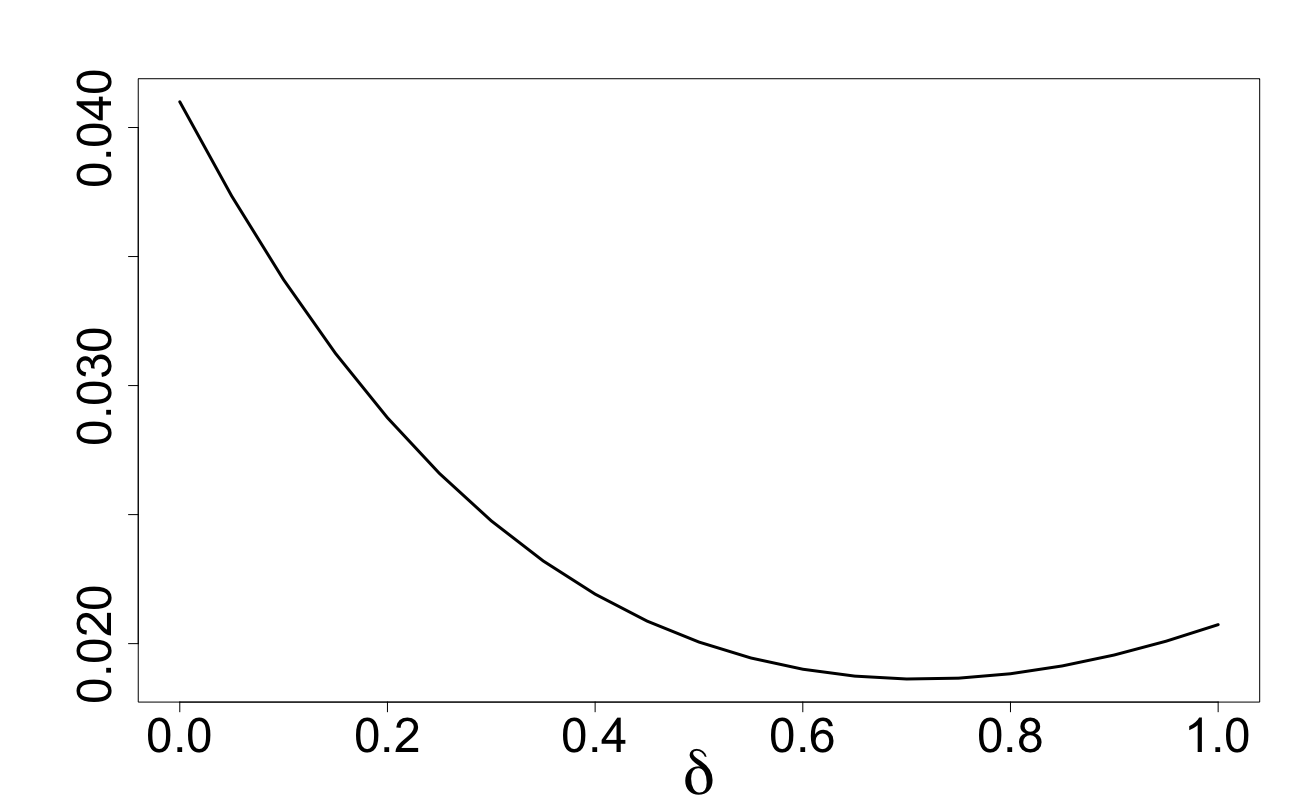}
  \caption{ $\mathbb{E}\theta(\mathbb{Z}_n)$ for Design S2: $d=20$, $n=1024$. }
    \label{quant_s2}
\end{minipage}%
\begin{minipage}{.5\textwidth}
  \centering
  \includegraphics[width=1\linewidth]{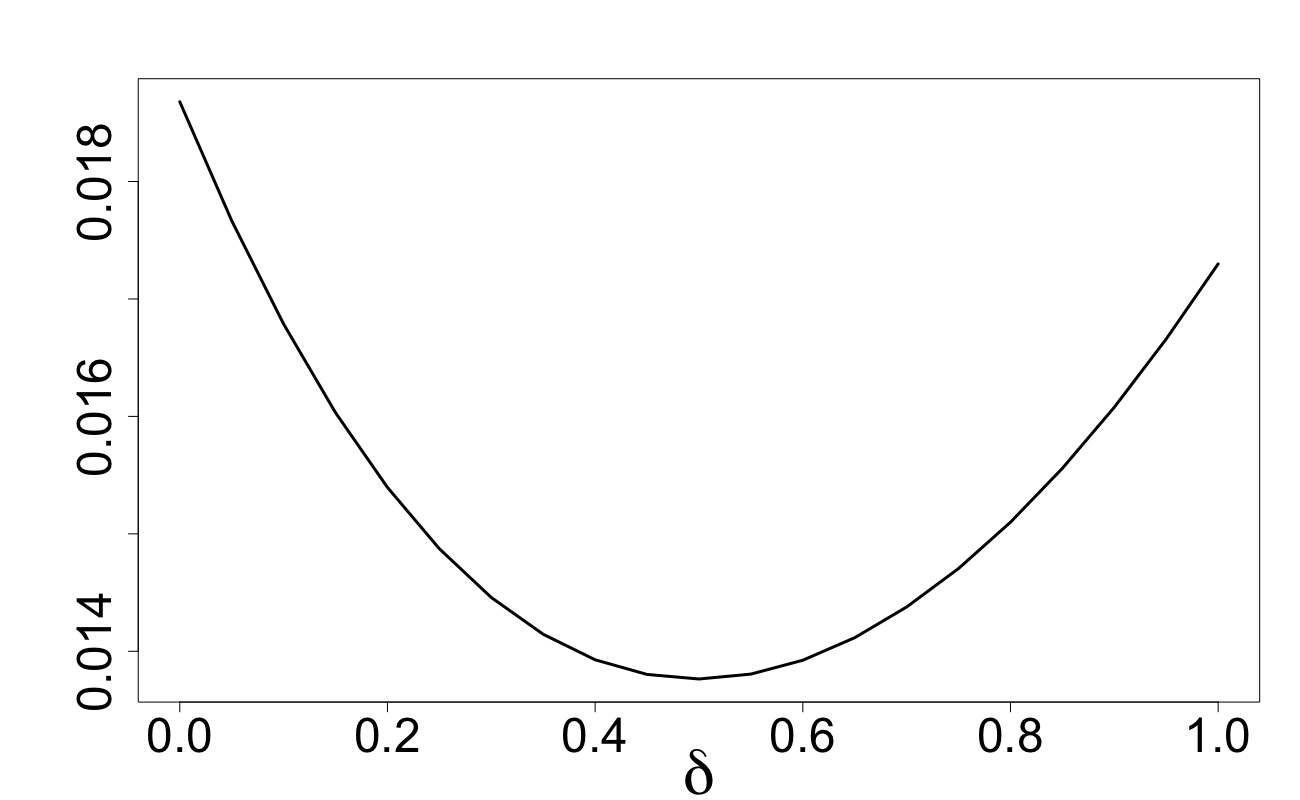}
  \caption{ $\mathbb{E}\theta(\mathbb{Z}_n)$ for Design S2: $d=50$, $n=1024$. }
    \label{quant_s2_end}
\end{minipage}
\end{figure}
\clearpage

From the above figures, we arrive at the following conclusions:
\begin{itemize}
\item The $\delta$-effect for the simplex is much less prominent than for the cube.
\item Between Designs S1 and S2,  the $\delta$-effect is more apparent for Design S2; for example, compare Figure~\ref{quant_s1_end} with  Figure~\ref{quant_s2_end}.
\end{itemize}

\section{Appendix: An auxiliary lemma}
\label{sec:appB}

{\bf Lemma 1.} \label{lem:1}
{\it Let $\delta>0$,  $u  \in \mathbb{R} $ and $\eta_{u,\delta}$ be a r.v. $\eta_{u,\delta} = (\xi-u)^2$, where
r.v.~$\xi \in [-\delta,\delta]$ has  Beta$_\delta(\alpha,\alpha)$ distribution with density
\be
\label{eq:beta1}
p_{\alpha,\delta}(t)= \frac{(2\delta)^{1-2\alpha}}{\mbox{Beta$(\alpha,\alpha)$}} [\delta^2-t^2]^{\alpha-1}\, , \;\;-\delta<t<\delta\, ,\alpha>0;
\ee
Beta$(\cdot,\cdot)$ is the Beta-function. The r.v. $\eta_{u,\delta} $ is concentrated on the interval $ [(\max(0, \delta-|u|) )^2,(\delta+|u| )^2]$.
Its first three  central moments  are:
\bea
\label{eq:inters1c}
\mu_{u,\delta}^{(1)} &=& E\eta_{u,\delta} =u^2+ \frac {{{\delta}}^{2}}{2\,{\alpha}+1} \, ,\\
\label{eq:inters1c2}
\mu_{u,\delta}^{(2)}&=&{\rm var} (\eta_{u,\delta}) =
{\frac {4\delta^{2}}{2\,{\alpha}+1}}
 \left[u^2+
 {\frac {{{\delta}}^{2}{\alpha}}{ \left( 2\,{\alpha}+1
 \right)  \left( 2\,{\alpha}+3 \right) }}
  \right]  \, ,\\
\label{eq:inters1c3}
\mu_{u,\delta}^{(3)}&=&E \left[\eta_{u,\delta} - E\eta_{u,\delta}\right]^3 =
{\frac {48{\alpha}\,{{\delta}}^{4}}{ \left( 2\,{\alpha}+1
 \right) ^{2} \left( 2\,{\alpha}+3 \right) }} \left[  u^2+ {\frac {{{\delta}}^{2} \left( 2\,{\alpha} -1\right) }{3
 \left( 2\,{\alpha}+5 \right)  \left( 2\,{\alpha}+1 \right) }}
 \right]
\, .
\eea

}

In the limiting case $\alpha=0$, where the r.v.~$\xi $ is concentrated at two points $\pm \delta$ with equal weights, we obtain:
$
\mu_{u,\delta}^{(1)} = E\eta_{u,\delta} =u^2+ \delta^{2} \,
$
and
{
\be \label{eq:mom_0}
\mu_{u,\delta}^{(2k)} = \left[2\delta u\right]^{2k},\;\;\; \mu_{u,\delta}^{(2k+1)}=0,\;\;{\rm for}\;\; k=1,2, \ldots
\ee
}
\vspace{-1cm}

\bibliographystyle{plain}
\bibliography{large_dimension}

\end{document}